%% file: StochSimMC.tex
\documentclass[12pt]{article}

\input{headlines}

\def\blue{\textcolor[rgb]{0.00,0.00,1.00}}

\usepackage{algpseudocode}

\title{\large\sffamily\bfseries Lecture Notes\\
\LARGE{Stochastic Simulation and Monte Carlo Method}\\
}
\author{\Large\sffamily Davoud Mirzaei\\ Uppsala University}
\date{\sffamily November 1, 2024 ~~ (2nd Edition)}
\begin{document}

\setlength{\abovedisplayskip}{3pt}
\setlength{\belowdisplayskip}{3pt}
\setlength{\abovedisplayshortskip}{0pt}
\setlength{\belowdisplayshortskip}{0pt}
\maketitle

%Custom colors for different environments
\definecolor{contcol1}{HTML}{72E094}
\definecolor{contcol2}{HTML}{24E2D6}
\definecolor{convcol1}{HTML}{C0392B}
\definecolor{convcol2}{HTML}{8E44AD}

\begin{tcolorbox}[title=Contents, fonttitle=\huge\sffamily\bfseries\selectfont,interior style={left color=contcol1!10!white,right color=contcol2!10!white},frame style={left color=contcol1!30!white,right color=contcol2!30!white},coltitle=black,top=2mm,bottom=2mm,left=2mm,right=2mm,drop fuzzy shadow,enhanced,breakable]
\makeatletter
\@starttoc{toc}
\makeatother
\end{tcolorbox}

\vspace*{10mm}

\thispagestyle{empty}
\newpage
\pagenumbering{arabic}

\input{part_DeterStoch}

\input{part_MC1}

\input{part_generation}

\input{part_MC2}

\input{part_StochProc}

\input{part_StochProcGen}

\input{part_SSA}

\input{part_MCMC}

\input{workouts_4}

\appendix

\input{Appendix}
\vsp

\end{document}

%% file: headlines.tex
\usepackage[english]{babel}
\usepackage{graphicx}
\usepackage{framed}
\usepackage[normalem]{ulem}
\usepackage{indentfirst}
\usepackage{amsmath,amsthm,amssymb,amsfonts,dsfont,mathdots}
\usepackage[italicdiff]{physics}
\usepackage[T1]{fontenc}
\usepackage{mathtools,extarrows}

\usepackage{setspace}
\usepackage{wrapfig}
\usepackage{lmodern,mathrsfs}
\usepackage[inline,shortlabels]{enumitem}
\setlist{topsep=2pt,itemsep=2pt,parsep=3pt,partopsep=2pt}
\usepackage[dvipsnames]{xcolor}
\usepackage[utf8]{inputenc}
\usepackage[a4paper, top=0.8in,bottom=0.8in, left=0.8in, right=0.8in, footskip=0.3in, includefoot]{geometry}
\usepackage[most]{tcolorbox}
\usepackage{tikz,tikz-3dplot,tikz-cd,tkz-tab,tkz-euclide,pgf,pgfplots}
\pgfplotsset{compat=newest}
\usepackage{multicol}
\usepackage[bottom,multiple]{footmisc} %ensures footnotes are at the bottom of the page, and separates footnotes by a comma if they are adjacent

\usepackage[backend=bibtex,style=numeric]{biblatex}
\usepackage{graphicx}
 \graphicspath{{figures/}} %Setting the graphicspath

\numberwithin{equation}{section}

\usepackage{framed,color}
\definecolor{shadecolor}{rgb}{0.9412,0.9725, 1}
\usepackage{caption}
\usepackage{tcolorbox}

\usepackage{multirow}
\usepackage{epstopdf}
\usepackage{listings}
\usepackage[]{mcode}
\usepackage{algorithm}
\usepackage{latexsym}
\usepackage{showidx}
\usepackage{latexsym}
\usepackage{amssymb}

\newcommand{\ignore}[1]{}
\usepackage{todonotes}

\usepackage{tikz}

\newcommand{\hide}[1]{} %To hide large blocks of code without using % symbols

\renewcommand{\qed}{\hfill\blacksquare}

\newcommand{\E}{\mathbb{E}}

\newcommand{\N}{\mathbb{N}}

\newcommand{\R}{\mathbb{R}}

\newcommand{\x}{\boldsymbol{x}}
\newcommand{\y}{\boldsymbol{y}}

\newcommand{\vv}{\boldsymbol{v}}

\newcommand{\Ss}{\mathcal{S}}

\newcommand{\bpi}{\boldsymbol{\pi}}
\newcommand{\e}{\boldsymbol{e}}

\def\ee{\mathrm{e}}

\newtheoremstyle{mystyle}{}{}{}{}{\sffamily\bfseries}{.}{ }{}
\newtheoremstyle{cstyle}{}{}{}{}{\sffamily\bfseries}{.}{ }{\thmnote{#3}}
\makeatletter
\renewenvironment{proof}[1][\proofname] {\par\pushQED{\qed}{\normalfont\sffamily\bfseries\topsep6\p@\@plus6\p@\relax #1\@addpunct{.} }}{\popQED\endtrivlist\@endpefalse}
\makeatother
\theoremstyle{mystyle}{\newtheorem{definition}{Definition}[section]}
\theoremstyle{mystyle}{\newtheorem{theorem}[definition]{Theorem}}
\theoremstyle{mystyle}{}
\theoremstyle{mystyle}{}
\theoremstyle{mystyle}{}
\theoremstyle{mystyle}{}
\theoremstyle{mystyle}{}
\theoremstyle{mystyle}{\newtheorem{workout}[definition]{Exercise}}
\theoremstyle{mystyle}{}
\theoremstyle{mystyle}{}

\usepackage{amsmath}
\newtheorem{remark}{Remark}[section]
\newtheorem{example}{Example}[section]

\theoremstyle{cstyle}{}

\newtheoremstyle{warn}{}{}{}{}{\normalfont}{}{ }{}
\theoremstyle{warn}

\newcommand{\warningsign}[1]{\tikz[scale=#1,every node/.style={transform shape}]{\draw[-,line width={#1*0.8mm},red,fill=yellow,rounded corners={#1*2.5mm}] (0,0)--(1,{-sqrt(3)})--(-1,{-sqrt(3)})--cycle;
\node at (0,-1) {\fontsize{48}{60}\selectfont\bfseries!};}}

\tcolorboxenvironment{definition}{boxrule=0pt,boxsep=0pt,colback={red!10},left=8pt,right=8pt,enhanced jigsaw, borderline west={2pt}{0pt}{red},sharp corners,before skip=10pt,after skip=10pt,breakable}
\tcolorboxenvironment{proposition}{boxrule=0pt,boxsep=0pt,colback={Orange!10},left=8pt,right=8pt,enhanced jigsaw, borderline west={2pt}{0pt}{Orange},sharp corners,before skip=10pt,after skip=10pt,breakable}
\tcolorboxenvironment{theorem}{boxrule=0pt,boxsep=0pt,colback={blue!10},left=8pt,right=8pt,enhanced jigsaw, borderline west={2pt}{0pt}{blue},sharp corners,before skip=10pt,after skip=10pt,breakable}
\tcolorboxenvironment{lemma}{boxrule=0pt,boxsep=0pt,colback={blue!10},left=8pt,right=8pt,enhanced jigsaw, borderline west={2pt}{0pt}{blue},sharp corners,before skip=10pt,after skip=10pt,breakable}
\tcolorboxenvironment{corollary}{boxrule=0pt,boxsep=0pt,colback={violet!10},left=8pt,right=8pt,enhanced jigsaw, borderline west={2pt}{0pt}{violet},sharp corners,before skip=10pt,after skip=10pt,breakable}
\tcolorboxenvironment{proof}{boxrule=0pt,boxsep=0pt,blanker,borderline west={2pt}{0pt}{CadetBlue!80!white},left=8pt,right=8pt,sharp corners,before skip=10pt,after skip=10pt,breakable}
\tcolorboxenvironment{remark}{boxrule=0pt,boxsep=0pt,blanker,borderline west={2pt}{0pt}{Green},left=8pt,right=8pt,before skip=10pt,after skip=10pt,breakable}
\tcolorboxenvironment{remarks}{boxrule=0pt,boxsep=0pt,blanker,borderline west={2pt}{0pt}{Green},left=8pt,right=8pt,before skip=10pt,after skip=10pt,breakable}
\tcolorboxenvironment{example}{boxrule=0pt,boxsep=0pt,blanker,borderline west={2pt}{0pt}{Black},left=8pt,right=8pt,sharp corners,before skip=10pt,after skip=10pt,breakable}
\tcolorboxenvironment{examples}{boxrule=0pt,boxsep=0pt,blanker,borderline west={2pt}{0pt}{Black},left=8pt,right=8pt,sharp corners,before skip=10pt,after skip=10pt,breakable}
\tcolorboxenvironment{cthm}{boxrule=0pt,boxsep=0pt,colback={gray!10},left=8pt,right=8pt,enhanced jigsaw, borderline west={2pt}{0pt}{gray},sharp corners,before skip=10pt,after skip=10pt,breakable}
\tcolorboxenvironment{workout}{boxrule=0pt,boxsep=0pt,colback={Cyan!0},left=8pt,right=8pt,enhanced jigsaw, borderline west={2pt}{0pt}{Cyan},sharp corners,before skip=10pt,after skip=10pt,breakable}
\tcolorboxenvironment{miniproject}{boxrule=0pt,boxsep=0pt,colback={gray!10},left=8pt,right=8pt,enhanced jigsaw, borderline west={2pt}{0pt}{gray},sharp corners,before skip=10pt,after skip=10pt,breakable}
\tcolorboxenvironment{labexercise}{boxrule=0pt,boxsep=0pt,colback={green!10},left=8pt,right=8pt,enhanced jigsaw, borderline west={2pt}{0pt}{green},sharp corners,before skip=10pt,after skip=10pt,breakable}

\newenvironment{talign*}{\let\displaystyle\textstyle\csname align*\endcsname}{\endalign}

\usepackage[explicit]{titlesec}
\titleformat{\section}{\fontsize{18}{30}\sffamily\bfseries}{\thesection}{18pt}{#1}
\titleformat{\subsection}{\fontsize{15}{18}\sffamily\bfseries}{\thesubsection}{15pt}{#1}
\titleformat{\subsubsection}{\fontsize{10}{12}\sffamily\large\bfseries}{\thesubsubsection}{8pt}{#1}

\titlespacing*{\section}{0pt}{5pt}{5pt}
\titlespacing*{\subsection}{0pt}{5pt}{5pt}
\titlespacing*{\subsubsection}{0pt}{5pt}{5pt}

\newcommand{\ds}{\displaystyle}

\DeclareMathAlphabet\mathbfcal{OMS}{cmsy}{b}{n}
\newcommand{\vsp}{\vspace{0.3cm}}
\newcommand{\pr}{\mathds{P}}
\newcommand{\Var}{\mathrm{Var}}
\newcommand{\Cov}{\mathrm{Cov}}
\newcommand{\td}{\mathrm{d}}

\newcommand{\Nor}{\mathcal{N}}
\newcommand{\Exp}{\mathcal{E}xp}
\newcommand{\Poi}{\mathcal{P}oi}
\newcommand{\Bin}{\mathcal{B}in}
\newcommand{\Ber}{\mathcal{B}er}
\newcommand{\Geo}{\mathcal{G}eo}
\newcommand{\Uni}{\mathcal{U}}
\newcommand{\Bet}{\mathcal{B}eta}
\newcommand{\Gam}{\mathcal{G}am}

%% file: part_DeterStoch.tex
These lecture notes are intended to cover some introductory topics in stochastic simulation for scientific computing courses offered by the IT department at Uppsala University, as taught by the author. Basic concepts in probability theory are provided in the Appendix \ref{sect:appendixA}, which you may review before starting the upcoming sections or refer to as needed throughout the text.
Some parts of our presentation here follow  \cite{DeGroot-Schervish:2007}, \cite{Ross:2002} and \cite{Rubinnstein-Kroese:2017}.

\section{Deterministic vs. Stochastic}

In the field of scientific computing, modeling and simulation are fundamental tools used to understand natural phenomena. Two main approaches are commonly employed for modeling and simulation: {\bf deterministic} and {\bf stochastic}. 
While both approaches aim to predict the behavior of systems, they differ  in how they handle uncertainty and randomness.
\vsp

\subsection{Deterministic models and methods}

A deterministic model is one in which the behavior of the system is entirely predictable and reproducible. Given a specific set of initial conditions and parameters, the outcome will always be the same, i.e. the future behavior can be predicted with complete certainty. Such models are often governed by mathematical equations (such as differential equations or algebraic relations) that describe precise relationships between different quantities. 
For example, consider Newton's laws of motion, which describe how objects move in response to forces. These laws are deterministic because, given the initial conditions (such as the position, velocity, and force acting on an object) the future motion of the object can be calculated exactly. All mathematical models expressed as ordinary and partial differential equations (ODEs and PDEs) are examples of deterministic models. 

Deterministic methods are employed to solve deterministic models\footnote{It is sometimes possible to use a stochastic method to solve a deterministic model. For example, a Monte Carlo method can be used to approximate the volume of an object, which is equivalent to solving a multi-dimensional integral.}. Methods such as Euler's method for solving ODEs, the trapezoidal rule for integration, and the Newton-Raphson method for solving non-linear equations are basic examples of deterministic methods. When it comes to more complicated models such as PDEs, the finite difference method (FDM), the finite element method (FEM), and the finite volume method (FVM) are examples of commonly used deterministic methods for solving such deterministic models.  

As a simple example, consider the decay of a radioactive material. This process can be described by a first-order differential equation, which relates the rate of decay to the amount of radioactive material present at a given time. The governing equation has the form
\[
\frac{dy}{dt} = -\lambda y(t)
\]
where \( y(t) \) represents the amount of radioactive material at time \( t \), and \( \lambda \) is the decay constant. The initial amount of material is given by \( y(0) = y_0 \). This model is deterministic because, given the initial amount of radioactive material \( y_0 \) and the decay constant \( \lambda \), the future amount of material can be calculated precisely at any time. Indeed, the solution to this equation is
\[
y(t) = y_0 e^{-\lambda t}.
\]
This solution describes the exponential decay of the material over time. No matter how many times we perform this calculation, we will always obtain the same result for a given set of input parameters and initial conditions.

However, solving differential equations analytically is not always possible or practical. In such cases, numerical methods are used to approximate the solutions. One common solver is the Runge-Kutta method. An adaptive version of this method has been implemented in Python library  $\texttt{scipy.integrate.solve\_ivp}$. Here, we call this library and compute a numerical solution to the radioactive decay equation for $\lambda = 0.5$ and $y_0 = 10^3$.
\begin{shaded}
\vspace*{-0.3cm}
\begin{verbatim}
import numpy as np
import matplotlib.pyplot as plt
from scipy.integrate import solve_ivp
lam, y0, FinalTime  = 0.5, 1000, 10   # rate, initial value, final time
def ODEfun(t,y):
    yprime = -lam*y
    return yprime
teval = np.linspace(0, FinalTime, 500)
sol = solve_ivp(ODEfun, [0,FinalTime], y0, t_eval = teval)
plt.figure(figsize = (6, 4))
plt.plot(sol.t,sol.y[0],linestyle = 'solid', color='blue')
plt.xlabel('$t$ (seconds)'); 
plt.ylabel('Amount of radioactive material, $y(t)$')
plt.title('Deterministic solution, radioactive decay: $y_0=1000,\lambda=0.5$')
\end{verbatim}
\vspace*{-0.3cm}
\end{shaded}

The plot is given on the left-hand side of Figure \ref{fig:deter-stoch}. 

\begin{center}
\includegraphics[scale=.5]{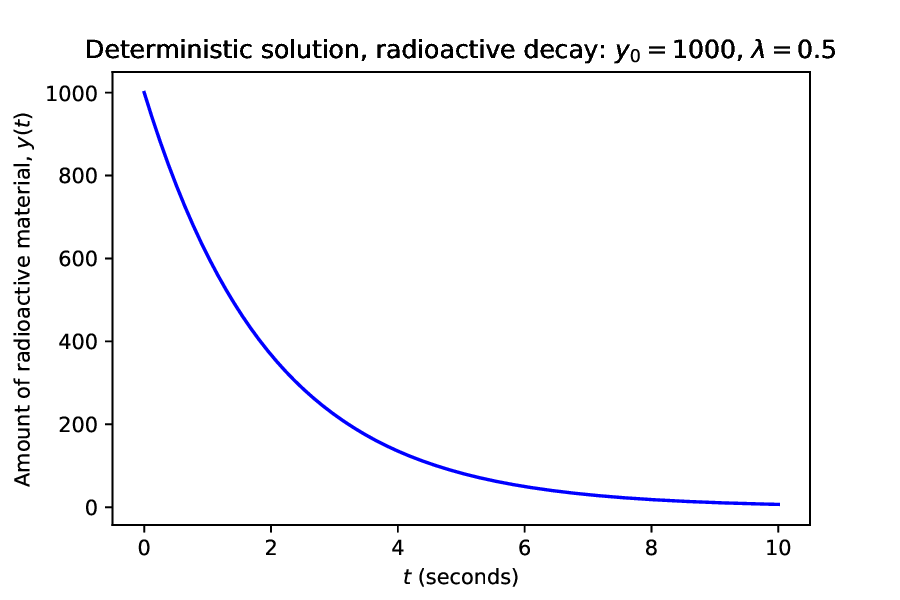}
\includegraphics[scale=.5]{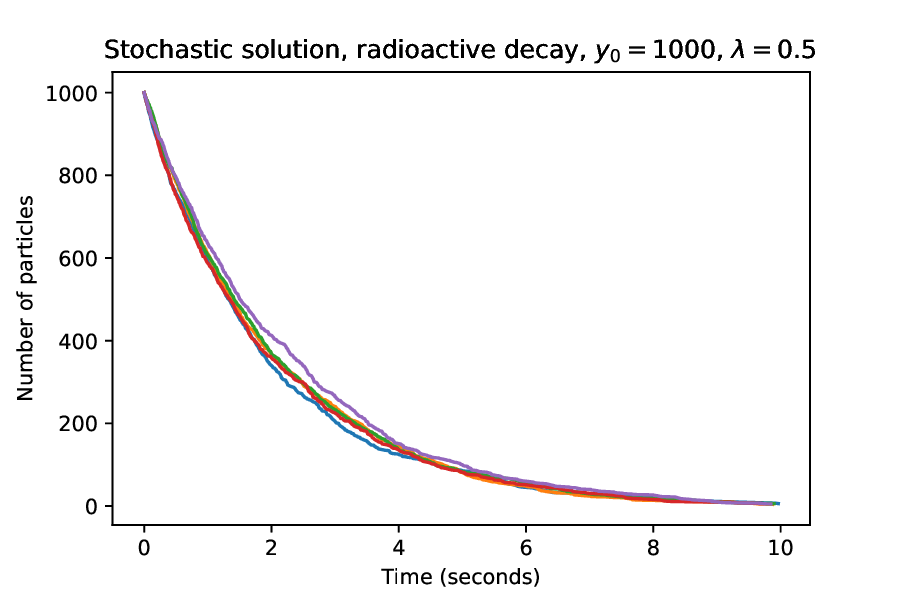}
\captionof{figure}{Numerical solution using the deterministic ODE-solver \texttt{RK45}
from Python library $\texttt{scipy.integrate.solve\_ivp}$ (left), and using the stochastic solver SSA (right)
}\label{fig:deter-stoch}
\end{center}

As we observe, similar to the analytical solution, the numerical scheme yields a deterministic outcome for the given input values. Randomness does not play a role here, and all results are fixed. For instance, the amount of material at time \( t = 4 \) is a precise value, \( y_4 = 135.3353 \) (rounded to four decimal places).
\vsp

\subsection{Stochastic models and methods}

Deterministic models are not always an accurate reflection of reality. For instance, in the above example the exact time at which an individual atom will decay is random. Similarly, in population dynamics, deterministic models cannot account for the random events such as environmental changes or disease outbreaks. 
Stochastic models incorporate randomness and uncertainty to have a better reflection of the behavior of real-world systems.
In a stochastic model, the same set of initial conditions can lead to different outcomes. Instead of predicting a single and precise result, these models describe a range of possible outcomes, each associated with a certain probability. These models are particularly useful in areas such as biology, finance, and quantum mechanics where small-scale fluctuations or random events play a significant role in determining the overall behavior of the system.

Let us revisit the example of radioactive decay. In the deterministic model, we used a differential equation to predict the amount of radioactive material at a given time. However, this model assumes that the decay process occurs continuously and deterministically, which is not the case in reality. In fact, radioactive decay is a random process, and the time at which each atom decays is uncertain. 
In a stochastic model, we describe the system as a series of random events. Specifically, we model the decay process as a reaction:
\[
y \xrightarrow{\; \lambda\; } z
\]
where \( y \) represents the number of radioactive molecules, \( z \) represents the decay products, and \( \lambda \) is the {\em propensity} of decay. The difference between this model and the deterministic model is that we no longer assume a continuous, predictable decay process. Instead, we model the decay as a random event that occurs with a certain probability, and the waiting time for the next decay is also a random variable. 
This is one of the key features of stochastic models that they contain random variables and {\em probability distributions} to describe the likelihood of different outcomes. 

To simulate this stochastic model, we can use a method known as {\em Gillespie’s Algorithm}, or the {\em Stochastic Simulation Algorithm (SSA)}. See section \ref{sect:ssa} for details. Each time we run the simulation, we will obtain a different result which reflects the randomness of the decay process. 

See the right panel in Figure \ref{fig:deter-stoch} for five simulation outcomes, with the respective values of \(y\) at $t=4$ being \(122\), \(228\), \(134\), \(135\), and \(139\) for each simulation. 
In this example, if we run the simulation many times and compute the average behavior of the system, we will find that the result closely approximates the deterministic solution. However, this is not always the case. When we model a phenomenon using both a stochastic model and a deterministic model (as we did with radioactive decay), the mean of stochastic solutions does not always approximate the deterministic solution. In some cases, while the deterministic solution converges to an equilibrium state, random noise and fluctuations may force all stochastic solutions to deviate from the equilibrium. 
\vsp 

%\begin{center}
%\includegraphics[scale=.5]{stoch_ex0}
%\captionof{figure}{Numerical solution using the stochastic solver SSA}%\label{fig:ssa_sol}
%\end{center}

\subsection{Which one?}

The choice between a deterministic and a stochastic model depends on the nature of the system being studied and the goals of the simulation. Deterministic models are ideal for systems that behave predictably, where small-scale fluctuations have little impact on the overall behavior of the system. These models are computationally efficient and provide precise and repeatable results. They are well-suited for systems where accuracy is important.
However, deterministic models can be limited in their ability to describe some real-world systems that are subject to randomness and uncertainty. In such cases, stochastic models provide a more realistic description of the system behavior. By incorporating randomness, stochastic models can capture the variability and uncertainty. 

One important consideration when choosing between a deterministic and a stochastic model is the size of the system. In large systems (e.g. a model with many particles), random fluctuations tend to average out, and the overall behavior of the system can be described accurately by a deterministic model. However, in small systems (e.g. a model with a few number of species), random events can have a significant impact on the system behavior, and make a stochastic model more appropriate.

Another consideration is the computational cost of the simulation. Deterministic models are typically more computationally efficient than stochastic models, as they require fewer simulations to obtain a precise result. Stochastic models, on the other hand, require many simulations to accurately estimate the probability distribution of different outcomes which makes them more computationally intensive.

%% file: part_MC1.tex
\section{Monte Carlo method -- I}\label{sect:mc1}
We begin by introducing the Monte Carlo (MC) algorithm as a stochastic method. This method was invented by {\em John von Neumann} and {\em Stanislaw Ulam} during World War II to improve decision making under uncertain conditions. The name {\em Monte Carlo} is motivated
by the randomness similar to games in the Monte Carlo casino.

In this section, we try to present the basic concepts behind the MC method. Since MC and other stochastic methods rely on random points, the next section is devoted to some techniques for generating random variables. Following that, we will revisit and further consider the MC algorithm in section \ref{sect:mc2}. 
\vsp 

\subsection{Let's play a game}

We play the Snakes-and-Ladders game using a 6-sided dice. 
The game board is shown in Figure \ref{fig:game-board}. 
The game rules are:
\begin{itemize} 
\item {Start from space 1, roll the dice and move forward the number of spaces shown on the dice.}
\item 
{If you land at the base of a ladder move up to the top of the ladder.}
\item
{If you land at the head of a snake slide down to the bottom of the snake.}
\item
{To finish the game you need to roll the exact number to get you to the last space $100$. For example, if you are at space 99 then you should toss the dice until getting 1 to finish.}
\end{itemize}
\vsp 

\begin{center}
\includegraphics[scale=.17]{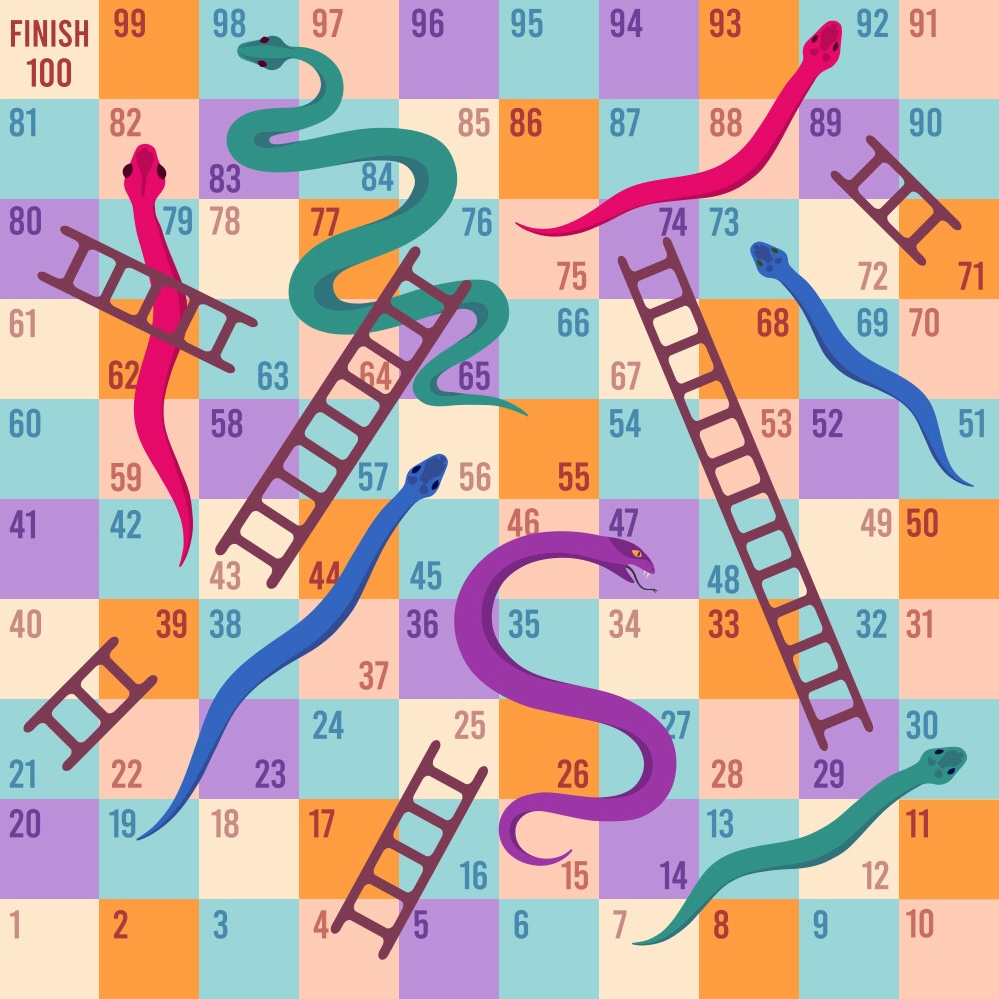}
\captionof{figure}{A Snakes-and-Ladders game board (image from \texttt{www.vectorstock.com}).}\label{fig:game-board}
\end{center}
\vsp 

Let $X$ represent the number of dice rolls required to reach the finish space. Our goal is to determine the expected number of rolls required to finish, or in the mathematics language, to compute the expectation of $X $, denoted as \( \mathbb{E}(X) \).

To answer this, we can ideally find a mathematical expression for the {\em probability density function (pdf)} of the random variable $X$ and then compute $\mathbb{E}(X)$ using the formula for expectation provided in Definition \ref{def:expectation} in the Appendix. However, finding a closed-form expression for the pdf of  $X$ is either impossible or extremely difficult. Instead, we will go for a simple numerical solution.

Assume that one individual plays the game and finishes it after, say, $x_1 = 51$ rolls. 
This single {\em observation} ({\em realization}) is insufficient to conclude that the game will always end after $51$ rolls. 
 However, if $ N$ individuals play the game, we can collect $ N $ observations $x_1, x_2, \dots, x_N $ of $X $. The {approximate expected number of rolls} can then be estimated by calculating the average (mean) of these $N$ observations:
$$
\E(X) \approx \frac{1}{N}(x_1 + x_2 + \cdots + x_N).
$$
This is a Monte Carlo solution: perform the experiment many times and compute the average of the results. As the number of observations $N $ increases, our estimation of the expected number becomes more accurate.

In practice, rather than playing the game manually, we can write a program to simulate it. To simulate the rolling of the dice, we generate random numbers from {\em discrete uniform distribution}  $\mathcal{DU}\{1, 2, \dots, 6\} $, where each outcome has a probability of $1/6$. Refer to the Appendix for definitions and notations of well-known probability distributions. Additionally, the code must account for the effect of snakes and ladders, which can move the player forward or backward on the board.

A simulation performed with $ N = 10,000$ observations shows that the expected number of rolls, \( \mathbb{E}(X) \), is approximately $44$. 
This means that in average people finish this game after approximately $44$ rolls. 
Other possible questions we can answer via the above MC implementation are:
\begin{itemize}
\item What is the probability of finishing the game by exactly $30$ rolls? To answer this, we look at the MC observations, count the number of 30s, and divide it by the total number of observations: 
$$
\pr(X = 30) \approx \frac{\# 30s}{N} \; {\approx 0.02}
$$

\item What is the probability of finishing the game by at most $30$ rolls?
$$
\pr(X \leqslant 30) \approx \frac{\# 30s + \# 29s + \cdots + \#1s}{N}\;
{\approx 0.34 }
$$
\end{itemize}
This approach allows us to approximate the probability based on the frequency of occurrences in the Monte Carlo simulation.
\vsp 
\subsection{A general structure}

There is no single, universally accepted definition of the Monte Carlo method. Its formulation varies based on the underlying mathematical model and the specific type of solution being sought. However, a rather general structure is outlined in Algorithm \ref{alg:mc} below.

\begin{algorithm}
\caption{A general structure of Monte Carlo}\label{alg:mc}
\begin{algorithmic}
\Require Number of observations $N$
\For {$k$ form $1$ to $N$}
    \State Perform one {\em stochastic simulation}
    \State Set $result[k]$ = result of the simulation
\EndFor
\State $FinalResult =$ mean($result$) or other statistical calculations
\end{algorithmic}
\end{algorithm}

The {\em stochastic simulation} can
differ depending on the problem. It could be an observation of a random variable or an observation of a stochastic process. 
In the next section, we will see how to generate random variables from various probability distributions.
Then in section \ref{sect:stoch-proc} we study the random processes.

%% file: part_generation.tex
\section{Random variable generation}\label{sect:generation}
As we pointed out, a typical stochastic simulation requires a set of random numbers, random variables, or a series of stochastic processes.
In this and the next sections we deal with the computer generation of such entities.

We start with generating a uniform random point. Today's random numbers are generated by simple computer algorithms instead of physical devices such as
coin flipping, dice rolling, roulette spinning, and card shuffling, or even modern physical generation methods such as those based on the universal background radiation or the noise of a PC chip.
Random number generation algorithms are usually fast, require little storage space, and can readily reproduce a
given sequence of random numbers. Although such sequences are generated by a deterministic algorithm,
they fulfill main statistical properties of true random sequences.
For this reason the generated numbers are sometimes called {\em pseudorandom} numbers.
Here, we do not pursue the details of algorithms for generating uniform random numbers and just use the following Python function which uses the
\verb+uniform+ function form the \verb+numpy.random+ library.
\begin{shaded}
\vspace*{-0.3cm}
\begin{verbatim}
 U =  np.random.uniform(a, b, size = N)
\end{verbatim}
\vspace*{-0.3cm}
\end{shaded}
\noindent
This generates $N$ uniform random points in interval $[a,b]$. Another possible command is 
\begin{shaded}
\vspace*{-0.3cm}
\begin{verbatim}
 U =  np.random.rand(N)
 U = a + (b-a)*U 
\end{verbatim}
\vspace*{-0.3cm}
\end{shaded}
\noindent
which first generates $N$ uniform random points in standard interval $[0,1]$ and then transfers them into interval $[a,b]$ using the linear map $x \mapsto a + (b-a)x$.

Now, we review some general methods for generating one-dimensional
random variables from a prescribed distribution.
\vsp 

\subsection{Inverse transform method}

Let $X$ be a random variable with pdf $f$ and cdf $F$. If $F$ is continuous and increasing then $F^{-1}$ has the usual definition 
$$
F^{-1}(y):=\{x: F(x) = y\}.
$$
To cover all cases including discrete and nondecreasing cdf functions, the
inverse function $F^{-1}$ can be defined as
\begin{equation}\label{FinvDef}
  F^{-1}(y) = \inf\{x:F(x)\geqslant y\} , \quad 0\leqslant y\leqslant 1.
\end{equation}
See Figure \ref{fig:Finverse} for an illustration. 
\begin{center}
\includegraphics[scale=.45]{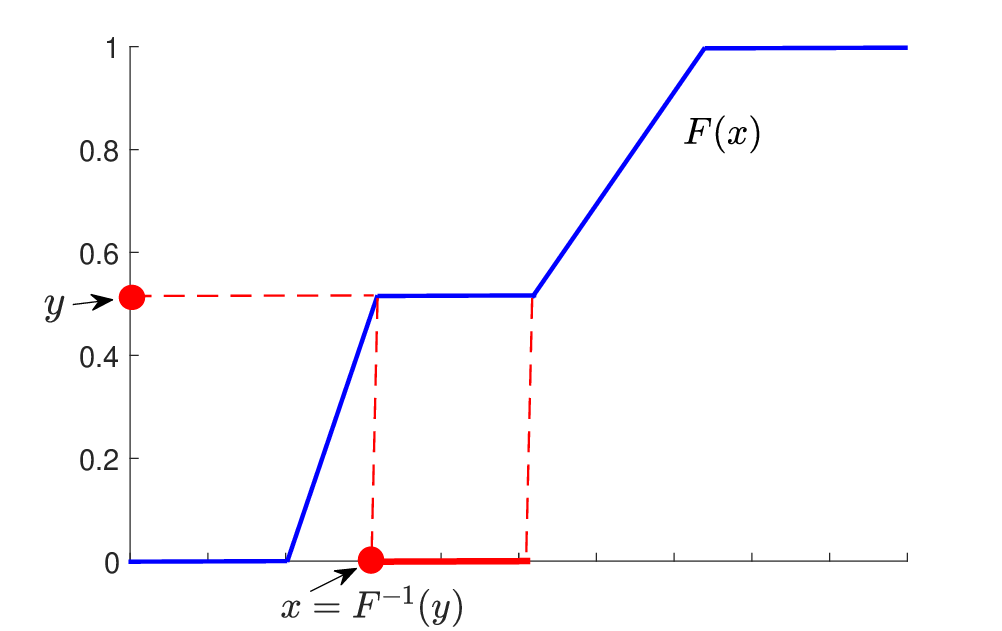}
\captionof{figure}{The inverse of a non-decreasing function}\label{fig:Finverse}
\end{center}

Now we have the following useful result.
\begin{theorem}\label{thm:FinvU}
If $U\sim \Uni(0,1)$ then $X = F^{-1}(U)\sim f$
\end{theorem}
\proof
Since $F$ is invertible and $\pr(U \leqslant u) = u$, we have
$$
\pr(X\leqslant x) = \pr(F^{-1}(U)\leqslant x) = \pr(U\leqslant F(x)) = F(x),
$$
which completes the proof.
$\qed$
\vsp

Theorem \ref{thm:FinvU} proposes a simple algorithm to generate a random variable $X$ with cdf $F$ (or pdf $f$): Generate $U  \sim \Uni(0, 1)$ and set
$X = F^{-1}(U)$. An illustration is given in Figure \ref{fig:itm_illust}, where the uniform variable $U$ on the $y$-axis is transferred to $f$-distributed variable $X$ on the $x$-axis via the cdf $F$.  

\begin{center}
\includegraphics[scale=.28]{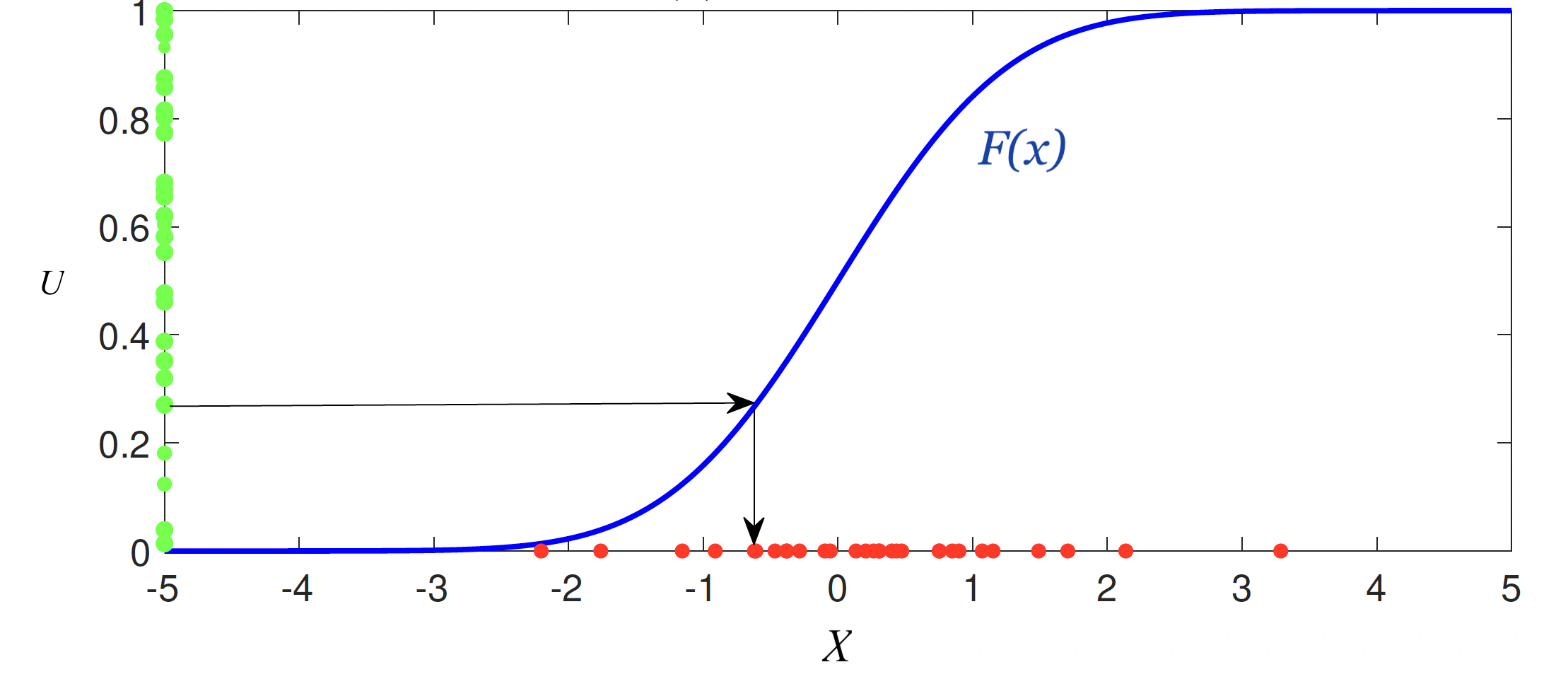}
\captionof{figure}{The inverse transform method}\label{fig:itm_illust}
\end{center}

\begin{example}\label{ex:inv-tran}
To generate a random point from pdf
\begin{equation*}\label{pdf_f2x}
f(x)= \begin{cases}
2x,& x\in[0,1]\\ 0,& \mbox{otherwise},
\end{cases}
\end{equation*}
first we obtain its corresponding cfd
$$
F(x)= \begin{cases}
0, &  x\in(-\infty,0),\\
x^2,& x\in[0,1]\\ 1,& \mbox{otherwise},
\end{cases}
$$
then we generate a uniform variable $U$ and finally we set $X = F^{-1}(U)=\sqrt{U}$.
\end{example}
\vsp

\subsubsection*{Sampling from exponential distribution}
If $X \sim {\Exp}(\lambda)$, then
its pdf $f$ is given by  $f(x)=\lambda \ee^{-\lambda x}$ and
its cdf $F$ by
$$
F(x) = \int_0^x \lambda \ee^{-\lambda y} \td y = 1-\ee^{-\lambda x}, \quad x\geqslant 0.
$$
The inverse of $F$ is $F^{-1}(x) = -\frac{1}{\lambda}\ln (1-x)$. Thus,
to sample from the exponential distribution,
we assume $U\sim \Uni(0,1)$ and set
\begin{equation}\label{exp-gen-for}
X = -\frac{1}{\lambda}\ln (U) \sim \Exp(\lambda).
\end{equation}
Keep in mind that $U\sim\Uni(0,1)$ implies $1-U\sim\Uni(0,1)$.

Here we write a Python function to generate $N$ random variable with exponential distribution.

\begin{shaded}
\vspace*{-0.3cm}
\begin{verbatim}
def RandExp(lam,N):
    # lam: distribution parameter, N: number of requested samples
    U =  np.random.rand(N)  # generate N uniform numbers in [0,1)
    X = -1/lam*np.log(1-U)  # use inverse transform to generate X
    return X
\end{verbatim}
\vspace*{-0.3cm}
\end{shaded}

The following code snippet plots the histogram of $N=500$ generated points with
parameter $\lambda= 0.5$.
\begin{shaded}
\vspace*{-0.3cm}
\begin{verbatim}
import numpy as np
import matplotlib.pyplot as plt
plt.figure(figsize = (5,3))
lam, N = 0.5, 500
X = RandExp(lam,N)
plt.hist(X, bins = 30, histtype = 'bar', color = 'red', density = 'true')
x = np.linspace(0,15,200)
f = lam*np.exp(-lam*x)
plt.plot(x,f,linestyle = '-', color = 'blue')
plt.title('Histogram of $X$ and the pdf $f(x)$')
plt.xlabel('$X$')
plt.ylabel('Frequency %')
\end{verbatim}
\vspace*{-0.3cm}
\end{shaded}

In Figure \ref{fig:exp_gen} the histograms for $\lambda=0.5$ with $N=500$ and $N=5000$ are shown.
For comparison the pdf of the exponential distribution is also plotted in the figures.

 \begin{center}
 \includegraphics[scale=0.68]{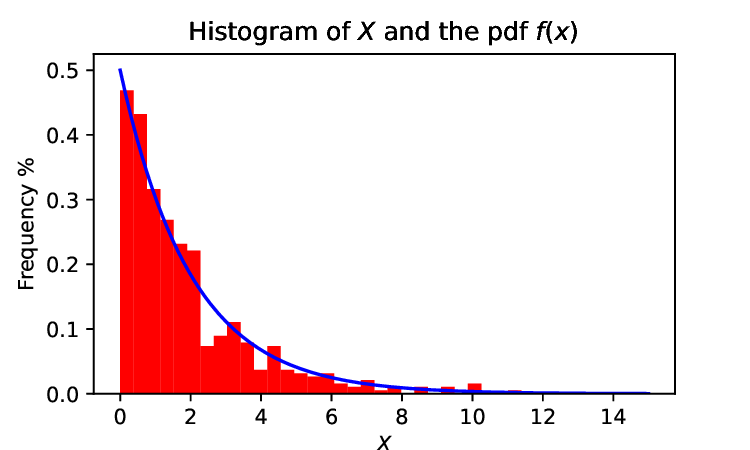}\includegraphics[scale=0.68]{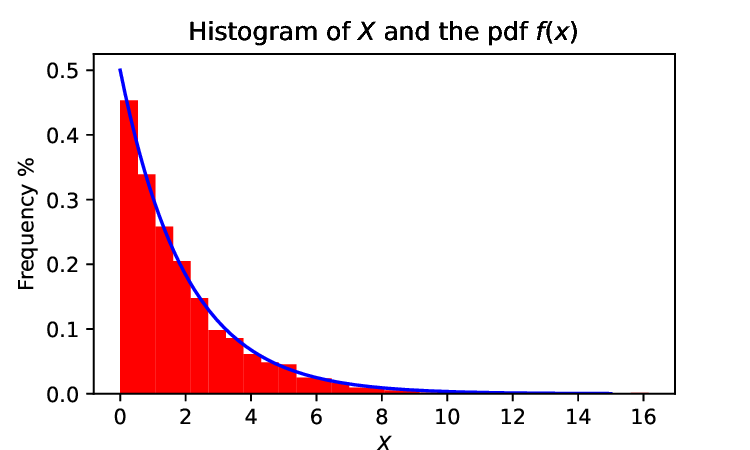}
 \captionof{figure}{Histograms of generated random variables with exponential distribution $\Exp(0.5)$ using the inverse transform method. $N=500$ (left), $N=5000$ (right)}\label{fig:exp_gen}
 \end{center}

%\begin{center}
% \includegraphics[scale=0.75]{stoch4}
% \captionof{figure}{Probability density function of the exponential distribution with $\lambda=0.5$.}\label{fig:exp_pdf}
% \end{center}

\subsubsection*{Sampling from normal distribution}

In the Appendix \ref{sect:appendix_normal_gen} we have shown how a change of variables can help to generate random points from the normal distribution using the inverse transform method. 
However, in our codes in the sequel, we use some built-in functions in Python. We either use
\begin{shaded}
\vspace*{-0.3cm}
\begin{verbatim}
 U =  np.random.normal(mu, sigma2, N)
\end{verbatim}
\vspace*{-0.3cm}
\end{shaded}
\noindent
which generates $N$ normal random points with mean $\texttt{mu}$ and variance $\texttt{sigma2}$, or use 
\begin{shaded}
\vspace*{-0.3cm}
\begin{verbatim}
 U =  np.random.randn(N)
 U = mu + sigma*U 
\end{verbatim}
\vspace*{-0.3cm}
\end{shaded}
\noindent
which first generates $N$ standard normal points (from $\Nor(0,1)$) and then transfers them into a new set of points with distribution $\Nor(\mu,\sigma^2)$ using the fact that if $Z\sim \Nor(0,1)$ then $X = \mu + \sigma Z \sim \Nor(\mu,\sigma^2)$.  

\subsubsection*{Sampling from discrete distributions}

So far, we have only sampled from continuous distributions. However, the inverse transform method can easily be applied to sample from discrete distributions as well. Let $X$ be a discrete distribution with
$$
\pr(X=x_k)  = p_k, \quad k=1,2,\ldots,m \quad \sum_{k=1}^{m} p_k =1.
$$
Without lose of generality, let $x_1<x_2< \cdots<x_m$. The cdf of $X$ is  given by
$$
F(x)  = \sum_{k:x_k\leqslant x}{p_k}.
$$
The plot of $F$ looks like a step-wise function and is shown in Figure \ref{fig:discF}.

\begin{center}
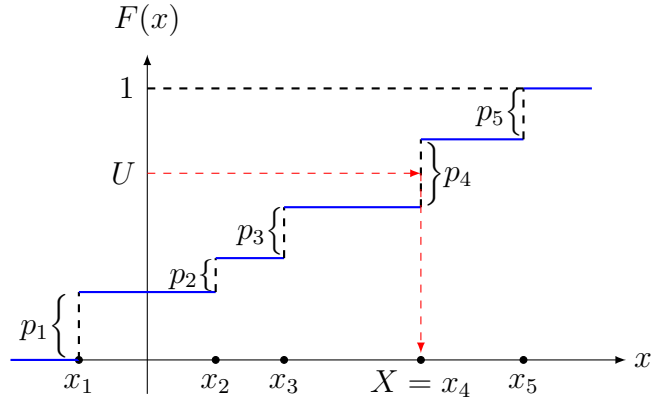

\begin{tikzpicture}[>=latex,scale=0.9]
%\tkzDefPoint(-1,0){xL}
%\tkzDefPoint(7,0){xR}
%\tkzDefPoint(7){xR}
\draw[->] (-1,0) -- (7,0);
\node at (7.25,0) {$x$};
\draw[->] (0,-.5) -- (0,4.5);
\node at (0,5) {$F(x)$};
\node at (-1,-.35) {$x_1$};
\draw[fill=black] (-1,0) circle (0.5mm);
\node at (1,-.35) {$x_2$};
\draw[fill=black] (1,0) circle (0.5mm);
\node at (2,-.35) {$x_3$};
\draw[fill=black] (2,0) circle (0.5mm);
\node at (4,-.35) {$X=x_4$};
\draw[fill=black] (4,0) circle (0.5mm);
\node at (5.5,-.35) {$x_5$};
\draw[fill=black] (5.5,0) circle (0.5mm);
\node at (-.3,4) {$1$};

\draw[thick,blue] (-2,0) -- (-1,0);
\draw[dashed,thick] (-1,0) -- (-1,1);
\node at (-1.5,.5) {$p_1\Big\{$};
\draw[thick,blue] (-1,1) -- (1,1);
\draw[dashed,thick] (1,1) -- (1,1.5);
\node at (0.65,1.255) {$p_2\{$};
\draw[thick,blue] (1,1.5) -- (2,1.5);
\draw[dashed,thick] (2,1.5) -- (2,2.25);
\node at (1.65,1.9) {$p_3\big\{$};
\draw[thick,blue] (2,2.25) -- (4,2.25);
\draw[dashed,thick] (4,2.25) -- (4,3.25);
\node at (4.4,2.75) {$\Big\} p_4$};
\draw[->,dashed,red] (0,2.75) -- (4,2.75);
\node at (-.35,2.75) {$U$};
\draw[->,dashed,red] (4,2.75) -- (4,0.1);
\draw[thick,blue] (4,3.25) -- (5.5,3.25);
\draw[dashed,thick] (5.5,3.25) -- (5.5,4);
\node at (5.15,3.65) {$p_5\big\{$};
\draw[thick,blue] (5.5,4) -- (6.5,4);

\draw[dashed,thick] (0,4) -- (5.5,4);
\end{tikzpicture}
\captionof{figure}{Inverse transform method for a discrete random variable.}\label{fig:discF}
\end{center}

To sample from such discrete random variable, according to definition of $F^{-1}$ in \eqref{FinvDef}, first we generate a uniform distribution
$U\sim \Uni(0,1)$ and then we find the smallest positive integer $k$ such that $U\leqslant F(x_k)$.
Finally $X=x_k$ is reported. Equivalently
$$
X =\begin{cases}
x_1, & \qquad\;\,\, 0< U\leqslant p_1\\
x_2, & \qquad\, p_1< U\leqslant p_1+p_2\\
x_3,  & p_1+p_2 < U\leqslant p_1+p_2+p_3\\
\; \vdots & \qquad \qquad \;\;\vdots\\
\end{cases}
$$

A Python function is given below. The inputs are  the sorted vector \verb+x+ containing the states $x_k$, the corresponding probability vector \verb+p+, and \verb+N+ the number of samples requested. The output is the vector \verb+X+ containing \verb+N+ random points   distributed with $\mathcal{DD}\{[x_1,\ldots,x_m],[p_1,\ldots,p_m]\}$. 

\begin{shaded}
\vspace*{-0.3cm}
\begin{verbatim}
def RandDisct(x, p, N):
    # x: sorted states, p: probabilities, N: number of requested samples
    cdf = np.cumsum(np.array(p))  # compute the cumulative vector
    U = np.random.rand(N)         # generate N uniform numbers in [0,1)
    idx = np.searchsorted(cdf, U) # search U values in cdf intervals
    X = np.array(x)[idx]
    return X      
\end{verbatim}
\vspace*{-0.3cm}
\end{shaded}
\vsp

\begin{example}
Using the code snippet below, we roll a dice $10$ times and report the result. 
\begin{shaded}
\vspace*{-0.3cm}
\begin{verbatim}
x = [1,2,3,4,5,6]
p = [1/6,1/6,1/6,1/6,1/6,1/6]
X = RandDisct(x, p, 10)
print('dice rolls = ', X)
\end{verbatim}
\vspace*{-0.3cm}
\end{shaded}
The result of a run is 
\begin{shaded}
\vspace*{-0.3cm}
\begin{verbatim}
dice rolls =  [5 2 4 4 4 5 6 2 1 3]
\end{verbatim}
\vspace*{-0.3cm}
\end{shaded}
\end{example}
\vsp

To sample from the Bernoulli distribution $\Ber(a)$ for $a\in[0,1]$ we can call the \texttt{RandDisc} function with $\texttt{x = [0,1]}$ and 
$\texttt{p = [1-a,a]}$. We can also use the following independent function to sample a vector of Bernoulli variables with probability $p\in[0,1]$. This function generates a uniform variable $U\sim \Uni(0,1)$, then sets $X=1$ if $U\leqslant p$, and $X=0$ otherwise.

\begin{shaded}
\vspace*{-0.3cm}
\begin{verbatim}
def RandBer(p, N):
    # p: the probability, N: number of requested samples
    X = np.zeros(N)          # set X = [0,0,...,0] initially
    U = np.random.rand(N)    # generate N uniform numbers in [0,1)
    idx = np.where(U <= p)   # find indices for which U is less than p
    X[idx] = 1               # change the corresponding values in X to 1      
    return X       
\end{verbatim}
\vspace*{-0.3cm}
\end{shaded}
\vsp

\begin{remark}
The \verb+RandDisct+ function can be used to generate samples from other discrete distributions. For example,
if $X\sim {\Bin}(p,n)$ with pdf
$$
f(x) = {n\choose x}p^x(1-p)^{n-x},\quad x = 0,1,\ldots,n,
$$
then it is enough to call \verb+X = RandDisct(x,p,N)+ for $\texttt{x} = [0,1,\ldots,n]$ and $\texttt{p} = [f(0),f(2),\ldots,f(n)]$.

There exist other approaches to generate binomial samples. For example one can generate $n$ iid random variables $X_1,\ldots,X_n$ from
${\Ber}(p)$ and set $X = X_1+\cdots+X_n$. 
\end{remark}
\vsp 
\begin{remark}
The \texttt{numpy.random} module provides a variety of built-in functions for generating random samples from commonly used distributions. For discrete distributions, you can use the function \texttt{numpy.random.choice} instead of the \texttt{RandDisct} function.
\end{remark}
\vsp

\subsection{Acceptance-Rejection method}\label{sect:accept-reject}

Usually, the inverse of the cdf $F$ is not available explicitly, and a numerical inversion might be costly and inefficient.
This will make the application of the inverse transform
method limited. The {\em acceptance-rejection} method, introduced by Stan Ulam and John von Neumann, is a more general method that can be used instead.

Suppose that we want to sample from a bounded pdf $f$ which is defined on some finite interval $[a,b]$ and is zero outside this interval.
Suppose further that we have an efficient method for sampling from another random variable with pdf $g$. We follow \cite{Rubinnstein-Kroese:2017}  and for simplicity we first assume that $g(x)=1$ on $[a,b]$ and 
$$
c = \sup \{f(x):x\in [a,b]\}.
$$
The graph of $f$ is clearly under (dominated by) the graph of $cg(x)\equiv c$. See the left panel in Figure \ref{fig:accept_reject}.

\begin{center}
 \includegraphics[scale=0.7]{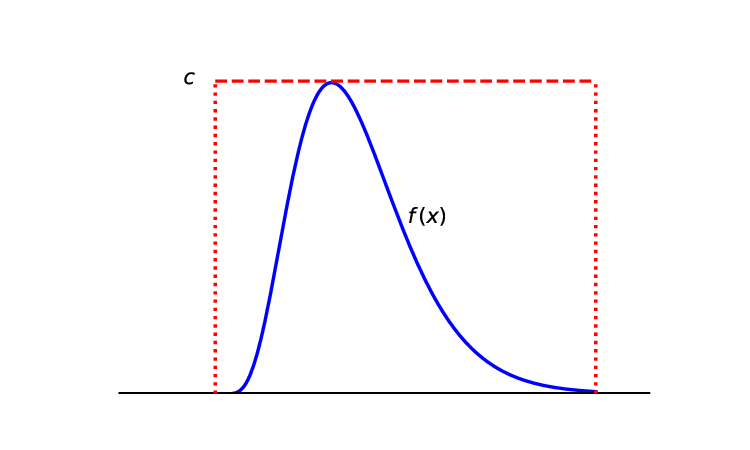}\includegraphics[scale=0.7]{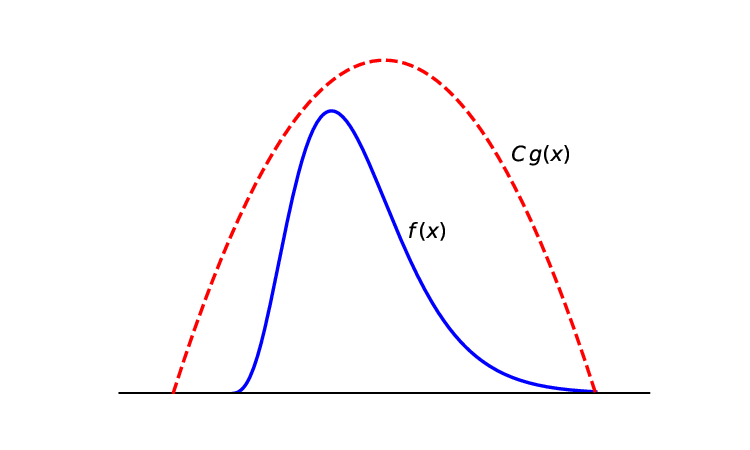}
 \captionof{figure}{Bounding a pdf $f(x)$ by a function $\phi(x)=Cg(x)$.}\label{fig:accept_reject}
 \end{center}

Now, we can generate a random variable $X\sim f$ by using the following algorithm.
\begin{algorithm}
\caption{Acceptance-Rejection Algorithm 1}\label{alg:acceptreject1}
\begin{algorithmic}
\Require Distribution $f$ on interval $[a,b]$, Constant $c$ 
\State 1. Generate $X \sim \Uni (a,b)$
\State 2. Generate $Y \sim \Uni (0,c)$ independent of $X$
\State 3. If $Y\leqslant f(X)$, accept $X$. Otherwise return to step 1.
\Ensure Random point $X$ from distribution $f$
\end{algorithmic}
\end{algorithm}

Since $X$ and $Y$ are uniformly distributed, the pair $(X, Y)$ is uniformly distributed on rectangle $[a, b] \times [0, c]$. In this rectangle, points $(X, Y)$ that lie under the graph of $f$ are accepted according to the criterion in step 3 of  Algorithm \ref{alg:acceptreject1}, while the others are rejected. As a result, the points that are accepted are uniformly distributed in the region under the graph of $f$. This means that the distribution of the accepted values $X$ is $f$.
The proof will come soon. See Figure \ref{fig:accept-reject-points} for an illustration.

\begin{center}
\includegraphics[scale=.25]{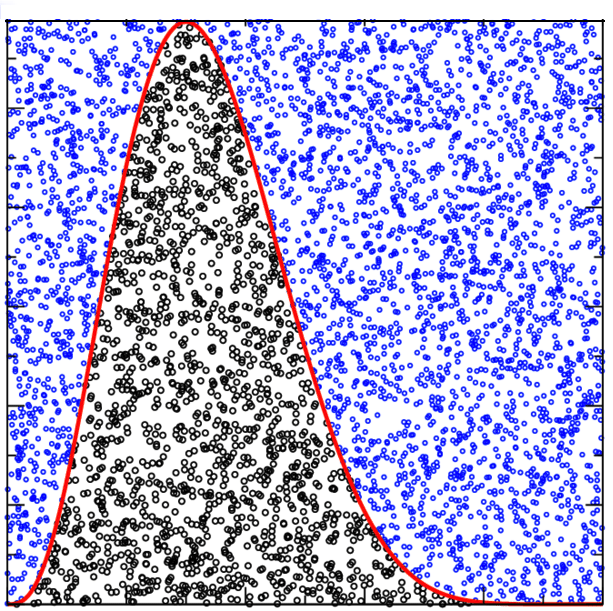}
\captionof{figure}{The graph of pdf $f$ (red curve), the bivariate uniformly distributed points on the rectangle (black and blue points together), and accepted (black) and rejected (blue) points. The $x$-component of accepted points are $f$-distributed.}\label{fig:accept-reject-points}
\end{center}

Sometimes, as we observe from Figure \ref{fig:accept-reject-points}, this algorithm produces many rejected points, which slows down the process. To overcome this inefficiency, we can replace the constant $c$ with a general function $\phi(x) = Cg(x)$, where $g$ is a pdf from which random variables can be easily generated. The pdf $g$ and the constant $C \geqslant 1$  must be chosen carefully so that the graph of $f$ lies under the graph of $\phi$ while ensuring they are as close as possible 
to minimize the number of rejected points thereby improve the efficiency of the algorithm. See the right panel of Figure \ref{fig:accept_reject}. The distribution $g$ is called the {\em proposal distribution}.
The algorithm is given below.

\begin{algorithm}
\caption{Acceptance-Rejection Algorithm 2}\label{alg:acceptreject2}
\begin{algorithmic}
\Require Distribution $f$, Proposal distribution $g$, Constant $C$  
\State 1. Generate $X \sim g$
\State 2. Generate $Y \sim \Uni (0,Cg(X))$ independent of $X$
\State 3. If $Y\leqslant f(X)$, accept $X$. Otherwise return to step 1.
\Ensure Random point $X$ from distribution $f$
\end{algorithmic}
\end{algorithm}

The random variable generated from this procedure has indeed the desired pdf $f$. Because if we
denote the area under graph $\phi(x)=Cg(x)$ by $\mathcal A$ and the area under graph $f(x)$ by $\mathcal B$, then the steps 1 and 2 of the above procedure imply that the random variable $(X,Y)$ is uniformly distributed on $\mathcal A$. To prove this let $h(x,y)$ be the joint pdf of $(X,Y)$. Then we have
$$
h(x,y) = h(y|x) g(x), \quad (x,y)\in \mathcal A.
$$
Item 2 shows that $h(y|x)= \frac{1}{Cg(x)}$ for $y\in[0,Cg(x)]$. Therefore, $h(x,y) = 1/C$ for $(x,y)\in \mathcal A$ which proves that
$(X,Y)$ is uniformly distributed on $\mathcal A$. This shows that an accepted variable $(\tilde X,\tilde Y)$ is uniformly distributed on $\mathcal B$. Since the area of $\mathcal B$ is unity, the pdf of $(\tilde X,\tilde Y)$ is $1$. The marginal pdf of $Z=\tilde X $ on $\mathcal B$ is 
$$\int_{0}^{f(x)}1\, \td y = f(x),$$
which completes the proof.
The efficiency of the algorithm is quantifies by
\begin{equation*}
  \pr((X,Y)\mbox{ is accepted}) = \frac{\mbox{area of }\mathcal{B}}{\mbox{area of }\mathcal{A}} = \frac{1}{C},
\end{equation*}
which means that for a greater efficiency $g$ should be as close as possible to $f$ in order to be able to choose a constant $C$ close to $1$.

In item 2 of Algorithm \ref{alg:acceptreject2} we have $Y\sim \Uni(0,Cg(x))$ which can be rewritten as $Y=U Cg(X)$ where $U\sim \Uni(0,1)$. Using this, item 3 can be replaced by $U\leqslant f(X)/(Cg(X))$. Theses result in a new version for the algorithm.

\begin{algorithm}
\caption{Acceptance-Rejection Algorithm 3}\label{alg:acceptreject3}
\begin{algorithmic}
\Require Distribution $f$, Proposal distribution $g$, Constant $C$  
\State 1. Generate $X \sim g$
\State 2. Generate $U = \Uni (0,1)$ 
\State 3. If $U\leqslant f(X)/(Cg(X))$, accept $X$. Otherwise return to step 1.
\Ensure Random point $X$ from distribution $f$
\end{algorithmic}
\end{algorithm}

A Python code is given here. The input arguments \verb+f+ and \verb+g+ are desired and proposal distributions, \verb+C+ is the constant factor, and \verb+N+ is the number of samples we ask.
The input \verb+gGen+ is an independent function which draws random samples from $g$. Finally, \verb+arg_gGen+ are input arguments
required to execute \verb+gGen+. The output is the vector \verb+Z+ containing $N$ samples with pdf $f$.
\begin{shaded}
\vspace*{-0.3cm}
\begin{verbatim}
def RandAcceptReject(f, g, C, N, gGen, *arg_gGen):
    Z = np.zeros(N)
    for k in range(N):
        reject = True
        while reject:
            X = gGen(*arg_gGen)
            U =  np.random.rand()
            if U <= f(X)/(C*g(X)):
                Z[k] = X
                reject = False
    return Z
\end{verbatim}
\vspace*{-0.3cm}
\end{shaded}
\vsp

\begin{example}

Consider again the pdf $f$ from Example \ref{ex:inv-tran}. To draw samples from this pdf assume that $g(x)=1$ and $C=2$.
Since $f(x)/(Cg(x))=x$, we generate uniform variables $X\sim \Uni(0,1)$ and $U\sim \Uni(0,1)$, and accept $X$ if $U\leqslant X$. If not, we repeat the process until receiving an acceptance. Using the following code we generate $N = 500$ and $5000$ random points from $f$ and plot the histograms. See Figure \ref{fig:f2x_gen}.

\begin{shaded}
\vspace*{-0.3cm}
\begin{verbatim}
import numpy as np
import matplotlib.pyplot as plt
f = lambda x: 2*x
g = lambda x: 1
def gGen(a, b, N):
    U = np.random.uniform(a, b, N)
    return U
arg_gGen = 0, 1, 1  
C, N = 2, 500 
X = RandAcceptReject(f, g, C, N, gGen, *arg_gGen)
plt.figure(figsize = (5, 3))
plt.hist(X, bins = 30, histtype = 'bar', color = 'red', density = 'true')
x = np.linspace(0,1,200)
plt.plot(x,f(x),linestyle = '-', color = 'blue')
\end{verbatim}
\end{shaded}
\begin{shaded}
\vspace*{-0.3cm}
\begin{verbatim}
plt.title('Histogram of $X$ and the pdf $f(x)$')
plt.xlabel('$X$')
plt.ylabel('Frequency %')
\end{verbatim}
\vspace*{-0.3cm}
\end{shaded}
 \begin{center}
 \includegraphics[scale=0.68]{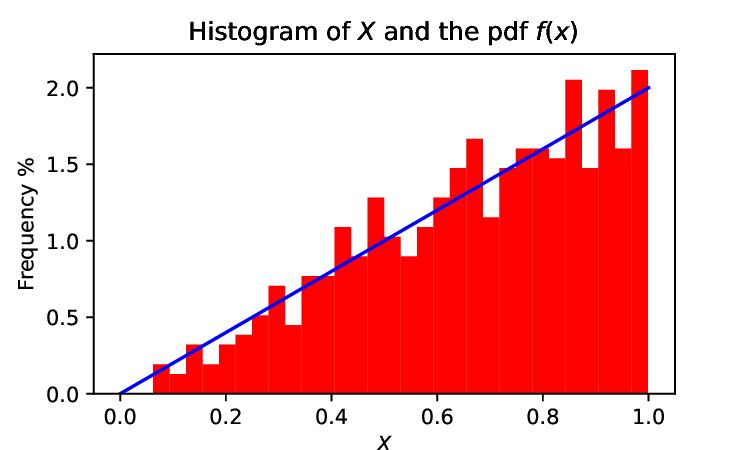}\includegraphics[scale=0.68]{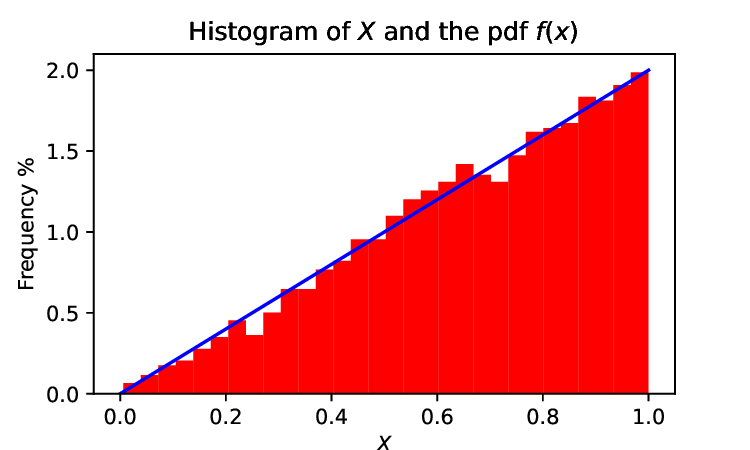}
 \captionof{figure}{Histograms of generated random variables with pdf $f$ in Example \ref{ex:inv-tran} using the acceptance-rejection algorithm. $N=500$ (left), $N=5000$ (right)}\label{fig:f2x_gen}
 \end{center}

\end{example}
\vsp

The acceptance-rejection method can also be used to draw samples from the standard normal distribution.
The positive portion of the normal standard pdf (with $\mu=0$ and $\sigma=1$) can be dominated by a constant $C$ times the pdf of the exponential distribution.
We can generate a positive random variable $X$ from the pdf
\begin{equation*}
  f(x) = \sqrt{\frac{2}{\pi}}\exp(-x^2/2), \quad x\geqslant 0
\end{equation*}
and then assign it a random sign. The sign can be sampled from the Bernoulli distribution. We assume that $g(x)=\exp(-x)$ which is the pdf of $\Exp(1)$, and $C = \sqrt{2e/\pi}$ to have
$f(x)\leqslant Cg(x)$. See Figure \ref{fig:normal_exp}.

\begin{center}
 \includegraphics[scale=0.68]{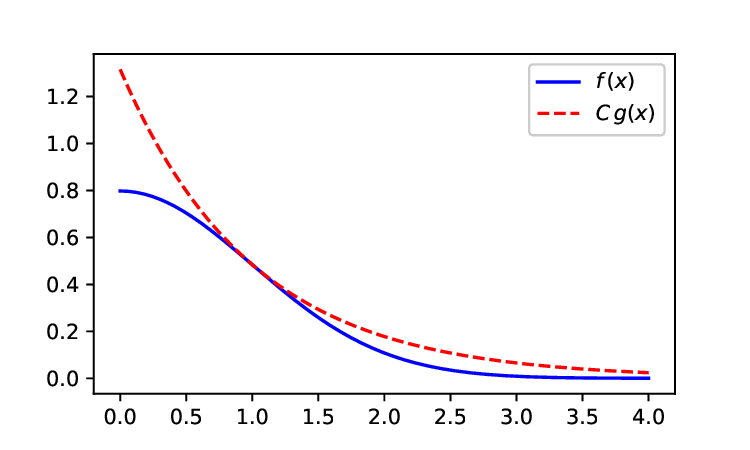}
 \captionof{figure}{Bounding the positive part of the standard normal distribution $f(x)$ by a constant times the exponential distribution $g(x)=\exp(-x)$.}\label{fig:normal_exp}
 \end{center}

The acceptance-rejection algorithm begins with generating a random variable $X\sim {\Exp}(1)$. The acceptance condition then is
$$
U\leqslant\frac{f(X)}{Cg(X)} = \exp(-(X-1)^2/2)),
$$
or equivalently
$$
-\ln U \geqslant \frac{(X-1)^2}{2}. \vsp 
$$
From \eqref{exp-gen-for} we know that $-\ln U\sim {\Exp}(1)$. Thus the last inequality can be rewritten as
$$
V_1\geqslant \frac{(V_2-1)^2}{2}
$$
where $V_1$ and $V_2$ are independent and both of $\Exp(1)$ distribution.

Finally we note that a short section on generating from multivariate distributions (in particular the multivariate normal distribution) is given in Appendix \ref{sect:randvectorgen}. 
\vsp 

%\begin{example}[Beta distribution]\label{ex:beta_gen}
%The pdf of the beta distribution ${\Bet}(\alpha,\beta)$ is given by
%$$
%f(x)=\frac{\Gamma(\alpha+\beta)}{\Gamma(\alpha)\Gamma(\beta)}x^{\alpha-1}(1-x)^{\beta-1}, \quad x\in [0,1], \quad \alpha,\beta\geqslant0.
%$$
%The case ${\Bet}(1,1)$ is identical with $\Uni(0,1)$.
%For general values $\alpha>0$ and $\beta>0$, it is difficult to use the inverse transform method directly.
%We first consider the case where either
%$\alpha=1$ or $\beta=1$. In that case, we can simply use the inverse transform method. For example if
%$\beta=1$ then the $\Bet(\alpha, 1)$ pdf is
%$$
%f(x) = \alpha x^{\alpha-1}, \quad x\in [0,1].
%$$
%The cdf then can be simply obtained as
%$$
%F(x) = \int_0^x \alpha  y^{\alpha-1}dy = x^{\alpha}, \quad x\in [0,1].
%$$
%with the inverse function
%$$
%F^{-1}(x) = x^{1/\alpha}.\vsp
%$$
%Therefore, for a uniform random variable $U$ in $[0,1]$ we have
%$$
%X=U^{1/\alpha} \sim {\Bet}(\alpha,1).
%$$
%For case $\alpha=1$ we can similarly show that
%$$
%X = 1-U^{1/\beta} \sim {\Bet}(1,\beta).
%$$
%For the general case we use the fact that if $Y\sim{\Gam}(\alpha,1)$ and $Z\sim{\Gam}(\beta,1)$ then
%$$
%X = \frac{Y}{Y+Z}\sim {\Bet}(\alpha,\beta).
%$$
%This requires the generation of random variables with Gamma distributions.
%
%\end{example}

%% file: part_MC2.tex
\section{Monte Carlo method -- II}\label{sect:mc2}

As we observed in section \ref{sect:mc1}, Monte Carlo uses random points to estimate the mean of (complicated) random variables/processes. For
continuous random variables it is equivalent to solving certain integrals. If \( X \) is a continuous random variable with pdf \( f(x) \), and \( g(X) \) is some function of \( X \), then \( g(X) \) becomes a new random variable. The expectation of \( g(X) \) is given by (refer to \eqref{def:expectation_g} in the Appendix)
\[
\E[g(X)] = \int_{-\infty}^{\infty} g(x) f(x) \, \td x.
\]
The Monte Carlo method is an stochastic tool to approximate such integrals numerically, particularly when it is hard to apply the available deterministic methods (e.g. in high-dimensional spaces or on complex domains). The process is called {\em Monte Carlo integration}. 
Consider the generic integral
\[
\int_a^b g(x) f(x) \, \td x \vsp
\]
where \( f(x) \) is a pdf associated with a random variable \( X \), and \( g(x) \) is some function, often referred to as the {\em performance function}. This integral can be interpreted as the expectation \( \E[g(X)] \) of the random variable \( g(X) \), where \( X \) is distributed according to the pdf \( f(x) \).

The expected value can be estimated by drawing random samples from  distribution \( f(x) \), evaluating \( g(x) \) at these random points, and then averaging the results. Specifically, the Monte Carlo procedure follows these steps:
\begin{itemize}
\item[1.] Generate \( N \) random samples \( x_1, x_2, \ldots, x_N \) from pdf \( f(x) \).
\item[2.] For each sample \( x_k \), compute the corresponding value \( g(x_k) \).
\item[3.] Approximate the integral by computing the average of these values:
   \[
   \frac{1}{N}\sum_{k=1}^{N} g(x_k) =: \overline{g}_N .
   \]
\end{itemize}
As \( N \) increases, this approximation converges to the true value of the integral, thanks to the law of large numbers. The convergence proof will come soon. 
\vsp 

\begin{example}
Consider the following 1D integral
\[
I = \int_0^1 g(x) \, dx.
\]
This integral can be interpreted as the expectation of \( g(X) \), where \( X \) is a uniformly distributed random variable on the interval \([0, 1]\). In other words
\[
I = \E[g(X)] = \int_{\blue{0}}^{\blue{1}} g(x) \cdot {\blue{1}} \, dx, \quad  X \sim \Uni(0,1).
\]
The Monte Carlo method generates \( N \) random points \( x_1, x_2, \ldots, x_N \) from distribution  \(\Uni(0, 1)\), and 
computes
$$
I \approx \frac{1}{N} \left( g(x_1) + g(x_2) + \cdots + g(x_N) \right)
$$
A code snippet is given below for $g = \sin x$. 
\begin{shaded}
\vspace*{-0.3cm}
\begin{verbatim}
import numpy as np
def g_fun(x):
    return np.sin(x)               # the integrand g(x) = sin(x)
int_exact = 1-np.cos(1)            # the exact value for comparison
int_mc = np.zeros(6)            
for k in range(6):
    N = 10**k                      # number of points from 1, 10,..., 10^5 
    X = np.random.uniform(0, 1, N) # generate uniform random points in [0,1]
    g = g_fun(X)                   # evaluate the integrand on X
    int_mc[k] = np.mean(g)         # take mean
print('int_mc = ', 
      np.round(np.abs(int_exact-int_mc),5)) # errors, rounded to 5 decimals
\end{verbatim}
\vspace*{-0.3cm}
\end{shaded}
\vsp
An execution gives 
\begin{shaded}
\vspace*{-0.3cm}
\begin{verbatim}
int_mc =  [0.03652 0.08578 0.02132 0.00362 0.00203 0.00063]
\end{verbatim}
\vspace*{-0.3cm}
\end{shaded}
A new execution will result in a new (different) error vector as the integral points are generated randomly. 
\begin{shaded}
\vspace*{-0.3cm}
\begin{verbatim}
int_mc =  [0.03011 0.00532 0.00749 0.01049 0.00279 0.00023]
\end{verbatim}
\vspace*{-0.3cm}
\end{shaded}
In any case, 
the accuracy improves as the number of random samples \( N \) increases. However, the convergence speed is slow. 

In general to estimate the integral of a function $g$ on finite interval $[a,b]$ the Monte Carlo integration is applied as below:
$$
\int_a^b g(x) dx =\blue{ (b-a)}\int_a^b g(x)\blue{\frac{1}{b-a}} dx \approx \underbrace{(b-a)}_{width} \underbrace{\frac{1}{N}\sum_{k=1}^N g(x_k)}_{mean~ hight}, \quad x_k\in \Uni(a,b)
$$

See also Figure \ref{fig:mc_int} for an illustration. 

\begin{center}
\qquad\qquad
\begin{tikzpicture}[>=latex,scale=.9]

\draw[->] (-2,0) -- (9,0);
\draw[->] (-1,-.5) -- (-1,4.5);
\node at (.5,-.35) {\small $a$};
\draw[fill=black] (.5,0) circle (0.5mm);
\node at (1.5,-.35) {\small $x_1$};
\draw[fill=black] (1.5,0) circle (0.5mm);
\node at (2,-.35) {\small $x_2$};
\draw[fill=black] (2,0) circle (0.5mm);
\node at (3,-.35) {\small $x_3$};
\draw[fill=black] (3,0) circle (0.5mm);
\node at (5,-.35) {\small $x_4$};
\draw[fill=black] (5,0) circle (0.5mm);
\node at (6.5,-.35) {\small $x_5$};
\draw[fill=black] (6.5,0) circle (0.5mm);
\node at (7.5,-.35) {\small $b$};
\draw[fill=black] (7.5,0) circle (0.5mm);

\draw [decorate,
    decoration = {brace,mirror},thick] (.5,-1) --  (7.5,-1); 
\node at (4,-1.5) {\small $width$};
    
\draw[thick,scale=1, domain=0:8, smooth, variable=\x, red] plot ({\x}, {sin(deg(\x))+3});
\node[red] at (8.6,4) {\small $g(x)$};

\draw[dashed,thick,blue] (1.5,0) -- (1.5,4); % x1
\draw[fill=blue] (1.5,4) circle (0.5mm);
\draw[dashed,thick,blue] (2,0) -- (2,3.9); % x2
\draw[fill=blue] (2,3.9) circle (0.5mm);
\draw[dashed,thick,blue] (3,0) -- (3,3.14); % x3
\draw[fill=blue] (3,3.14) circle (0.5mm);
\draw[dashed,thick,blue] (5,0) -- (5,2); % x4
\draw[fill=blue] (5,2) circle (0.5mm);
\draw[dashed,thick,blue] (6.5,0) -- (6.5,3.2); % x5
\draw[fill=blue] (6.5,3.2) circle (0.5mm);
\draw[thick,blue] (.5,0) -- (.5,3.48); % a
\draw[thick,blue] (7.5,0) -- (7.5,3.93); % b

\draw[->] (1.6,4.1) -- (3.3,4.7);
\draw[->] (2.1,4) -- (3.7,4.7);
\draw[->] (3.1,3.24) -- (3.9,4.5);
\draw[->] (4.9,2.1) -- (4.2,4.5);
\draw[->] (6.4,3.3) -- (4.5,4.5);
\node at (4,5) {\small $hights$};

\end{tikzpicture}

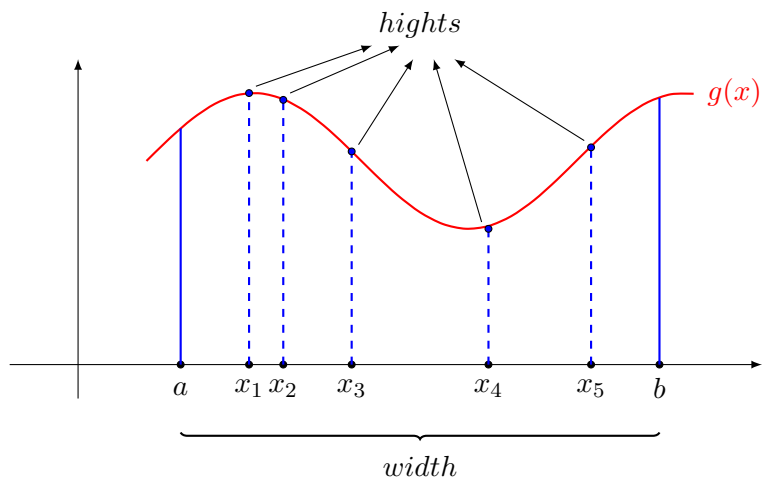
\captionof{figure}{Schematic of 1D Monte Carlo integration}\label{fig:mc_int}
\end{center}

\end{example} 
\vsp

\begin{example}
To estimate the value of integral
$$
 I = \int_{-\infty}^\infty (x^4-x+1) \ee^{-x^2/2} \td x
$$
using Monte Carlo method, we can write
$$
I =  \sqrt{2\pi}\int_{\blue{-\infty}}^{\blue{\infty}}\underbrace{(x^4-x+1)}_{g(x)} \underbrace{\blue{\frac{1}{\sqrt{2\pi}} \ee^{-x^2/2}}}_{f(x)} dx
$$
where $f(x)$ is the pdf of the standard normal distribution on $(-\infty,\infty)$. The estimate then is  
$$
I\approx \sqrt{2\pi}\times \frac{1}{N}\sum_{k=1}^N(x_k^4-x_k+1) 
$$
where $x_k$ are generated from $\Nor(0,1)$.

\end{example}
\vsp

\subsection{Convergence of Monte Carlo integration}

Before providing an analysis for convergence of the Monte Carlo method, we first compare its convergence rate with a deterministic method for computing a one-dimensional integral.  As a simple deterministic integration quadrature consider the mid-point (MP) rule:
\begin{align*}
\int_{0}^1 g(x) \td x& = h\left[g(x_1^*)+g(x_2^*)+\cdots +g(x_N^*)\right] + \mathcal O(h^2)\\
&=\frac{1}{N}\left[g(x_1^*)+g(x_2^*)+\cdots +g(x_N^*)\right] + \mathcal O(N^{-2})
\end{align*}
where $h = 1/N$ and $N$ is the number of integration points in the interval $[0,1]$, and $x_k^* = (x_{k-1}+x_k)/2$ are mid points. 
Compared to the Monte Carlo method, the integration points $x_k^*$ in the mid-point rule are a set of {\em equidistance} points. 
For a sufficiently smooth function, say $g\in C^2[a,b]$, the order of convergence of this method is $2$ --- if $h$ is halved (equivalently if $N$ is doubled) then the error is quartered. The error plot for both methods are given in Figure \ref{fig:mpmc-compare}. For Monte Carlo, the results of ten simulations are depicted and the approximate average rate is computed, which is close to $0.5$. 

\begin{center}
\includegraphics[scale=.55]{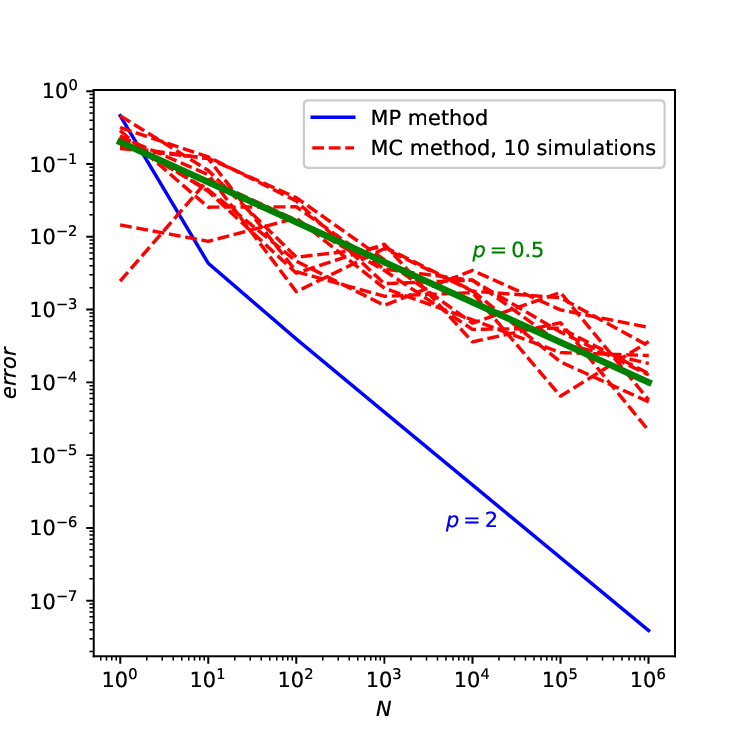}
\captionof{figure}{Convergence orders of Monte Carlo (MC) and mid-point (MP) methods for integration.}\label{fig:mpmc-compare}
\end{center}

The convergence of the mid-point rule is of order \( h^2 \), or equivalently \( N^{-2} \). When we extend this rule to a two-dimensional integral on a rectangular domain, the convergence remains at order \( h^2 \), which translates to \( N^{-1} \) since, in two dimensions, \( h \) is proportional to \( N^{-1/2} \). In general, for a \( d \)-variate integral, the convergence of the mid-point method is of order \( N^{-2/d} \) in the {\em maximum norm}. 
This reduction in the convergence order with increasing dimensions also occurs in other deterministic methods, such as the trapezoidal method and Simpson's method.
However, as we will see, the convergence of the Monte Carlo method is {\em almost surely} of order \( N^{-1/2} \), {\em independent of the dimension}. This indicates that while deterministic methods, if practically applicable, are preferred for low-dimensional integrals, the Monte Carlo method is often a better choice for high-dimensional integrals.

In order to determine how accurate the Monte Carlo solution is, we 
assume that $X$ is a random variable with pdf $f$ and $X_1, X_2, \ldots, X_N$ is a sample (a set of independent and identically distributed (iid) random variables) from $X$. For a function $g$, assume that $g(X)$, as another random variable, has the mean $\mu$ and variance $\sigma^2$, i.e.,  
$$
\mu = \E_f[g(X)], \quad \sigma^2 = \Var_f[g(X)]. 
$$
Since, $X_k$ are iid, all $g(X_k)$ have the same mean $\mu$ and variance $\sigma^2$. 
The new random variable
\begin{equation}\label{mc:basic_estimate}
Y = \frac{1}{N}\sum_{k=1}^{N}g(X_k), \quad X_k\sim f
\end{equation}
is an {\em unbiased} estimator for $\mu = \E[g(X)]$ in the sense that 
$$
\E(Y)=\frac{1}{N}\sum_{k=1}^{N}\E[g(X_k)] = \frac{N\mu}{N} = \mu.
$$ 
What can we say about the variance of $Y$? Using the {\em central limit theorem} (see Appendix \ref{sect:limit-thm}), for sufficiently large values of $N$ we have 
$$
{Y}\sim \Nor(\mu,\sigma^2/N) \;\;\mbox{   or   }\;\; \sqrt N\frac{{Y}-\mu}{\sigma}\sim \Nor(0,1).
\vsp
$$
This shows that the variance of $ Y$ is $\sigma^2/N$, proving that $Y$ approaches $\mu$ in probability with convergence rate $\mathcal O(1/\sqrt N)$.
More precisely, 
assume that $\Phi$ denotes the standard normal cdf and $z_\gamma$ denotes the $\gamma$-quantile of $\Nor (0,1)$, i.e.,
$\Phi(z_\gamma) = \gamma$. This means that if $Z\sim \Nor(0,1)$ then
$$
\pr(-z_{1-\alpha/2}\leqslant Z \leqslant z_{1-\alpha/2}) = 1-\alpha.
$$
Thus, we can write
$$
\pr\left(-z_{1-\alpha/2}\leqslant \frac{\sqrt N({Y}-\mu)}{\sigma} \leqslant z_{1-\alpha/2}\right) = 1-\alpha.
$$
or equivalently
$$
\pr\left(\mu-z_{1-\alpha/2}\frac{\sigma}{\sqrt N}\leqslant  Y \leqslant \mu+z_{1-\alpha/2}\frac{\sigma}{\sqrt N}\right) = 1-\alpha. \vsp
$$
In other words, with probability $(1-\alpha)100\%$ the {\em confidence interval} for $ Y$ is
$$
\left[\mu-z_{1-\alpha/2}\frac{\sigma}{\sqrt N}\; ,\; \mu+z_{1-\alpha/2}\frac{\sigma}{\sqrt N}\right].
$$
\begin{figure}[ht!]
\begin{center}
  \includegraphics[scale=0.13]{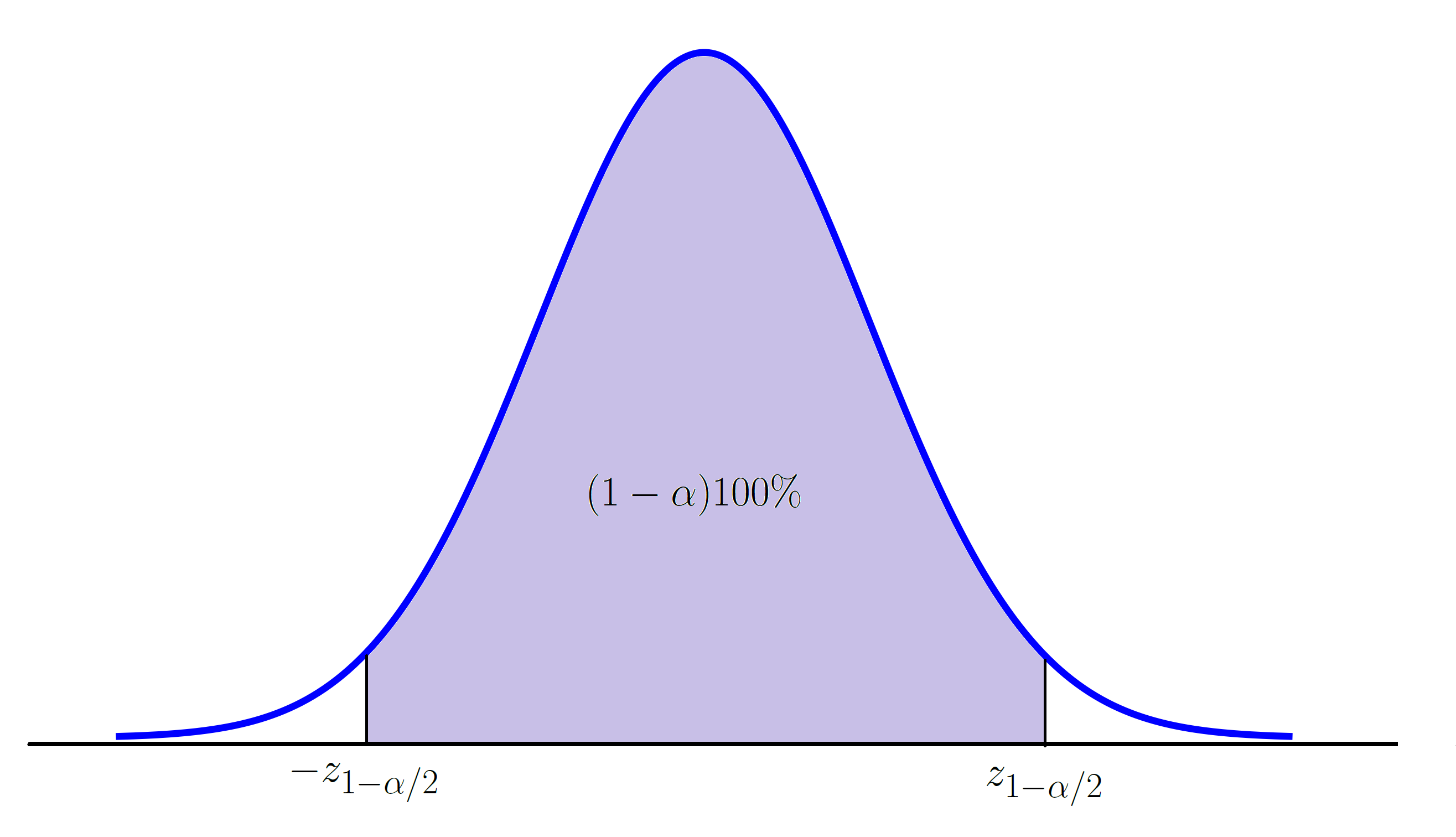}
\end{center}
\caption{Confidence interval in the standard normal distribution}\label{fig:confid_interval}
\end{figure}
See Figure \ref{fig:confid_interval}.
For example, if $\alpha = 0.05$ then $z_{1-\alpha/2}=z_{0.975}\doteq 1.96$, which means that with $95\%$ probability $ Y$ falls in interval
$$
\left[\mu-1.96\frac{\sigma}{\sqrt N}\; ,\; \mu+1.96\frac{\sigma}{\sqrt N}\right].
$$
This probability will increase to $99\%$ if we replace the factor $1.96$ in the confidence interval by $z_{0.995}\doteq 2.576$,
and to $0.999\%$ by $z_{0.9995}\doteq 3.29$. In general, the accuracy of the estimator ${Y}$ is determined by its standard deviation, i.e., $\sigma/\sqrt N$.

We again emphasis that $\mu$ and $\sigma^2$ are the mean and the variance of random variable $g(X)$.
Our aim was to estimate $\mu$ but for error estimation the value of $\sigma$ is also required.
Usually, $\sigma^2$ is unknown, but can be estimated with the {\em sample variance}
\begin{equation*}
  s_N^2 = \frac{1}{N-1}\sum_{k=1}^{N}(g(x_k)-\overline g_N)^2, \quad \overline g_N = \frac{1}{N}\sum_{k=1}^N g(x_k)
\end{equation*}
where $x_k$ are generated from pdf $f$. The sample variance $s_N^2$
 tends to $\sigma^2$ by the law of large numbers. Consequently, for large values of $N$, the {\em approximate} confidence interval
for $ Y$ is
$$
\left[\mu-z_{1-\alpha/2}\frac{s_N}{\sqrt N}\; ,\; \mu+z_{1-\alpha/2}\frac{s_N}{\sqrt N}\right].
$$
We must be aware that the confidence interval can be trusted as far as $s_N^2$ is a proper estimate for the variance $\sigma^2$.
\vsp

%we assume that $f(x) = 0.5\ee^{-0.5x}$ (the pdf of the exponential distribution with $\lambda=0.5$) and $H(x)= 2(x^2-x)$. Thus, we generate Exp$(0.5)$ variables $X_1,\ldots,X_N$ (for example using \verb+ExpGen+ function) and estimate $I$ by
%$$
%\overline I_N = \frac{1}{N} \sum_{k=1}^{N}2(X_k^2-X_k).
%$$
%The results for different values of $N$ are given in Table \ref{tb:mc_exp}. The exact value of the integral is $12$.
%\begin{center}
%\begin{tabular}{llllll}
%  \hline
%  % after \\: \hline or \cline{col1-col2} \cline{col3-col4} ...
%   & $N$ & $\quad{\overline} I_N$ & $|I-{\overline} I_N|$ & $\quad s_N^2$ & $s_N/\sqrt{N}$\\
%   \hline
%   & $10^1$ & $5.08495$ & $6.91504$ & $163.53$ & $4.04391$ \\
%   & $10^2$ & $8.73476$ & $3.26523$ & $431.56$ & $2.07740$ \\
%   & $10^3$ & $11.2824$ & $0.71751$ & $835.18$ & $0.91388$ \\
%   & $10^4$ & $12.3201$ & $0.32017$ & $1154.6$ & $0.33980$ \\
%   & $10^5$ & $11.9608$ & $0.03911$ & $1004.3$ & $0.10021$ \\
%   & $10^6$ & $12.0040$ & $0.00401$ & $1044.8$ & $0.03232$ \\
%   \hline
%   & Exact & $12$ &            & $1040$ &  \\
%  \hline
%\end{tabular}
%\captionof{table}{Monte Carlo integration.}\label{tb:mc_exp}
%\end{center}
%\end{example}

\begin{example}\label{ex:mc_normalcdf}
To estimate the cdf of the standard normal distribution, i.e.,
\begin{equation}\label{cdf-normalPhi}
\Phi(t) = \int_{-\infty}^{t}\frac{1}{\sqrt{2\pi}}\ee^{-x^2/2}\td x
\end{equation}
using the Monte Carlo method, we generate $N$ normal variable $X_1,\ldots,X_N\sim \Nor(0,1)$ and set
\begin{equation}\label{mc:basic_est}
\overline \Phi(t) = \frac{1}{N}\sum_{k=1}^{N}\mathbb I_{X_k\leqslant t},
\end{equation}
with (exact) variance $\Phi(t)[1-\Phi(t)]/N=:\sigma^2/N$, since the variables $\mathbb I_{X_k\leqslant t}$ are independent
Bernoulli variables with success probability $\Phi(t)$.
Here, $\mathbb I_A = 1$ or $0$ when $A$ is true or false, respectively.  The confidence interval tells that with $(1-\alpha)100\%$ probability the error is at most 
$z_{1-\alpha/2}\sigma/\sqrt N$ (or approximately $z_{1-\alpha/2}s_N/\sqrt N$). For example, for $t=0$ we have $\Phi(t)=0.5$ and
$\sigma^2=\Phi(t)[1-\Phi(t)]=1/4$. For this special case, to achieve a precision of three decimals with probability $95\%$, we set
$\alpha = 0.05$, giving $z_{1-\alpha/2}\doteq 1.96$, and choose $N$ such that
$$
\frac{z_{1-\alpha/2}\sigma}{\sqrt N} \doteq \frac{1.96}{2\sqrt{N}}\leqslant \frac{1}{2}10^{-3}.
$$
which gives 
$
N > 3.85\cdot 10^6
$.
Experimental results for different values of $N$ and $t$ are obtained by executing the following code.
\begin{shaded}
\vspace*{-0.3cm}
\begin{verbatim}
import numpy as np
from scipy.stats import norm
t, K = 0, 7
Phi, PhiBar = norm.cdf(t), np.zeros(K)
for k in range(K):
    N = 10**(k+1)
    X = np.random.normal(0,1,N)
    PhiBar[k] = 1/N*len(X[(X < t)])
Error = abs(PhiBar-Phi)/Phi
print('Phi(',t,') = ',Phi,'\n','PhiBar = ',PhiBar,'\n','Error = ',Error)
\end{verbatim}
\vspace*{-0.3cm}
\end{shaded}

Here, an output of this code for $t=0$ is given.
We again note that $\Phi(0)=0.5$.
The last estimation $0.5002081$ with relative error $0.0004162$ corresponds to $N = 10^7$ standard normal samples, which confirms our
error estimation above.
\begin{shaded}
\vspace*{-0.3cm}
\begin{verbatim}
Phi(0) =  0.5
PhiBar =  [0.7  0.46  0.516  0.4983  0.50105  0.500231  0.5002081]
Error  =  [0.4  0.08  0.032  0.0034  0.00210  0.000462  0.0004162]
\end{verbatim}
\vspace*{-0.3cm}
\end{shaded}
However, if we execute the code for $t=-4.5$, i.e., to estimate $\Phi(-4.5) \doteq 3.39767\times10^{-6}$, the following results will be obtained.
\begin{shaded}
\vspace*{-0.3cm}
\begin{verbatim}
Phi(-4.5) =  3.3976731247300535e-06
PhiBar =  [0.0e+00 0.0e+00 0.0e+00 0.0e+00 0.0e+00 4.0e-06 3.6e-06]
Error  =  [1.      1.      1.      1.      1.      1.8e-01 6.0e-02]
\end{verbatim}
\vspace*{-0.3cm}
\end{shaded}
In this experiment, we observe that for values of \( N \) up to \( 10^5 \), the estimated values remain at \( 0 \), resulting in errors of \( 100\% \). This outcome indicates that none of the generated random points fall before \( x = -4.5 \). For larger values of \( N \), such as \( 10^6 \) and \( 10^7 \), the estimations are nonzero but still quite poor.
The issue lies in the fact that we are attempting to estimate the probability of a very rare event. To have enough samples from \( f \) in such sub-domains ($x\leq -4.5$, the left tail of normal distribution), a huge number of samples is needed in the whole domain. This is a drawback of the basic Monte Carlo method. However, by applying a modification we can significantly improve the accuracy with fewer samples. This is discussed in the next subsection.
\end{example}
\vsp
\begin{workout} The irrational number $\pi$ is the volume of the unit ball in $\R^2$. We can compute the area of the first quadrant sector and multiply it by 4:
$$
\pi = 4\int_{0}^{1}\sqrt{1-x^2}\td x.
$$
Use Monte Carlo to estimate $\pi$ by approximating the above integral with different number of samples. Compare the errors with the Monte Carlo error bound. 
Plot the errors in log-log scale for $10$ executions.  

\end{workout}
\vsp
\begin{workout}
Estimate the integral
$$
I = \int_{0}^{\infty}\ee^{-0.5x}(x^2-x)\td x \vsp
$$
using the Monte Carlo integration with different values of $N$.
Identify the right distribution to sample from for this estimation. You can compute the exact value of this integral using the integration by parts. Compare the exact value and the Monte Carlo solutions. Are your results confirmed by the error bound of the Monte Carlo method? 
\end{workout}
\vsp
\subsection{Importance sampling}
The importance sampling approach is based on the principle that the expectation
of $g(X)$ with respect to density $f$ can be written in the alternative form
\begin{equation}\label{mc:importanceE}
\E_f[g(X)] = \int g(x)f(x)\td x  =\int g(x)\frac{f(x)}{\ell(x)}\ell(x)\td x = \E_\ell\left[g(X)\frac{f(X)}{\ell(X)}\right]
\end{equation}
where $\ell(x)$ is another density function, called the {\em importance sampling function} or
{\em envelope}.
Equation \eqref{mc:importanceE} suggests to estimate $\E_f[g(X)]$ by drawing
iid random variables  $X_1, \ldots ,X_N$ from $\ell$ (not $f$!) and use the estimator
\begin{equation}\label{mc:importanceEstimate}
  Y_{\mathrm{IS}} = \frac{1}{N}\sum_{k=1}^{N} g(X_k)\frac{f(X_k)}{\ell(X_k)},\quad X_k \sim \ell
\end{equation}
instead of basic estimator \eqref{mc:basic_estimate}. Here, the subscript `IS' stands for `Importance Sampling'.   
For this strategy to work, it must be easy to sample from $\ell$ and to
evaluate $f$, even when it is not easy to sample from $f$.
Indeed, \eqref{mc:importanceEstimate} does converge to
$\E_f[g(X)]$ for the same reason the basic Monte Carlo estimator converges,
whatever the choice of the distribution $\ell$ is, as long as $\mathrm{supp}(\ell)\subset \mathrm{supp}(g \times f)$. This assumption on the support of $\ell$ is important because a smaller support truncates the integral \eqref{mc:importanceE} and produces
a biased result.
The ratio of densities is denoted by
$$
w(x) = \frac{f(x)}{\ell(x)},
$$
and is called the {\em likelihood} ratio. This ratio needs only to be known up to a constant; let say $w(x)= c\tilde w(x)$. This is useful, for example, when the distribution $f$ or $\ell$ is known only up to a constant.
Since $\E_\ell[w(X)]=1$ we can write
\begin{equation*}
  \E_f[g(X)] = \E_\ell\left[g(X)w(X)\right] = \frac{\E_\ell\left[g(X)w(X)\right]}{\E_\ell\left[w(X)\right]},
\end{equation*}
which motivates the {\em weighted sample estimator}
\begin{equation}\label{eq:weightedimportance}
 \frac {\ds \frac{1}{N}\sum_{k=1}^{N} g(X_k)w(X_k)}{\ds \frac{1}{N}\sum_{k=1}^{N} w(X_k)} = \frac{\ds \sum_{k=1}^{N} w_k g(X_k)}{\ds \sum_{k=1}^{N} w_k}=: Y_{\mathrm{IS}}^w
\end{equation}
where $w_k = \tilde w(X_k)$. Note that the estimators $Y_{\mathrm{IS}}$ and $Y_{\mathrm{IS}}^w$ are mathematically the same but $Y_{\mathrm{IS}}^w$ can also be used in situations where either $f$ or $\ell$ is missing a normalizing constant. Both number $N$ and constant $c$ are cancelled from the numerator and denominator in \eqref{eq:weightedimportance}.
\vsp 

\begin{example}\label{ex:is_mc_normalcdf}
Coming back to Example \ref{ex:mc_normalcdf}, consider again the estimation of $\Phi(t)$, the cdf of the standard normal distribution.
As we observed, the basic Monte Carlo method fails to produce accurate estimation for $\Phi(-4.5)$, at least for values of $N$ smaller than $10^6$.
Here, we apply an importance sampling estimator. Let $f$ be the pdf of the standard normal distribution.  We consider a distribution $\ell$ with a support restricted to $(-\infty,-4.5]$ to remove
the unnecessary variation of the basic Monte Carlo estimator due to simulating zero values for $x>-4.5$. A good choice is to take
$\ell(x) = \tilde{\ell}(-x-4.5)$ for $x\in(-\infty,-4.5]$, where $\tilde \ell$ is the pdf of the exponential distribution $\Exp(1)$, i.e.
$$
\ell(x) = \exp(x+4.5),\quad x\leqslant -4.5
$$
See Figure \ref{fig:transfer_Exp}.
\begin{center}
  \includegraphics[scale =0.3]{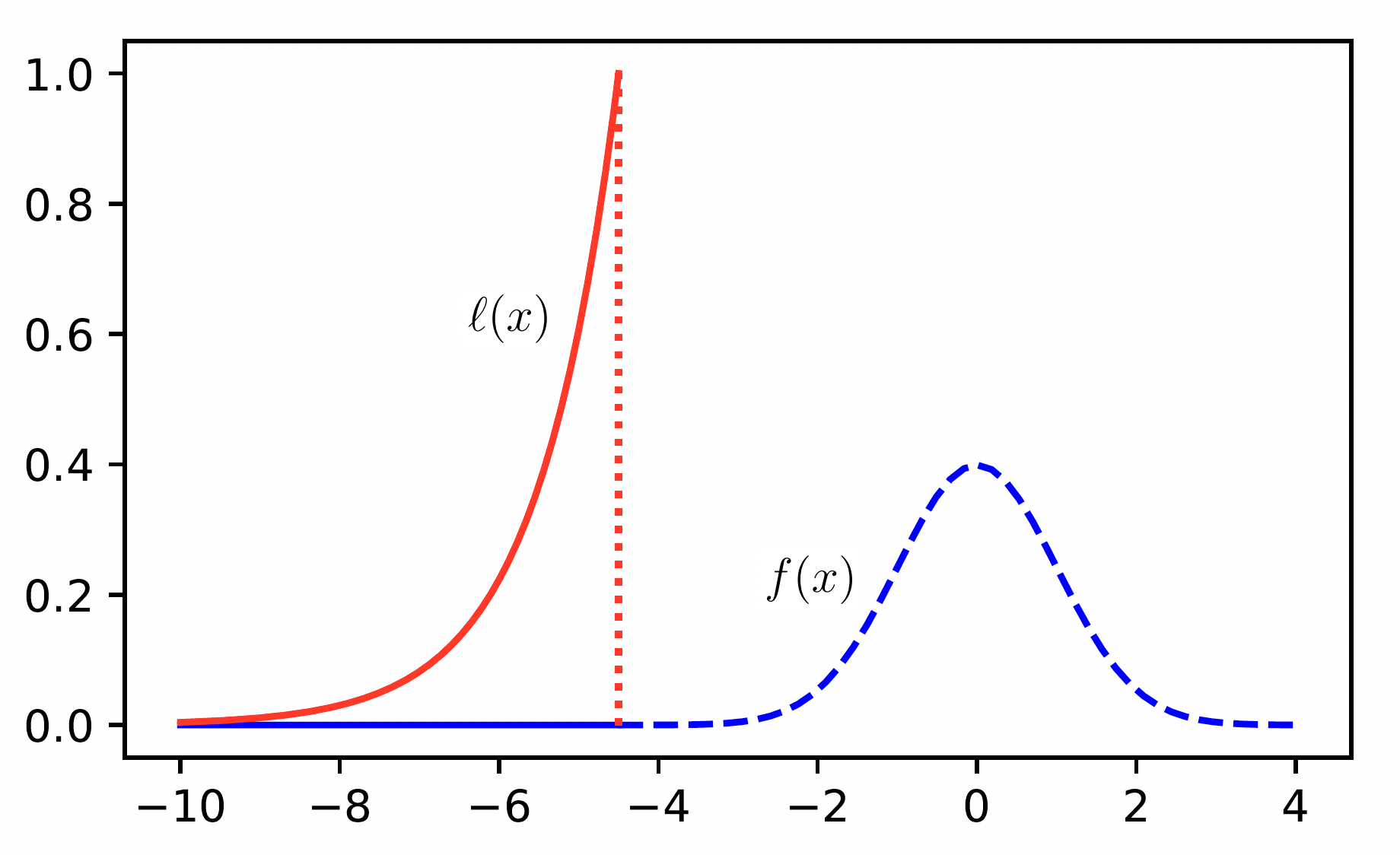}
  \captionof{figure}{The primary pdf $f$ and the importance sampling pdf $\ell$.}\label{fig:transfer_Exp}
\end{center}
The support of $\ell$ is the same as the integral domain in \eqref{cdf-normalPhi} for $t=-4.5$.
If we generate random variables $X_1,\ldots,X_N$ from $\ell$ then $g(X_k)=\mathbb I_{X_k\leqslant -4.5} = 1$ for all $k=1,\ldots,N$ and  the estimation \eqref{mc:importanceEstimate} becomes
\begin{equation}\label{mc:impsamp_est}
\overline \Phi_{\mathrm{IS}}(t)  = \frac{1}{N}\sum_{k=1}^{N} \frac{f(X_k)}{\ell(X_k)}, 
\vsp
\end{equation}
for $t=-4.5$. 
The code and the outputs are given below.
\begin{shaded}
\vspace*{-0.3cm}
\begin{verbatim}
t = -4.5
import numpy as np
from scipy.stats import norm, expon   # normal and exponential dist.
K = 7
Phi, PhiBar = norm.cdf(t), np.zeros(K)
for k in range(K):
    N = 10**(k+1)
    X = -RandExp(1,N) + t
    gX = 1
    fX = norm.pdf(X, loc = 0, scale = 1)
    ellX = expon.pdf(-(X-t), scale = 1)
    PhiBar[k] = 1/N*np.sum(gX*fX/ellX)
Error = abs(PhiBar-Phi)/Phi
print(' Phi(',t,')=',Phi,'\n','PhiBar=',PhiBar,'\n','Error =',Error)
\end{verbatim}
\vspace*{-0.3cm}
\end{shaded}
Outputs:
\begin{shaded}
\vspace*{-0.3cm}
\begin{verbatim}
Phi(-4.5)=3.398e-06
PhiBar=[2.068e-6 4.084e-6 3.332e-6 3.417e-6 3.388e-6 3.396e-6 3.397e-6]
Error =[3.912e-1 2.019e-1 1.946e-2 5.720e-3 2.958e-3 4.943e-4 1.252e-4]
\end{verbatim}
\vspace*{-0.3cm}
\end{shaded}
We observe more accurate estimations compared with the results obtained at the end of Example \ref{ex:mc_normalcdf} using the basic Monte Carlo method.
\end{example}
\vsp
To analyze the results of Examples \ref{ex:mc_normalcdf}  and \ref{ex:is_mc_normalcdf}  we can look at the variances of estimators \eqref{mc:basic_estimate} for basic Monte Carlo and \eqref{mc:importanceEstimate} for Monte Carlo with importance sampling.
Since $X_k$ are assumed to be independent random variables from either $f$ or $\ell$, the variance of $Y$ or $Y_{\mathrm{IS}}$
is the sum of the variances of the individual terms divided by $N$. The variances of individual variables in \eqref{mc:basic_estimate} and \eqref{mc:importanceEstimate} are
\begin{align*}
  &\Var_f[g(X)] = \E_f(g^2(X))-(\E_f[g(X)])^2 = \int g^2(x)f(x)\td x - \left(\int g(x)f(x)\td x\right)^2,\\
  &\Var_\ell\left[g(X)\frac{f(X)}{\ell(X)}\right] = \int g^2(x)\frac{f^2(x)}{\ell(x)} \td x - \left(\int g(x)f(x)\td x\right)^2,
\end{align*}
respectively.
The second terms (squares of expectations) are the same for both variances.  However, depending on the importance sampling function $\ell$, the first term in the second variance can be smaller than the first term in the first variance. To have this, we assume that $f^2(x)/\ell(x)\leqslant f(x)$ or equivalently, $f(x)\leqslant \ell(x)$ for all $x$ in the integration domain. This means that $\ell$ should not have a ``lighter tail'' than $f$. See the graphs of $f$ and $\ell$ in Figure \ref{fig:transfer_Exp}.
\vsp
\begin{workout}
Coming back to Examples \ref{ex:mc_normalcdf} and \ref{ex:is_mc_normalcdf}, show that for $t=-4.5$ we have
\begin{align*}
\Var_f(\overline \Phi)&=\Phi(-4.5)[1-\Phi(-4.5)]/N \doteq \frac{3.3977}{N}\times 10^{-6},\\
\Var_\ell\left(\overline \Phi_{\mathrm{IS}}\right)& = \left[\exp(-4)/(2N\sqrt \pi)\Phi(-4\sqrt{2}) \right]-[\Phi(-4.5)]^2\doteq \frac{3.8373}{N}\times 10^{-11}.
\end{align*}
What do you conclude from this comparison? 
\end{workout}
\vsp 
\begin{workout}
Execute the Python code of Example \ref{ex:is_mc_normalcdf} with $t=0$ instead of $t=-4.5$. Compare the outputs with those of Example \ref{ex:mc_normalcdf}. Report your observation and try to give an analysis.
\end{workout}
\vsp 

To find a proper importance sampling function $\ell$, one may try to minimize
$$
\Var_\ell\left[g(X)\frac{f(X)}{\ell(X)}\right]
$$
over all possible functions $\ell$. We can show that the solution of this minimization problem is
$$
\ell^*(x) = \frac{|g(x)|f(x)}{\ds \int|g(x)|f(x)\td x},
$$
and in particular case $g(x)\geqslant 0$,
$$
\ell^*(x) = \frac{g(x)f(x)}{\E_f[g(X)]}.
\vsp
$$
This sampling function results in a zero variance for its corresponding importance sampling estimator. However, the optimal density \( \ell^* \) is not practical because deriving \( \ell^* \) requires knowing \( \E_f[g(X)] \), which is the quantity we are trying to estimate it! Additionally, in some cases, the explicit form of the performance function \( g(x) \) may not be known in advance. In practice, we try to approximate \( \ell^* \) using sample values \( g(X_1), \ldots, g(X_N) \).  We do not pursue this further and refer you to more advanced texts in Monte Carlo simulation.

%% file: part_StochProc.tex
\section{Stochastic processes}\label{sect:stoch-proc}

A stochastic process is a family of random variables evolving in time. For example, each stochastic solution of the radioactive decay problem given on the right-hand side of Figure \ref{fig:deter-stoch} is a stochastic process. When it comes to compare with the deterministic case, the solution of a deterministic model (e.g. a differential equation) is a {\em function} while the solution of its stochastic counterpart is a {\em stochastic process}. In practice, we can have a set of stochastic processes as solutions, and a possibility is to use the Monte Carlo method to average them and provide an estimate of the expected behavior of the system.  
\vsp 

\begin{definition}
A family of random variables $\{X_t: t\in \mathcal T\}$, all defined on the same probability space, is called a {\em stochastic
process} or random process. If $\mathcal T=\{0,1,2,\ldots\}=:\N_0$ the sequence $X_0,X_1,\ldots$ is called a stochastic process with {\em discrete time parameter}, and if $\mathcal T=[0,\infty)$ the family is called a stochastic process with {\em continuous time parameter}.
 The first random variable $X_0$
is called the {\em initial state} of the process; and the random variable $X_t$ for a $t\in\mathcal T$
is called the {\em state} of the process at time $t$.
\end{definition}
\vsp 

We note that the index family $\mathcal T$ determines the discrete and continuous nature of the stochastic process.
Independent of this, random variables
$X_t$ (states) may have either discrete or continuous probability density functions. 
\vsp 

\begin{example}\label{ex:telephone}
Suppose that a certain business office has five telephone
lines, each one of which might be in use at any given time. At discrete times
$t=0,1,2,\ldots$ minutes the telephone lines are observed and the number of lines that are being used at each time is noted. Let $X_0$
denote the number of lines that are being used at the first time, let $X_1$ denote
the number of lines that are being used at the second time, $1$ minutes later; and in general
$$
X_t = \text{ the number of lines that are being used when they are observed time } t.
$$
Then $X_0,X_1,X_2,\ldots,$ is stochastic process with discrete time parameter.
The state $X_t$ of the process at any discrete time $t$ is the number of lines
being used at that time. Therefore, each state must be an integer between $0$ and $5$.
\end{example}
\vsp 

\begin{example}\label{ex:coins_in_purse1}
A coin purse contains 5 quarters (each worth 25\textcent), 5 dimes (each 10\textcent) and 5 nickels (each 5\textcent). 
Assume that we draw coins one by one and set on a table.  
Let 
$$
X_t = \text{ total value of coins set on the table after } t \text{ draws}.
$$
We see that $\{X_t, \;t=0,1,2,\ldots\}$ is a stochastic process with discrete time variable $t=0,1,2,\ldots$ with initial state $X_0=0$. 

Assume that in the first 6 draws, 3 nickels and 1 quarter and 2 dimes are drawn with
$$
\begin{array}{cccccc}
X_1 = 25, & X_2 = 30,& X_3 = 35, & X_4 = 45,  & X_5 = 50,  & X_6 = 60 .\\
\mbox{{\small quarter}}&\mbox{{\small nickel}}&\mbox{{\small nickel}}&\mbox{{\small dime}}&\mbox{{\small nickel}}&\mbox{{\small dime}}
\end{array}
$$ 
What is the probability of $X_7=65$ given the above information? We are left with 4 quarters, 3 dimes and 2 nickels. We simply see that $\pr(X_7=65\,|\, \text{all information above}) = 2/9 $ and 
$\pr(X_7=80\,|\, \text{all information above}) = 0$.  
This means that we can predict the future of the process (i.e. the state at time $t+1$) based on available information at all previous times $t=0,1,\ldots,t$.  
\end{example}

In a stochastic process with a discrete time parameter, the state of the process
varies in a random manner from time to time. To describe a complete probability
model for a particular process, it is necessary to specify the distribution for the
initial state $X_0$ and also to specify for each $t = 0, 1, \ldots$ the conditional distribution
of the subsequent state $X_{t+1}$ given $X_1, \ldots , X_t$. These conditional distributions are
equivalent to the collection of conditional probabilities of the form
\begin{equation*}
\pr(X_{t+1} = x_{t+1}|X_0 = x_0, X_1 = x_1, \ldots,  X_t=x_t).
\end{equation*}
\vsp 

\subsection{Markov processes}
Markov processes are stochastic processes whose futures are conditionally independent
of their pasts given their present values. On the other words, a stochastic process is Markov, if one can make predictions
for the future of the process based solely on it's present state.
One can say that a Markov process is {\em memoryless}. 
\vsp
\begin{example}\label{ex:coins_in_purse2}
The stochastic process \( X_t \) in Example \ref{ex:coins_in_purse1} is not a Markov process because any prediction about \( X_7 \) would require all information of previous states \( X_6, \ldots, X_1 \). 
Now we define a new process
\[
Y_t = (q_t, d_t, n_t),
\]
where \( q_t \), \( d_t \), and \( n_t \) are the counts of quarters, dimes, and nickels, respectively, on the table at time \( t \). For example, \( Y_0 = (0,0,0) \), \( Y_1 = (1,0,0) \), \( Y_2 = (1,0,1) \), \(\ldots\),  \( Y_6 = (1,2,3) \).
The process \( Y_t \) is Markov because the probability distribution of \( Y_{t+1} \) depends only on the current state \( Y_t \) and not the earlier states. 
\end{example}
\vsp 
A Markov process with a discrete index set, i.e. $\mathcal T = \N_0$
is called a {\em Markov chain}. The states of a Markov chain can be either discrete (countable) or continuous. For a
Markov chain $X = \{X_t: t\in \N_0\}$ with a discrete state space $\Ss$ we have
\begin{equation}\label{markov-property}
  \pr(X_{t+1}=x_{t+1}| X_0 = x_0, X_1 = x_1, \ldots, X_t=x_t) = \pr(X_{t+1}=x_{t+1}|X_t=x_t),
\end{equation}
for all $x_0,\ldots,x_{t+1}\in \Ss$ and $t\in  \N_0$. This property says that the conditional distributions of $X_{t+1}$  given $X_1,\ldots, X_t$ depend only on $X_t$ and not on the earlier states $X_1,\ldots , X_{t-1}$.

If the state space $\Ss$ is finite then the Markov chain is called {\em finite}.
Using the conditional probability rule \eqref{product_rule0}, we have 
\begin{equation}\label{trans01}
\pr(X_0=x_0,X_1=x_1)= \pr(X_0=x_0)\pr(X_1=x_1|X_0=x_0).
\end{equation}
For Markov chains, by applying the product rule \eqref{product_rule} and the Markov chain property \eqref{markov-property}, 
we can prove the following theorem.
\vsp 
\begin{theorem}\label{thm-markov0}
For a finite Markov chain the joint pdf for the first $t$ states is
  \begin{equation}\label{transtt1}
  \begin{split}
    \pr&(X_0=x_0,X_1=x_1,\ldots,X_t=x_t)= \\
   & \pr(X_0=x_0)\pr(X_1=x_1|X_0=x_0) \pr(X_2=x_2|X_1=x_1)\cdots \pr(X_t=x_t|X_{t-1}=x_{t-1}).
    \end{split}
  \end{equation}
\end{theorem}
\vsp
The proof is straightforward and is left as an exercise. 
Theorem \ref{thm-markov0} shows that the joint distribution of a finite Markov chain $X$ can be characterized by the distribution of the initial state $X_0$, i.e., $\pr(X_0=x_0)$, and the {\em one-step transition probabilities}
$$
\pr(X_{t+1}=x_{t+1}|X_{t}=x_{t}), \quad x_t,x_{t+1}\in \Ss.\vsp
$$

In a finite Markov chain we assume that $|\Ss|=m$. It will also
be convenient to name the $m$ states
using the integers $1,2,\ldots,m$ or $0,1,2,\ldots,m-1$. Then for each $t$ and $j$, $X_t= j$ will mean that the chain
is in state $j$ at time $t$. For example, if the states in Example \ref{ex:telephone} are the numbers of
phone lines in use at given times, we have $\Ss=\{0,1,2,3,4,5\}$ and $m = 6$.
\vsp

\begin{definition}
If the transition probabilities $\pr(X_{t+1}=j|X_{t}=i)$ for $i,j\in \Ss$ are independent of the time then the chain is called
{\em time-homogeneous}.
\end{definition}
\vsp 

For a time-homogeneous Markov chain there exist probabilities $p_{ij}$, independent of $t$, such that
\begin{equation}\label{tran-dists}
p_{ij}=\pr(X_{t+1}=j|X_{t}=i), \quad i,j\in \Ss,\quad \forall t\in \N_0.\vsp 
\end{equation}
We can put this probabilities into a matrix $P$ as
$$
P = \begin{bmatrix}
      p_{00} & p_{01} & p_{02} & \cdots \\
      p_{10} & p_{11} & p_{12} & \cdots \\
      p_{20} & p_{21} & p_{22} & \cdots \\
      \vdots & \vdots & \vdots & \ddots
    \end{bmatrix},\vsp
$$
which is called the {\em transition matrix}. If $|\Ss|=m$ then $P$ is an $m\times m$ matrix.
We note that all elements of $P$ are nonnegative and each row sums up to unity. i.e.,
$$
\sum_{j=0}^{\infty}p_{ij}=1,
$$
because in the language of multivariate distributions the transition probabilities $p_{ij}$ are indeed conditional
density functions of $X_{t+1}$
given $X_t$ which can also be denoted by
$
g(j|i) = p_{ij} 
$
for all $t$ and  $i, j$.
\vsp 
\begin{example}\label{ex:weather}
Consider the daily morning weather in a city. Assume there can be
three different states: (1) sunny, (2) cloudy, or (3) rainy.
Observations of the weather forecasting office show that
a sunny day is never
followed by another sunny day. Rainy or cloudy weather is equally probable
after a sunny day. A rainy or cloudy day is followed by $50\%$ probability
by another day with the same weather. If, on the other hand, the weather
is changing from cloudy or rainy weather, the following day will be sunny
only in half of the cases. Based on this observation the transition matrix is
$$
\begin{array}{rccc}
         & \mathrm{sunny}    & \mathrm{cloudy} & \mathrm{rainy} \\
  \mathrm{sunny}  & \big\lceil 0.00 & 0.50 & 0.50\big\rceil \\
  \mathrm{cloudy} & \big|      0.25 & 0.50 & 0.25\big| \\
  \mathrm{rainy}  & \big\lfloor0.25 & 0.25 & 0.50\big\rfloor
\end{array}=:P
\vsp
$$

The transition matrix is independent of time, which means that the transition rules remain unchanged during a period of time; for example during a month. Sometimes it is easier to have a {\em transition graph} instead of the transition matrix. The transition graph of the above example is shown in Figure \ref{fig:weather}.
\begin{center}
\includegraphics[scale=.2]{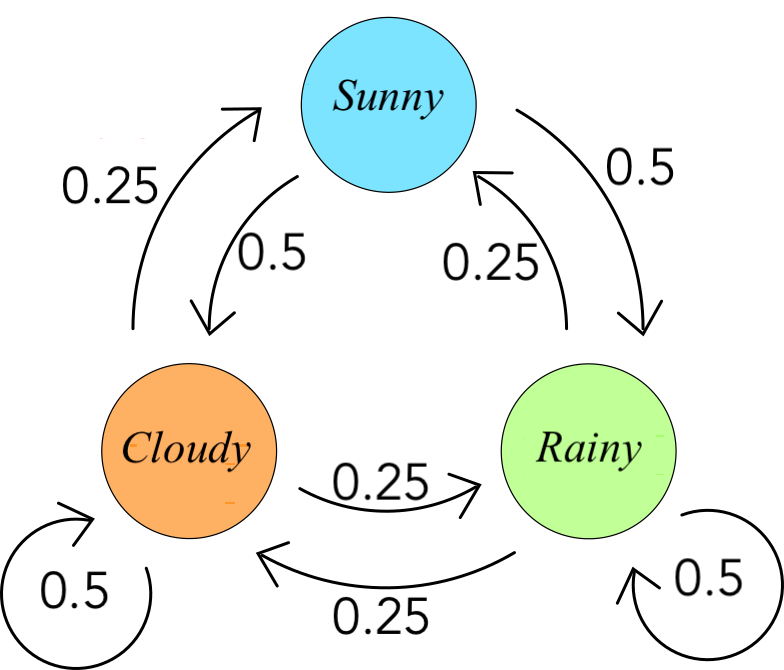}
\captionof{figure}{Transition graph of the daily weather example.}\label{fig:weather}
\end{center}

Looking at the graph, we observe that if, for example, today is cloudy then tomorrow is
sunny with $25\%$ probability, i.e., $\pr(X_{t+1} = sunny\,|\,X_{t} = cloudy) = 0.25$. 
The process is Markov. Why?
\end{example}
\vsp 
%If the distribution of the initial state $X_0$ and the transition matrix $P$ are known then we can completely specify the distribution of $X_t$ at any time $t$.

A vector consisting of nonnegative numbers
that sum to $1$ is called a probability vector. The {\em initial distribution} of the chain is also called the {\em initial probability vector}.

\vsp 
\begin{remark}
Following statistics texts, in this section a vector in $\R^n$ is considered as a $(1\times n)$ array (a row vector). This violates our notation in previous numerical linear algebra lectures where by a vector we meant a column vector.
\end{remark}
\vsp 
For a chain with $m$ possible states $\Ss=\{1,2,\ldots,m\}$, the initial probability vector is
$$
\bpi^{(0)} = \big [\pr(X_0=1),\, \pr(X_0=2), \ldots ,\pr(X_0=m)\big ],
$$
and the distribution (probability vector) of $X$ at time $t$ is denoted by
$$
\bpi^{(t)} = \big [\pr(X_t=1),\, \pr(X_t=2), \ldots ,\pr(X_t=m)\big ].
$$
At $t=1$, the distribution is identified according to \eqref{trans01}:
\begin{align*}
 \pi^{(1)}_j = \pr(X_1=j) &= \sum_{i=1}^{m}\pr(X_1=j,X_0=i) \\
 &=  \sum_{i=1}^{m}\pr(X_1=j|X_0=i)\pr(X_0=i) \\
 &=\sum_{i=1}^{m} p_{ij}\pi^{(0)}_i
\end{align*}
for $j=1,2,\ldots,m$. This can be written in a matrix-vector form as
$$
\bpi^{(1)} = \bpi^{(0)}P.
$$

In general, we can prove by induction that
\begin{equation*}
  \bpi^{(t)} =  \bpi^{(t-1)}P=\bpi^{(t-2)}PP = \cdots = \bpi^{(0)}P^{t}.
\end{equation*}
Note that, if $\bpi^{(t)}$ is a probability vector so is $\bpi^{(t+1)}$, because $0\leqslant p_{ij}\leqslant1$ and rows of $P$ sum up to $1$. We can also show that the $t$-step transition probabilities are
$$
\pr(X_t=j|X_0=i) = p_{ij}^{(t)},\quad i,j,=1,2,\ldots,m
$$
where $p_{ij}^{(t)}$ are entries of matrix $P^t$. The proof for the $2$-step transition is as follows
\begin{align*}
\pr(X_2=j|X_0=i) &= \sum_{k=1}^{m} \pr(X_1=k,X_2=j|X_0=i)\\
& = \sum_{k=1}^{m} \pr(X_1 = k|X_0=i) \pr(X_2 = j |X_1 = k, X_0= i)\\
& = \sum_{k=1}^{m} \pr(X_1 = k|X_0=i) \pr(X_2 = j |X_1 = k)\qquad (\mbox{Markov chain property})      \\
& = \sum_{k=1}^{m} p_{ik}p_{kj} = p_{ij}^{(2)}.
\end{align*}
The general case can be proved similarly.
\vsp 
\begin{example}\label{ex:weather2}
Consider the transition matrix in Example \ref{ex:weather}.
If the weather is cloudy in the day 0 then $\bpi^{(0)}=[0,1,0]$, and the probability vector $\bpi^{(1)}$ is computed by
$$
\bpi^{(1)} = \bpi^{(0)}P =
\begin{bmatrix}
  0 &
  1 &
  0
\end{bmatrix}
\begin{bmatrix}
                 0 &0.50 & 0.50 \\
                 0.25 & 0.5 &0.25  \\
                 0.25 & 0.25 & 0.5
               \end{bmatrix}
 =
\begin{bmatrix}
  0.25 &
  0.5 &
  0.25
\end{bmatrix}.
$$

This means that in the first day the weather is sunny with probability $0.25$, cloudy with probability $0.5$ and rainy with probability $0.25$.
For the second day we compute $\bpi^{(2)}$:
$$
\bpi^{(2)} = \bpi^{(1)}P =
\begin{bmatrix}
  0.25 &
  0.5 &
  0.25
\end{bmatrix}
\begin{bmatrix}
                0 &0.50 & 0.50 \\
                 0.25 & 0.5 &0.25  \\
                 0.25 & 0.25 & 0.5
               \end{bmatrix}
 \doteq
\begin{bmatrix}
  0.19 &
  0.44 &
  0.37
\end{bmatrix}.
$$
We can simply predict the weather for days after.
\end{example}
\vsp 

\subsection*{Stationary distribution}
In Example \ref{ex:weather2} (daily weather) let us continue to compute distribution vectors $\bpi^{(t)}$ for higher values of $t$,
with $\bpi^{(t+1)}=\bpi^{(t)}P$ or equivalently $\bpi^{(t+1)}=\bpi^{(1)}P^{t}$. The results with initial vector $\bpi^{(0)}=[0,1,0]$ are (by rounding the final numbers to $4$ decimal digits)
$$
\bpi^{(1)}=
\begin{bmatrix}
  0.25 \\
  0.50 \\
  0.25
\end{bmatrix}^T,\;
\bpi^{(2)}=
\begin{bmatrix}
    0.1875\\
    0.4375\\
    0.3750
\end{bmatrix}^T,\;
\bpi^{(3)}=
\begin{bmatrix}
    0.2031\\
    0.4063\\
    0.3906\\
\end{bmatrix}^T,\;
\ldots, \;
\bpi^{(6)}=
\begin{bmatrix}
    0.2000 \\
    0.4001 \\
    0.3999
\end{bmatrix}^T,\;
\bpi^{(7)}=
\begin{bmatrix}
    0.2000\\
    0.4000\\
    0.4000
\end{bmatrix}^T
$$
The matrix $P^7$ is
$$
P^7\doteq
\begin{bmatrix}
    0.2000  &  0.4000  &  0.4000\\
    0.2000  &  0.4000  &  0.4000\\
    0.2000  &  0.4000  &  0.4000
\end{bmatrix}.\vsp
$$
The same result with more or less number of iterations will obtained for other choices of initial distribution vector $\bpi^{(0)}$. The distribution of the Markov chain approaches the {\em stationary distribution} $\bpi=[0.2,0.4,0.4]$.
In fact, all rows of $P^t$ approaches this stationary distribution:
\begin{equation}\label{limergodic}
\lim_{t\to\infty}p_{ij}^{(t)} = \pi_j, \quad \forall i\in \Ss 
\end{equation}
with $\pi_j>0$. 
This, on the other hand, means that if we set the initial distribution to
$$
\bpi^{(0)} = [0.2,0.4,0.4],
$$
then we can show that $\bpi^{(1)}=\bpi^{(0)}P=\bpi^{(0)}=[0.2,0.4,0.4]$,
which means that $\bpi^{(0)}$ is also the distribution after one transition. Hence, it will
 be the distribution after two or more transitions.
\vsp 
\begin{definition}
The Markov chain $ X_0,X_1,\ldots$ with transition matrix $P=(p_{ij})$ is called ergodic if the limits
\eqref{limergodic}, i.e., $\pi_j$
\begin{enumerate}
\item exist for all $ j\in \Ss$,
\item are positive ($\pi_j>0$) and independent of $ i\in \Ss$
\item form a probability vector $ {\boldsymbol{\pi}}=(\pi_1,\ldots,\pi_m)$, i.e., $ \ds\sum_{j\in \Ss}\pi_j=1$.
\end{enumerate}
\end{definition}
\vsp 
If the limits \eqref{limergodic} exist then $\ds \lim_{t\to\infty}\bpi^{(t+1)} = \ds \lim_{t\to\infty}\bpi^{(t)}P$ which means 
$\bpi P = \bpi$.
\vsp 
\begin{definition}[Stationary Distribution]
Let $P$ be the transition matrix for a Markov chain. A probability
vector $\bpi$ that satisfies
\begin{equation}\label{station-dist}
\bpi P = \bpi
\end{equation}
is called a {\em steady-state distribution} for the Markov
chain.

\end{definition}
\vsp 

From \eqref{station-dist} we observe that for a finite Markov chain if the stationary distribution
$\bpi$ exists then it is the
eigenvector of $P^T$ corresponds to eigenvalue $1$, because 
\begin{equation}\label{Pteigval}
P^T \bpi^T = \bpi^T.
\end{equation}
On the other hand, since the rows of $P$ sum up to unity we have 
$$
P\e^T = \e^T
$$ 
where $\e = [1,1,\ldots,1]\in\R^m$ which shows that $1$ is an eigenvalue of
$P$ and thus an eigenvalue of $P^T$ because eigenvalues of $P$ and $P^T$ are the same.
This eigenvalue should correspond to at least one eigenvector for $P^T$ that is $\bpi$ from \eqref{Pteigval}.
However,
the uniqueness of this eigenvector is still not guaranteed. 

%To ensure uniqueness, the $P$ must be {\bf irreducible}.

Let us characterize the relations between states in the following way: 
For arbitrary but fixed states $ i,j\in \Ss$ we say that the state $ j$ is accessible from state $ i$ if $ p_{ij}^{(t)}>0$ for some $t\geqslant0$, we say that $i$ is accessible from (or leads to) state $j$ and write $i\rightarrow j$.
We say that $i$ and $j$ {\em communicate} if $i\rightarrow j$ and $j\rightarrow i$, and write $i\leftrightarrow j$. Using the
relation `$i\leftrightarrow j$', we can divide the state space $\Ss$ into equivalence classes such that all the states in
an equivalence class communicate with each other but not with any state outside
that class. If there is only one class (i.e., $\Ss$ itself), the Markov chain is said to
be {\em irreducible}. Besides irreducibility we need a second property of the transition probabilities, namely  aperiodicity.
The period $d_i$ of the state $ i\in \Ss$ is given by 
$$ 
d_i={\rm gcd}\{t\ge 1:p_{ii}^{(t)}>0\}
$$ 
where ``gcd'' denotes the greatest common divisor. We define $d_i=\infty$ if $ p_{ii}^{(t)}=0$ for all $ t\ge 1$.
A state $ i\in \Ss$ is said to be aperiodic if $d_i=1$.
The Markov chain $ \{X_t\}$ and its transition matrix $P=(p_{ij})$ are called {\em aperiodic} if all states of $ \{X_t\}$ are aperiodic.
We can show that the periods $ d_i$ and $ d_j$ coincide if the states $ i,j$ belong to the same equivalence class of communicating states.
Thus, if the Markov chain $ \{X_t\}$ is irreducible then all its states have the same period. 
Here we give the statements of two fundamental theorems without proofs. 
For more details refer to \cite{Rubinnstein-Kroese:2017}. 
\vsp 

\begin{theorem}
The Markov chain $ X_0,X_1,\ldots$ is ergodic if and only if it is irreducible and aperiodic. 
\end{theorem}
\vsp 

\begin{theorem}
For an irreducible and aperiodic Markov chain $ X_0,X_1,\ldots$ (ergodic Markov chain) with transition matrix $P$,
the stationary distribution $\bpi$ is {\em uniquely} determined by solving the eigenvalue problem \eqref{station-dist}. The eigenvalue $\lambda_1=1$ is the  dominant
eigenvalue and $|\lambda_k|< |\lambda_1|=1$ for $\lambda = 2,3,\ldots,m$. 
\end{theorem}
\vsp

\subsection{Random walk on the integers}\label{sect:randomwalk}
Assume that your initial state is $X_0 = 0$ and at each time $t$ you walk by steplength $1$ either to the right or to the left of your current position $X_t$ with probabilities $p$ and $q$, respectively, where 
$p\in(0,1)$ is a real number and $q=1-p$. The process $X_t$ is called a 
{simple random walk} which has the state space $\Ss = \mathbb Z$ (set of integer numbers) and a transition graph shown in Figure \ref{fig:randomwalk}.

\begin{center}
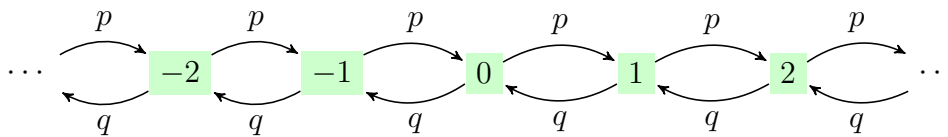

\begin{tikzpicture}[->,>=stealth',shorten >=1pt,auto,node distance=2cm, 
                    semithick]
\node[fill=green!20,draw=none,text=black] (A)            {$-2$};
\node[fill=green!20,draw=none,text=black,right of=A] (B) {$-1$};
\node[fill=green!20,draw=none,text=black,right of=B] (C) {$0$};
\node[fill=green!20,draw=none,text=black,right of=C] (D) {$1$};
\node[fill=green!20,draw=none,text=black,right of=D] (E) {$2$};
\node[fill=green!0, draw=none,text=black, left of=A] (Z) {$\cdots$};
\node[fill=green!0, draw=none,text=black,right of=E] (F) {$\cdots$};

% Arrows to the left
\path (Z) edge [bend left] node[above] {$p$} (A);
\path (A) edge [bend left] node[above] {$p$} (B);
\path (B) edge [bend left] node[above] {$p$} (C);
\path (C) edge [bend left] node[above] {$p$} (D);
\path (D) edge [bend left] node[above] {$p$} (E);
\path (E) edge [bend left] node[above] {$p$} (F);

% Arrows to the right
\path (A) edge [bend left] node[below] {$q$} (Z);
\path (B) edge [bend left] node[below] {$q$} (A);
\path (C) edge [bend left] node[below] {$q$} (B);
\path (D) edge [bend left] node[below] {$q$} (C);
\path (E) edge [bend left] node[below] {$q$} (D);
\path (F) edge [bend left] node[below] {$q$} (E);
\end{tikzpicture}
 \captionof{figure}{Transition graph for random walk on integers}\label{fig:randomwalk}
\end{center}
The random walk $X_t$ can be characterized as 
$$
X_{t+1} = X_t + \tilde B_t, \quad t = 0,1,2,\ldots,
$$
where
$\tilde B_t$ has a Bernoulli distribution with state $\{-1,1\}$ and probabilities $q$ and $p$. Indeed $\tilde B_t = 2B_t-1$ where $B\sim \Ber (p)$, the standard Bernoulli distribution. The above relation shows that $X_{t+1}$ depends solely on $X_t$, i.e. the random walk is a Markov process. It is also a Markov chain because the index set for $t$ is discrete.  
We also observe that 
$$
X_t = \sum_{j=0}^{t-1} \tilde B_j,
$$
which means that the distribution of $X_t$ is binomial\footnote{The sum of Bernoulli distributions is a binomial distribution.}.
Since $\E(\tilde B_j) = 2p-1$ and $\Var(\tilde B_j) = 4p(1-p)$, we simply have
$$
\E(X_t) = t(2p-1), \quad \Var(X_t) = 4tp(1-p). 
$$
For special case $p=q=1/2$ we have $\E(X_t) = 0$ and $\Var(X_t) = t$.

The transition matrix of the random walk process has the form
$$
P = \begin{bmatrix}
      \ddots &\vdots & \vdots & \vdots & \vdots & \iddots\\
      \cdots & 0 & p & 0 & 0 & \cdots \\
      \cdots & q & 0 & p & 0 & \cdots \\
      \cdots &0  & q & 0 & p &\cdots \\
      \cdots &0  & 0 & q & 0 &\cdots \\
      \iddots &\vdots &\vdots &\vdots & \vdots & \ddots
    \end{bmatrix}.
$$
Since $X$ starts at $0$, i.e.,
$\pr(X_0=0)=1$ and $\pr(X_0=j)=0$ for all $j\in\mathbb Z\setminus\{0\}$, we have
\begin{align*}
\bpi^{(0)} = [\cdots,0,0,1,0,0,\cdots].
\end{align*}
Then $\pi^{(1)}=\pi^{(0)}P$ has the form
\begin{align*}
\bpi^{(1)} = [\cdots,0,q,0,p,0,\cdots],
\end{align*}
which means that after one step the chain moves to right with probability $p$ and to left with probability $q$.
\vsp 

\subsection{Gaussian processes}
Gaussian processes are generalizations of multivariate normal distributions. 
If you are not familiar with multivariate (normal) distributions, see sections \ref{sect:jointdist}-\ref{sect:multi-normal} in the Appendix. 
\vsp 
\begin{definition}
The stochastic process $ \{X_t,\,t\in\mathcal T\}$ is called a {\em Gaussian process} if all its finite dimensional distributions are normal (Gaussian). In the other words, $\{X_t,\,t\in\mathcal T\}$ is a Gaussian process if for any choice of $n$ and $t_1,t_2,\ldots,t_n\in\mathcal T$ we have
\begin{equation*}
  (X_{t_1},X_{t_2},\ldots,X_{t_n})\sim \Nor(\mu,\Sigma)
\end{equation*}
for some expectation vector $\mu$ and covariance matrix $\Sigma$ both depend on the choice of $t_1,\ldots,t_n$.
\end{definition}
\vsp

Equivalently,
$\{X_t,\,t\in\mathcal T\}$ is a Gaussian process if any linear combination
$$
\sum_{k=1}^{n}c_k X_{t_k}
$$
has a normal distribution. A Gaussian process is fully determined by its expectation function $\mu(t) = \E(X_t)$ for $t\in\mathcal T$ and
covariance function $\Sigma(s,t)= \Cov(X_s,X_t)$ for $s,t\in\mathcal T$.

An important Gaussian process is the {\em Wiener process} or the standard {\em Brownian motion}, which can also be considered as a continuous version of the random walk on the integers.
\begin{definition}
A Gaussian process $\{W_t,\,t\in\mathcal T\}$ with $\mu(t)=0$ for all $t\in\mathcal T$ and $\Sigma(s,t)=s$ for all
$0\leqslant s\leqslant t$ is called a {Wiener process}.
\end{definition}
We can prove that, a Wiener process is a Markov process with a continuous sample path that is nowhere differentiable. Moreover, in the Wiener process the increments $W_t-W_s$ on intervals $[s,t]$ are independent and normally distributed. More precisely, from the definition and for two choices $s,t\in \mathcal T$ with $0\leqslant s\leqslant t$ we have 
$$
\begin{bmatrix}
W_s\\ W_t
\end{bmatrix}\sim \Nor\left(\begin{bmatrix}
0\\ 0
\end{bmatrix},\begin{bmatrix} s&s\\s & t \end{bmatrix}\right).
$$
The Cholesky factorization of the covariance matrix is 
$$
\Sigma = \begin{bmatrix} s&s\\s & t \end{bmatrix} = \begin{bmatrix} \sqrt s&0\\ \sqrt s & \sqrt{t-s} \end{bmatrix} \begin{bmatrix} \sqrt s&\sqrt s\\ 0 & \sqrt{t-s} \end{bmatrix}=:BB^T
$$
which results in 
$$
\begin{bmatrix}
W_s\\ W_t
\end{bmatrix}= 
\begin{bmatrix}
0\\ 0
\end{bmatrix} + 
\begin{bmatrix} \sqrt s&0\\ \sqrt s & \sqrt{t-s} \end{bmatrix} 
\begin{bmatrix}
Z_1\\ Z_2
\end{bmatrix} 
\vsp
$$
where $Z_1$ and $Z_2$ are two univariate standard normal distributions\footnote{If $W\sim \Nor(\mu,\Sigma)$ then $W = \mu + B Z$ where $Z\sim \Nor({0},I)$ and $B$ is the Cholesky factor of $\Sigma$.}.
Solving the system gives $W_s = \sqrt s Z_1$ and $W_t = W_s + \sqrt{t-s} Z_2$. The later shows 
\begin{equation}\label{wiener_property}
W_t-W_s \sim \Nor(0,t-s), \quad \mbox{for all }t\geqslant s\geqslant0.
\end{equation}
This proves that in the Wiener process the increments $W_t-W_s$ on intervals $[s,t]$ are independent and normally distributed.
\vsp

%% file: part_StochProcGen.tex
\section{Stochastic process generation}\label{sect:stoch_proc_gen}

In this section, we provide a brief overview of algorithms for generating a few standard stochastic processes. However, stochastic processes are not limited to the standard examples presented here. For any specific stochastic model, users can use standard methods for sampling random points and random processes to generate a specific stochastic process that is a solution to their model. In section \ref{sect:ssa} we will see an example. 
\vsp
\subsection{Generating Markov chains}
Assume that $X = \{X_0, X_1,\ldots,X_n\}$ is the first $n+1$ random variables of a finite Markov chain with initial distribution $\bpi^{(0)}$, state space $\Ss=\{1,2,\ldots,m\}$, and transition matrix $P$.
We follow the simulation process given at the beginning of section \ref{sect:randvectorgen} for dependent variables by the use of conditional distributions.
We first generate $X_0$ from distribution $\bpi^{(0)}$. If $X_0=i$ is generated, we then generate $X_1$ from the conditional distribution of $X_1$ given $X_0=i$. In the other words, we generate $X_1$ from the $i$-th row of $P$. Let $X_1=j$ be generated. From here on, we use the Markov property, so we generate $X_2$ from the conditional distribution of $X_2$ given $X_1=j$, i.e. we  generate $X_2$ from the $j$-th row of $P$. This process is continued until $X_n$ is generated.

In the Python function below, the input variables \verb+InitDist+, \verb+TransMat+ and \verb+ChainLen+ play the roles of
$\bpi^{(0)}$, $P$ and $n$, respectively. The outputs are the integer vector (state vector) \verb+X+ of size $n$, and the stationary distribution \verb+StationDist+ of length $m$. The stationary distribution is approximated by
$m$ iterations of the {\em power method} for the dominant eigenvector of $P$.
Note that, the \verb+RandDisct+ function (the discrete random variable generator) is called to generate random numbers from the rows of $P$.

\begin{shaded}
\vspace*{-0.3cm}
\begin{verbatim}
def MarkovChainGen(InitDist, TransMat, ChainLen):
    # This function generates a Markov chain of length 'ChainLen' from 
    #   Transition matrix 'TransMat' with initial distribution 'InitDist'
    m = len(InitDist)
    states = range(m)
    X = np.zeros(ChainLen)
    p = StationDist = InitDist
    for j in range(ChainLen):
        i = RandDisct(states,p,1)[0]
        X[j] = i
        p = TransMat[i,:]
        StationDist = TransMat@StationDist
    return X, StationDist
\end{verbatim}
\vspace*{-0.3cm}
\end{shaded}
\vsp
\begin{example}
Consider again the weather forecasting model in Example \ref{ex:weather}. 
We start with the initial state of sunny on the day 0. We want to estimate the probability of a rainy day on the fifth day. 
We can generate multiple weather sequences of length $5$ (Markov chains of length $5$) and use Monte Carlo to  estimate the probability of a rainy day on the fifth day. The code is given below.

\begin{shaded}
\vspace*{-0.3cm}
\begin{verbatim}
# Transition matrix
P = np.array([[0.00, 0.50, 0.50],  # From sunny to (sunny, cloudy, rainy)
              [0.25, 0.50, 0.25],  # From cloudy to (sunny, cloudy, rainy)
              [0.25, 0.25, 0.50]]) # From rainy to (sunny, cloudy, rainy)
# Map state names to indices for the matrix
StateMap = {"sunny": 0, "cloudy": 1, "rainy": 2}
InitDist = [1,0,0]  # [sunny, cloudy, rainy]
RainyCount = 0
N = 1000
for i in range(N):    # Run N simulations
    X, S = MarkovChainGen([1,0,0], P, 5)
    if X[-1] == StateMap["rainy"]: # Check if the fifth day ended in rainy
        RainyCount += 1

PrRainy = RainyCount / N  # Estimate probability
print('Estimated probability of a rainy fifth day = ', PrRainy)

\end{verbatim}
\vspace*{-0.3cm}
\end{shaded}

An execution gives the answer $\texttt{0.402}$ which shows that by 
about $40\%$ probability the fifth day is rainy. 
\end{example}
\vsp 

\subsection{Random walk on the integers}
Random walk on the integer is Markov chain as it was described in section \ref{sect:randomwalk}. Thus, the \verb+MarkovChainGen+ function can be used to generated this process. However, a simpler algorithm without forming the transition matrix $P$ is based on the relation $X_{t+1}=X_t+\tilde B_t$.
\begin{shaded}
\vspace*{-0.3cm}
\begin{verbatim}
def RandWalkGen(p, X0, t):
    B = np.append(X0,2*RandBer(p,t)-1)
    X = np.cumsum(B)
    return X
\end{verbatim}
\vspace*{-0.3cm}
\end{shaded}

A typical sample path for the case $p=q=1/2$, $X_0=0$ and $t=100$ is plotted in Figure \ref{fig:randomwalk_path}. The horizontal axis represents the time index $t=0,1,2,\ldots,100$ and the vertical axis the states $j\in \mathbb Z$.

\begin{center}
\includegraphics[scale=0.63]{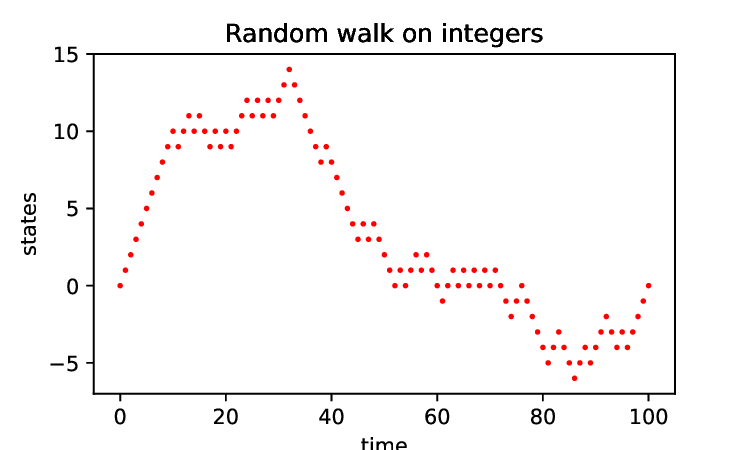}
\captionof{figure}{A random walk path on the integers with $p=1/2$.}\label{fig:randomwalk_path}
\end{center}
\vsp 
\begin{example}[gambler's ruin]
Suppose a gambler starts with a certain amount of money, say $\$K$, and in each round, either wins or loses a fixed amount, say $\$1$, with  probabilities $p$ and $1-p$, respectively. The game continues until the gambler either reaches a target amount $\$T$ or loses all the money (being {\em ruined}). 
The target amount \(T\) is an {\em absorbing state} on the right side and the amount $0$ is an {absorbing state} on the left side. 

This scenario is a random walk. 
The goal is to determine the probability that the gambler reaches the target amount \(T\) before being ruined. 
We generate many sample paths and use Monte Carlo to estimate this probability. 
The simulation should track whether the gambler ends up with \(X_t = 0\) or \(X_t = T\) in multiple paths. By averaging the outcomes, we estimate the probability of ruin. The code is given below  for some values of parameters. 

\begin{shaded}
\vspace*{-0.3cm}
\begin{verbatim}
T = 100   # target amount
K = 30    # initial money
p = 0.5   # probability of winning each round
N = 1000  # number of MC simulations
RuinCount = 0
for _ in range(N):
    money = K
\end{verbatim}
%\vspace*{-0.3cm}
\end{shaded}

\begin{shaded}
%\vspace*{-0.3cm}
\begin{verbatim}
    while money > 0 and money < T:
        B = 2*RandBer(p,1)-1
        money += B
    if money == 0:
        RuinCount += 1
PrRuin = RuinCount / N  # Estimated probability of ruin
print('Estimated Probability of Ruin = ', PrRuin)
\end{verbatim}
\vspace*{-0.3cm}
\end{shaded}

An execution shows the probability $\texttt{0.703}$. This problem has indeed an exact solution 
\[
\pr(\text{ruin}) = 
\begin{cases}
    \frac{1 - \left(\frac{q}{p}\right)^k}{1 - \left(\frac{q}{p}\right)^T} & \text{if } p \neq q, \\
    \frac{T - k}{T} & \text{if } p = q = \frac{1}{2},
\end{cases}
\]
which can be used to assess the accuracy of the above Monte Carlo algorithm. 
\end{example}
\vsp
A random walk on integers can model several real-world scenarios.
In a financial modeling, a random walk can be used to simulate investment portfolios, where each step represents a gain or loss in asset value.
In population dynamics, a random walk can model population survival, where each step could represent birth or death events.
In a biological model, a random walk can describe gene mutations, where each step indicates a genetic drift towards different traits or alleles.

The random walk can be extended to $d$ dimensions by replacing the Bernoulli distribution with a general discrete distribution for updates.
We indeed have 
$$
X_{t+1} = X_t + D_t
$$
where $D_t\sim \mathcal D\mathcal D([x_1,\ldots,x_{2d}],[p_1,\ldots,p_{2d}])$. As an example in 2 dimensions, we have $x_1=[-1,0]$, $x_2=[1,0]$, $x_3=[0,1]$, and $x_4=[0,-1]$, which means that at each time step the process moves  to left, right, up or down with probabilities $p_1$, $p_2$, $p_3$ and $p_4$ respectively. 

\begin{center}
\includegraphics[scale=.4]{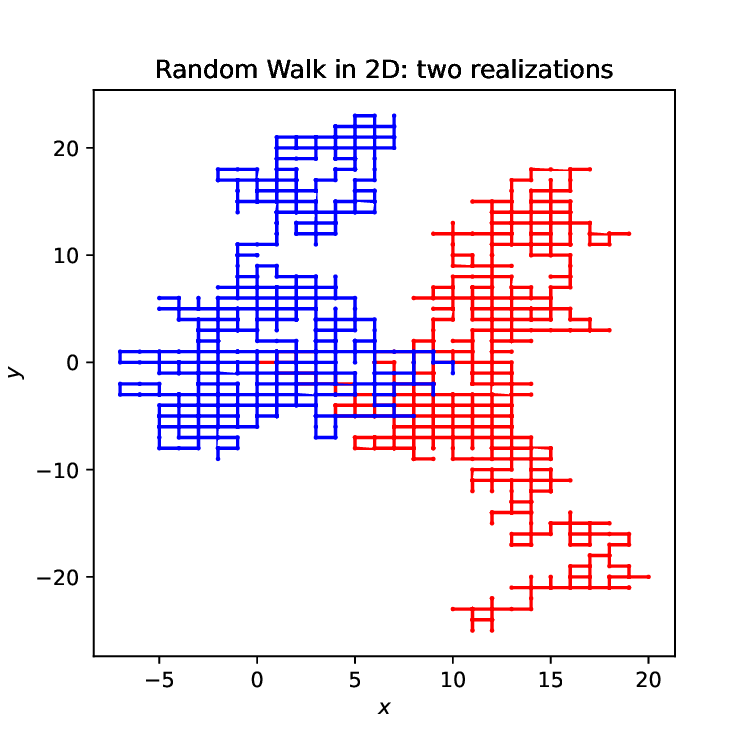}\includegraphics[scale=.42]{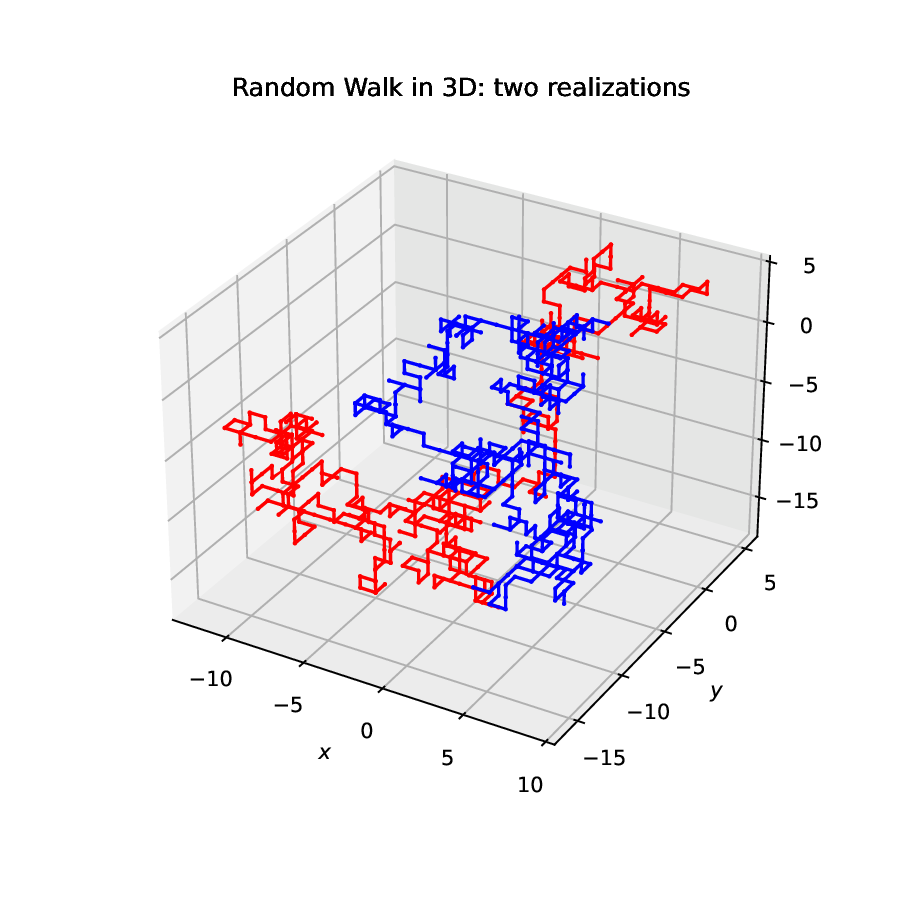}
\captionof{figure}{Random walk realizations in 2 and 3 dimensions.}\label{fig:randomwalk2D}
\end{center}

Two random walk realizations for $X_0 = [0,0]$ and $p_k=1/4$ are shown on the left side of Figure \ref{fig:randomwalk2D}.
The same extension applies in 3 dimensions with 6 possible different directions for a new update. A couple of realizations are shown on the right side of Figure \ref{fig:randomwalk2D} for $X_0=[0,0,0]$ and $p_k = 1/6$. 
In both 2 and 3 dimensions, the code is left as an exercise to the reader. 
\vsp

\subsection{Generating Gaussian processes}\label{sect:gauss_process_gen}
In a Gaussian process $\{X_t:\,t\in\mathcal T\}$, each finite dimensional vector $(X_{t_1},\ldots,X_{t_n})$ has a multivariate normal distribution. This means that any multivariate normal sampler can be used to generate realizations of a Gaussian process at prescribed times $t_1,\ldots,t_n$ provided that the mean vector $\mu = (\mu(t_1),\ldots,\mu(t_n))$ and the covariance matrix $\Sigma=(\Sigma(t_k,t_j))$ for $k,j=1,2\ldots,n$ are given.

For a Wiener process (or Brownian motion) $\{W_t:\,t\in\mathcal T\}$, as an special case, the algorithm can be simplified. For this process we have $\mu(t)=0$ and
$\Sigma(s,t)=s$ for $0\leqslant s\leqslant t$. For given times $t_1,\ldots,t_n$ with $0<t_1<t_2<\cdots <t_n$, the covariance matrix for variable $(W_{t_1},\ldots,W_{t_n})$ has the form 
$$
\Sigma =
\begin{bmatrix}
  t_1 & t_1 & t_1 & \cdots & t_1 \\
  t_1 & t_2 & t_2 & \cdots & t_2 \\
  t_1 & t_2 & t_3 & \cdots & t_3 \\
  \vdots & \vdots & \vdots & \ddots & \vdots \\
  t_1 & t_2 & t_3 & \cdots & t_n
\end{bmatrix}.
%\begin{bmatrix}
%  \sqrt{t_1} & \sqrt{t_1} & \sqrt{t_1} & \cdots & \sqrt{t_1} \\
%  0 & \sqrt{t_2-t_1} & \sqrt{t_2-t_1} & \cdots & \sqrt{t_2-t_1} \\
%  0 & 0 & \sqrt{t_3-t_2} & \cdots & \sqrt{t_3-t_2} \\
%  \vdots & \vdots & \vdots & \ddots & \vdots \\
%  0 & 0 & 0 & \cdots & \sqrt{t_n-t_{n-1}}
%\end{bmatrix}=:BB^T \vsp
$$
The Cholesky factor for $\Sigma$ is 
$$
B =
\begin{bmatrix}
  \sqrt{t_1} & 0 & 0 & \cdots & 0 \\
  \sqrt{t_1} & \sqrt{t_2-t_1} & 0 & \cdots & 0 \\
  \sqrt{t_1} & \sqrt{t_2-t_1} & \sqrt{t_3-t_2} & \cdots & 0 \\
  \vdots & \vdots & \vdots & \ddots & \vdots \\
  \sqrt{t_1} & \sqrt{t_2-t_1} & \sqrt{t_3-t_2} & \cdots & \sqrt{t_n-t_{n-1}}
\end{bmatrix}.\vsp
$$
Then $W = \mu + BZ = { 0} + BZ$ where $Z$ is a vector of iid random variables with distributions $\Nor(0,1)$, i.e.,
$$
W =
\begin{bmatrix}
  W_{1} \\
  W_2 \\
  \vdots \\
  W_n
\end{bmatrix} = B
\begin{bmatrix}
  Z_{1} \\
  Z_2 \\
  \vdots \\
  Z_n
\end{bmatrix}
=
\begin{bmatrix}
  W_0 + \sqrt{t_1-t_0} Z_1 \\
  W_1 + \sqrt{t_2-t_1}Z_2 \\
  \vdots \\
  W_{n-1} + \sqrt{t_n-t_{n-1}}Z_n
\end{bmatrix}\vsp
$$
for $t_0=0$ and $W_0=0$. This means that 
$$
W_{k+1} = W_k + \sqrt{t_{k+1}-t_k}\, Z_{k+1}, \quad Z_{k+1} \sim \Nor(0,1), \quad W_0 = 0.
$$
The Python function is given below. 
\begin{shaded}
\vspace*{-0.3cm}
\begin{verbatim}
def BrownianMotionGen(t_vec, dim):
    # This functions generates a dim-dimensional Brownian motion on 
    #    time samples t_vec = [t_1, t_2, ..., t_n]
    n = len(t_vec)
    W = np.zeros([dim,n])
    for k in range(n-1):
        Z = np.random.normal(0,1,dim)
        W[:,k+1] = W[:,k] + np.sqrt(t_vec[k+1]-t_vec[k])*Z
    return W
\end{verbatim}
\vspace*{-0.3cm}
\end{shaded}
However, the function works also for higher dimensions. 
In higher dimensions each step of the process updates by a vector whose components are drawn from a multivariate normal distribution. 
We note that in higher dimensions the term {\em Brownian motion} is commonly used while in the 1-dimensional case both terms {\em Wiener process} and {\em Brownian motion} are used.

In Figure \ref{fig:wiener_path} two sample paths of the 1D, 2D and 3D Brownian motions in time interval $[0,1]$ are shown at times $t_k = k\Delta t$ for $\Delta t=10^{-3}$.

\begin{center}
\includegraphics[scale=0.55]{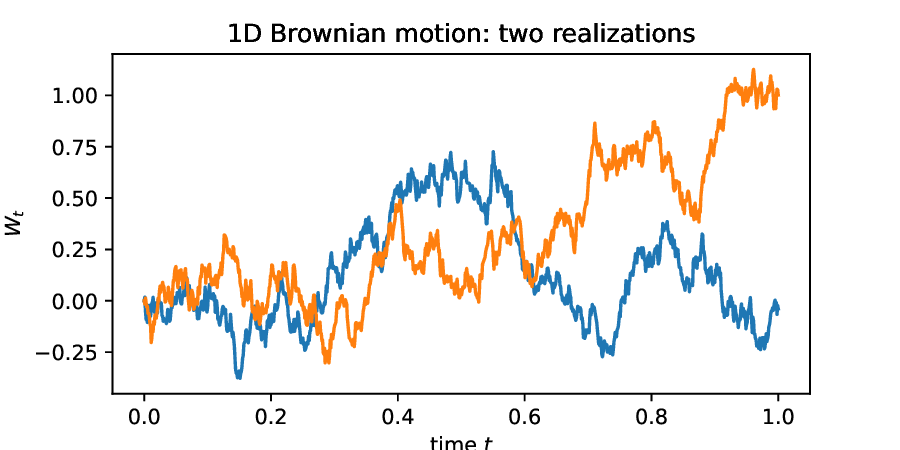}\\
\includegraphics[scale=0.5]{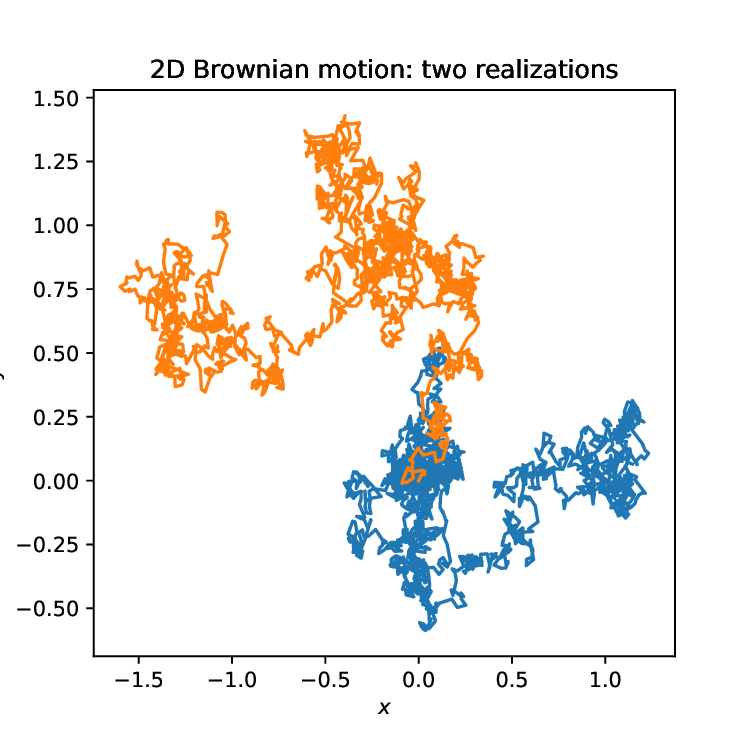}\includegraphics[scale=0.5]{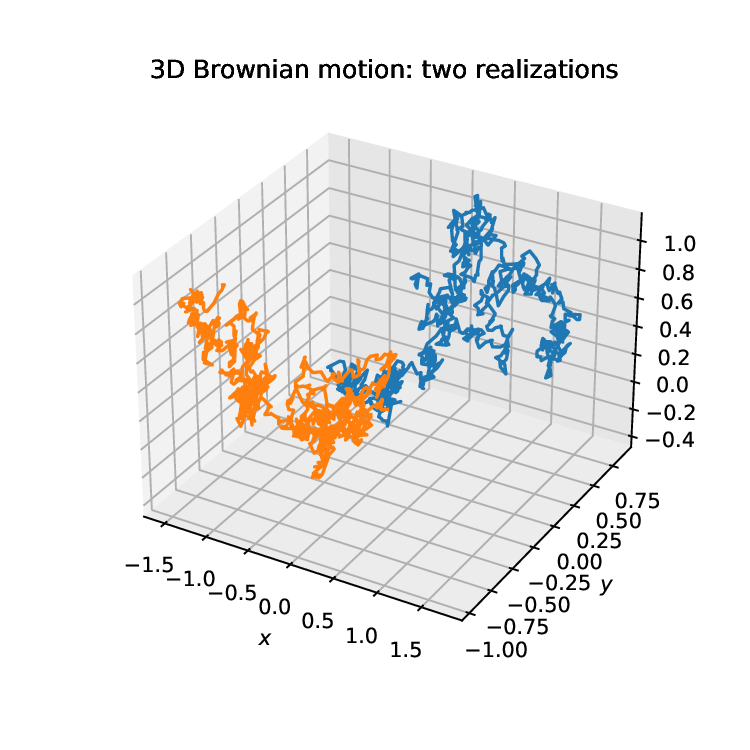}
\captionof{figure}{Two realizations of the Wiener process on time interval $[0,1]$ in 1D (top), 2D (down-left) and 3D (down-right). }\label{fig:wiener_path}.
\end{center}

\begin{example}[Brownian particles]
The motion of particles suspended in a medium, such as liquid or gas, was first observed by the Scottish botanist Robert Brown in 1827. While examining pollen from the plant {\em Clarkia pulchella} under a microscope, he noticed the erratic movement of pollen grains in water where large number of water
molecules ``push'' the particles in different
directions. This particle motion was later understood to be a stochastic process, now referred to as Brownian motion, as we discussed it above.

Now, consider the movement of a single Brownian particle in the water that is confined to a fixed cubic box, say \([-1, 1]^3\). Our goal is to estimate the {\em expected time} required for a Brownian particle, initially positioned at the center of the box, to hit one of the box's boundaries.

To solve this problem, we can use the Monte Carlo method. We simulate a large number of 3D Brownian motion paths, say \(10^4\) realizations. For each simulated path, we record the time taken for the particle to first reach one of the box boundaries. Then we compute the mean of these recorded times to estimate the expected hitting time. The following code provides an implementation for this Monte Carlo simulation approach. 

\begin{shaded}
\vspace*{-0.3cm}
\begin{verbatim}
import numpy as np
dt, N = 0.005, 10000
HitTime = np.empty(N)
for j in range(N):
    k = 0
    W = np.array([0,0,0])
    while True:
      W = W + np.sqrt(dt)*np.random.normal(0,1,3)
      if any(abs(W) >= 1):
        break
      k += 1
    HitTime[j] = dt*k
HitTime_mean = np.mean(HitTime)
HitTime_std = np.std(HitTime)
err = 1.96*HitTime_std/np.sqrt(N)
print('Expected Hitting Time = ', HitTime_mean)
print('Error with 95% probability = ', err)
\end{verbatim}
\vspace*{-0.3cm}
\end{shaded}
In the code, the hitting time is computed by the product of time step $\Delta t$ and number of steps until the boundary is hit.
We choose a small time step $\Delta t = 0.005$ and a large value $N = 10^4$ (number of simulations).
Finally we computed both the mean and the error of estimation with $95\%$ probability using the formula $1.96 s/\sqrt N$ where $s$ is the sample standard deviation. 
The output for a run reported below.
\begin{shaded}
\vspace*{-0.3cm}
\begin{verbatim}
Expected Hitting Time =  0.4856
Error with 95% probability =  0.0061
\end{verbatim}
\vspace*{-0.3cm}
\end{shaded}
\end{example}
\vsp
%\subsubsection*{Generating diffusion processes}
Another class of processes are {\em diffusion processes} which are Markov with a continuous time parameter and
continuous sample pathes, like as Wiener processes. Wiener processes form a basis for defining and generating diffusion processes.
In fact, a diffusion process is defined as the solution of the following {\em stochastic differential equation} (SDE)
\begin{equation}\label{sde_diffproc}
  \mathrm \td X_t = a(t,X_t)\mathrm \td t + b(t,X_t)\mathrm \td W_t
\end{equation}
where $\{W_t: t\geqslant0\}$ is Wiener process and $a(t,x)$ and $b(t,x)$ are some deterministic functions, usually refereed to as
{\em drift} and {\em diffusion} coefficients, respectively. The resulting process for spacial case $a(t,x)=\mu$ and $b(t,x)=\sigma$ is obtained as
$$
X_t = \mu t  + \sigma W_t
$$
and is called a {\em Brownian motion}. The more special case with $\mu=0$ and $\sigma=1$ gives the Wiener process or the standard Brownian motion $X_t = W_t$.
The case $a(t,x) = \mu x$ and $b(t,x) = \sigma x$ results in the {\em geometric Brownian motion}.

For generating a diffusion process (approximately) we can use the explicit {\em Euler–Maruyama method}\footnote{The explicit Euler–Maruyama method is the stochastic equivalent of the explicit Euler method for solving SDEs.} to discretize \eqref{sde_diffproc}
as
\begin{equation*}
  Y_{k+1} = Y_k + a(t_k,Y_k) \Delta t + b(t_k,Y_k)\sqrt{\Delta t} Z_{k+1}, \quad k=0,1,2,\ldots,
\end{equation*}
where $\Delta t$ is a time step, $t_k = k\Delta t$, and $Z_1,Z_2,\ldots$ are independent random variables with $\Nor(0,1)$ distributions.
The process $\{Y_k,k=0,1,\ldots\}$ approximates the exact process $\{X_t,t\geqslant 0\}$ in the sense that $Y_k\approx X_{k\Delta t}$.
The initial variable $Y_0$ should be generated from the distribution of $X_0$. The Python function is given below.
\begin{shaded}
\vspace*{-0.3cm}
\begin{verbatim}
def DiffusionProcessGen(drift, diffusion, tspan, dt, X0, *args): 
    # This function generates a diffusion process by solving the SDE
    #       dX_t = a(t,X_t) dt + b(t,X_t)dW_t,   X_0 = X0
    # for drift coefficient a(t,x) and diffusion coefficient b(t,x)
    # on interval tspan = [t_0, t_end] at equidistance times t_0, t_1, ...,t_n
    # using the explicit Euler–Maruyama method    
    t = tspan[0]
    Y = {}
    Y[0] = X0; j = 0;
    while t <= tspan[1]:
        dW = np.sqrt(dt)*np.random.randn()
        Y[j+1] = Y[j] + drift(t,Y[j],*args)*dt + diffusion(t,Y[j],*args)*dW
        t += dt; j += 1
    Y = [Y[i] for i in range(len(Y))]
    return Y
\end{verbatim}
\vspace*{-0.3cm}
\end{shaded}
\vsp 

%A sample path for the geometric Brownian motion with $\mu = 1$, $\sigma=0.25$, $X_0=1$, and $\Delta t=10^{-3}$ is generated by 
%\begin{shaded}
%\vspace*{-0.3cm}
%\begin{verbatim}
%import matplotlib.pyplot as plt
%def drift(t,x):     # drift coefficient a(x,t)
%    return 1*x
%def diffusion(t,x): # diffusion coefficient b(x,t)
%    return 0.25*x
%Y = DiffusionProcessGen(drift, diffusion, [0,2], .001, 1)
%plt.figure(figsize = (5, 3))
%plt.plot(np.linspace(0,2,len(Y)),Y, linestyle = '-', color='red')
%plt.xlabel('time $t$'); plt.ylabel('$Y_t$'); 
%plt.title('Geometric Brownian motion')
%\end{verbatim}
%\vspace*{-0.3cm}
%\end{shaded}
%\noindent
%and is given in Figure \ref{fig:gbm_path}
%with $\mu = 1$, $\sigma=0.25$, $X_0=1$, $\Delta t=10^{-3}$ and $N = 2000$.

%\begin{center}
%\includegraphics[scale=0.7]{stoch15}
%\captionof{figure}{A typical realization of the geometric Brownian motion %on $[0,2]$.}\label{fig:gbm_path}.
%\end{center}
%\vsp
\begin{example}\label{ex:stock_value}
One application of Brownian motions in combination with the Monte Carlo method is in the pricing of financial options. 
Here we give an example for pricing a {\em European Call Option}. 

In finance, a European call option is a contract that gives the holder the right, but not the obligation, to buy a stock at a fixed price (the strike price \( K \)) at a specified future date (the expiration time \( T \)). The {\em Black-Scholes} model assumes that the stock price follows the SDE
\[
\td S_t = \mu S_t \, \td t + \sigma S_t \, \td W_t,
\]
where \(S_t\) is the stock price at time $t$ (in year), 
\( \mu \) is the drift term representing the average rate of return, \( \sigma \) is the volatility (standard deviation of returns), and \( W_t \) is a Wiener process. As we observe, the stock price \(S_t\) is indeed a geometric Brownian motion. 

The payoff of the call option at the final time is 
\[ 
C_T = \max(S_T - K, 0) 
\] 
for the given strike price \( K \). 
We know from Feynman-Kac that the value of the call option at earlier times $t<T$ is given by
$$
C_t = \E(e^{-r(T-t)}C_T| S_t)
$$
where $r$ is the risk-free interest rate. 
 This expectation is taken under the appropriate risk-neutral measure, which sets the drift $\mu$
equal to the risk-free rate $r$.  
The option price at time $t = 0$ (generally representing the present year) for a given initial stock price $S_0$ is 
$$
C_0  = \E(e^{-rT}C_T) = e^{-rT} \E(C_T).
$$
To apply the Monte Carlo method to estimate $C_0$, we generate $N$ (a large number) simulation paths for \(S_t\) up to time \( t = T \) and for each path we calculate the payoff at expiration, i.e. $C_T$. Finally we compute the mean of $C_T$ and multiply it by $e^{-rT}$. 

As an example, we assume that the expiration time is $T = 0.5$ a year, the current
asset price is $S_0 = 102$, the volatility is $\sigma = 30\%$, the interest rate is $r = 4\%$, and the strike price is \( K = 100\).  
We use $\Delta t = 0.001$ and generates $N = 10^4$ Monte Carlo simulations.
We also estimate the standard error. 

\begin{shaded}
\vspace*{-0.3cm}
\begin{verbatim}
r = 0.04        # risk-free interest rate
mu = r          # drift coefficient (average rate of return) and interest 
sigma = 0.3     # volatility 
S0 = 102        # initial stock price
K = 100         # strike price 
T = 0.5         # expiration time
dt = 0.001      # steplength for Euler-Maruyama method
def drift(t,x,*args):      # drift 
    m = args[0]
    return m*x
\end{verbatim}
%\vspace*{-0.3cm}
\end{shaded}
\begin{shaded}
%\vspace*{-0.3cm}
\begin{verbatim}
def diffusion(t,x,*args):  # diffusion 
    s = args[1]
    return s*x
N = 10**4            # number of MC simulations
CT = np.empty(N)
for j in range(N):
    S = DiffusionProcessGen(drift, diffusion, [0,T], dt, S0, mu, sigma)
    ST = S[-1]
    CT[j] = max(ST-K , 0)
C0 = np.exp(-r*T)*np.mean(CT)
std = np.std(np.exp(-r*T)*CT)
err = 1.96*std/np.sqrt(N)
print("The call value is {0}' +/- {1} with 95% probability".format(V0, err))
\end{verbatim}
\vspace*{-0.3cm}
\end{shaded}

An execution gives the following output:
\begin{shaded}
\vspace*{-0.3cm}
\begin{verbatim}
The call value is 10.0314 +/- 0.2898 with 95% probability
\end{verbatim}
\vspace*{-0.3cm}
\end{shaded}

In the Monte Carlo loop, we generated \( N \) paths of a geometric Brownian motion with an initial value \( S_0 = 102 \). Ten trajectories are shown on the left side of Figure \ref{fig:geo-brown}. Note that only the final value of each path, \( S_T \), was used to compute the payoff of the option.
The histogram of random variable $S_T$ is also shown on the right-hand side of figure \ref{fig:geo-brown} which shows that the distribution of \( S_T \) is a log-normal distribution. This is consistent with the theoretical property of geometric Brownian motion. Note that a random variable \( X \) is said to have a log-normal distribution if \( \log(X) \) has a normal distribution. 

\begin{center}
\includegraphics[scale=.6]{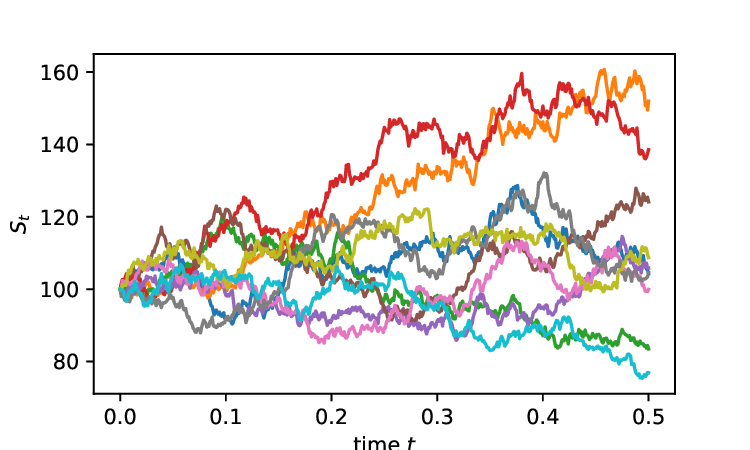}
\includegraphics[scale=.6]{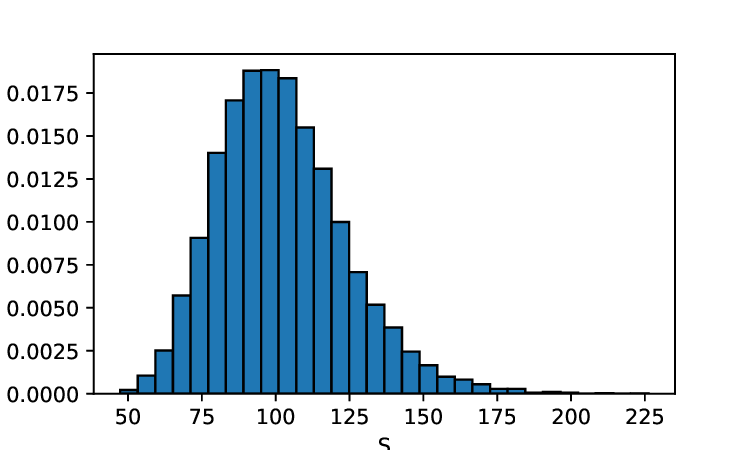}
\captionof{figure}{Ten realizations of the geometric Brownian motion $S_t$ for $t\in [0,0.5]$ (left), the histogram of the variable $S_T$ for $T = 0.5$ (right)}\label{fig:geo-brown}
\end{center}

\end{example}
\vsp

%% file: part_SSA.tex
\section{Stochastic Simulation Algorithm (SSA)}\label{sect:ssa}

The Stochastic Simulation Algorithm (SSA), also known as the {\em Gillespie algorithm}, proposed by Daniel T. Gillespie in 1977\footnote{Daniel T. Gillespie, {\em Exact Stochastic Simulation of Coupled Chemical Reactions}, The Journal of Physical Chemistry, Vol. 81, No. 25, 1977.}, is a powerful method employed in systems biology to predict the behavior of complex biological systems. 
Other areas of application are in epidemiology and ecology. 

\subsection{Simulation of a simple epidemic model}

To illustrate the SSA algorithm, we consider a susceptible-infected-recovered (SIR) model. Such models describe the spread of a virus (for example influenza, covid, etc.) within a population. The population is divided into three groups (states): susceptible individuals, infected individuals, and recovered individuals. Each state is represented by the variables \( S \), \( I \), and \( R \), respectively. 
The relationships between these states can be depicted in a flow diagram (transition graph), where transitions occur based on specific rates. See Figure \ref{fig:sir_model}. 

\begin{center}
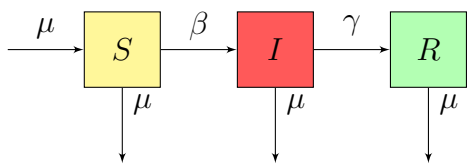

\begin{tikzpicture}[node distance=1cm,auto,>=latex',every node/.append style={align=center},int/.style={draw, minimum size=1cm}]
    \node [int,fill=yellow!50] (S) {$S$};
    \node [int, right=of S,,fill=red!65] (I) {$I$};
    \node [int, right=of I,fill=green!30] (R) {$R$};
    \coordinate[below=of S] (outS);
    \coordinate[left=of S] (inS);
    \coordinate[below=of I] (outI);
    \coordinate[above=of I] (abovI);
    \coordinate[below=of R] (outR);
    \path[->, auto=false] (S) edge node {$\beta$ \\[.2em]} (I)
                          (I) edge node {$\gamma$ \\[.2em]} (R)
                          (S) edge node[right] {$\mu$ \\[.02em]} (outS)
                          (I) edge node[right] {$\mu$ \\[.02em]} (outI)
                          (R) edge node[right] {$\mu$ \\[.02em]} (outR);
   \path[<-, auto=false] (S) edge node {$\mu$ \\[.02em]} (inS);                                                    
\end{tikzpicture}
\captionof{figure}{The transition graph of a simple SIR model.}\label{fig:sir_model}
\end{center}

The model has some parameters: \( \mu \) is the birth and death rate, \( \beta \) is the infection rate, and \( \gamma \) is the recovery rate. The unit of \( S \), \( I \), and \( R \) is individual, and the unit of rates is one over time unit, for example \(\frac{1}{\text{day}}\).

To simplify the mathematical formulation, we assume that the variables \( S(t) \), \( I(t) \), and \( R(t) \) are continuous functions over time \( t \), and the simulation is valid in a defined interval \([0, t_{\text{final}}]\). Additionally, it is assumed that the rates of birth and death are equal, and all newborns are susceptible. The model does not account for pathogen-induced mortality (death because of the virus), and it is assumed that recovered individuals remain immune throughout the period of the epidemic. Furthermore, the infection rate depends on the ratio of infected individuals to the total population, expressed as \( \beta \frac{I}{N} \). This assumption is quite reasonable, as a higher number of infected individuals increases interactions between the susceptible group and the infected group which leads to an elevated infection rate over time.

Given the assumption above, we can write a system of ordinary differential equations (ODEs) to describe the dynamics of the SIR model. We denote the total population at time \( t \) by 
\[ N(t) = S(t) + I(t) + R(t). \] 
The differential equations governing the dynamics are as follows:

\begin{align*}
\frac{dS}{dt} &= \mu N - \mu S - \beta \frac{I}{N} S \\
\frac{dI}{dt} &= \beta \frac{I}{N} S - \mu I - \gamma I \\
\frac{dR}{dt} &= \gamma I - \mu R
\end{align*}
with  initial conditions 
$$ 
  S(0) = S_0, \quad I(0) = I_0 , \quad R(0) = R_0. \vsp 
$$  
This initial value problem is a deterministic model for the given phenomenon. 
To solve this system of equations, we can employ a deterministic numerical methods (an ODE solver) such as a Runge-Kutta method. 
Here we call the \texttt{RK45} from the \texttt{scipy.integrate.solve\_ivp} module. 

\begin{shaded}
\vspace*{-0.3cm}
\begin{verbatim}
import numpy as np
import matplotlib.pyplot as plt
from scipy.integrate import solve_ivp

mu, bet, gam = 1e-4, 0.25, 0.05   # rates
Initial = [198,2,0]               #   [S(0), I(0) R(0)]
FinalTime = 120                   # final time of simulation

def ODEfun(t,y):
    yprime = np.zeros(3); 
    S,I,R = y
    N = np.sum(y)
    yprime[0] = mu*N - bet*S*I/N-mu*S
    yprime[1] = bet*S*I/N -(mu+gam)*I
    yprime[2] = gam*I - mu*R
    return yprime

teval = np.linspace(0, FinalTime,500)
sol = solve_ivp(ODEfun, [0,FinalTime], Initial, t_eval = teval)

plt.figure(figsize = (6, 4))
plt.plot(sol.t,sol.y[0],linestyle = 'solid', color='blue', label = '$S$')
plt.plot(sol.t,sol.y[1],linestyle = 'solid', color='red', label = '$I$')
plt.plot(sol.t,sol.y[2],linestyle = 'solid', color='green', label = '$R$')
plt.xlabel('time $t$'); plt.ylabel('Individuals')
plt.title('Deterministic solution using RK45')
plt.legend(loc='center right')
\end{verbatim}
\vspace*{-0.3cm}
\end{shaded}

The results are given in Figure \ref{fig:ssa_deter} which illustrates the dynamics of the susceptible, infected, and recovered populations over time.
The number of infected individuals increases from initially two to approximately 125 and reaches its peak around the 20th day. Following this peak, the epidemic enters a decline phase, with the number of infections decreasing rapidly until the epidemic eventually comes to an end.

\begin{center}
\includegraphics[scale=.5]{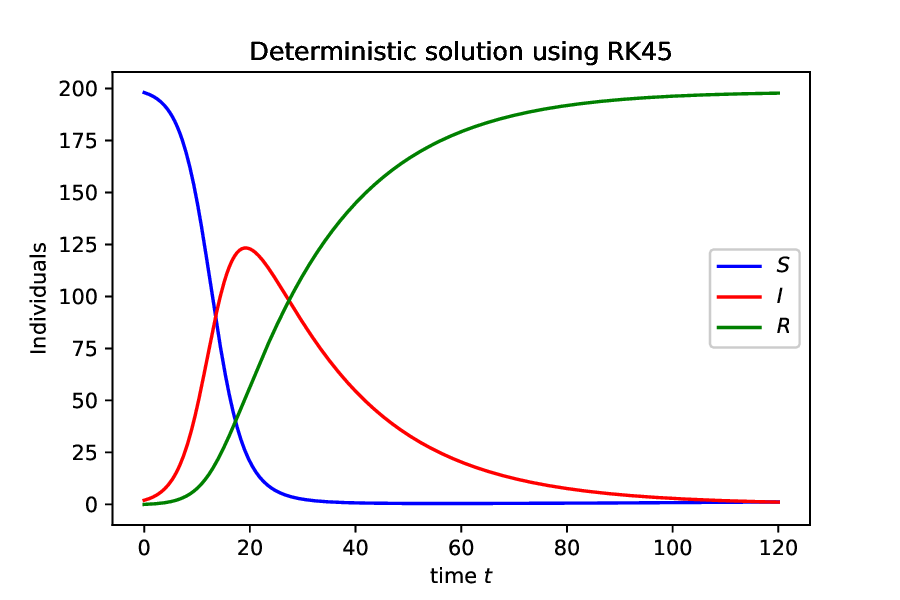}
\captionof{figure}{Solution of SIR model using \texttt{RK45}.  }\label{fig:ssa_deter}
\end{center}

While the above deterministic model provides valuable insights, it overlooks the discrete nature of populations and the inherent uncertainties present in the model. A stochastic approach reformulates the model by replacing continuous variables \( S(t) \), \( I(t) \), and \( R(t) \) with discrete values. We define a series of reactions (processes) with their associated {\em propensity functions} which are indeed the number of individuals moving into or out of the groups:
\[
\begin{array}{r r c ll ll}
1.&\emptyset &\xrightarrow{\;\quad \;} & S &&&  w_1 = \mu N \\
2.&S& \xrightarrow{\;\quad \;}& I &&&  w_2 = \beta \frac{I}{N} S \\
3.&I& \xrightarrow{\;\quad \;}& R &&&  w_3 = \gamma I \\
4.&S &\xrightarrow{\;\quad \;}& \emptyset &&&  w_4 = \mu S \\
5.&I &\xrightarrow{\;\quad \;}& \emptyset &&&  w_5 = \mu I \\
6.&R &\xrightarrow{\;\quad \;}& \emptyset &&&  w_6 = \mu R
\end{array}
\]
%\vspace{.2cm} 

\noindent
The reactions include transitions from susceptible to infected individuals, from infected to recovered individuals, and the natural birth and death processes affecting all the groups.
The propensity functions represent the tendency (or likelihood) of each reaction to occur.
We now assume a state vector
$$
 {\y}(t) = (y_1(t), \ldots, y_n(t)) 
$$
 containing the state variables. For instance, for the SIR model we have \( n = 3 \) states and \( m = 6 \) reactions, with the propensity functions corresponding to each reaction:
\begin{align*}
&\y = (S, I, R),\\
&w_1 = \mu N, \; w_2 = \beta\frac{I}{N}S ,\; w_3 = \gamma I, \; w_4 = \mu S,\; w_5 = \mu I, \; w_6 = \mu R .
\end{align*}
Also we compute the {\em total propensity} function  
$$
a(\y) = \sum_{j=1}^m w_j(\y)
$$
which represents the overall tendency of the system to evolve and undergo changes.
The positive propensity functions $w_j$ represent the relative likelihood of each reaction occurring compared to the others. To convert these into probabilities, we define
$$
p_j := \frac{w_j}{a}, \quad j = 1,\ldots, m
$$
where $a$ is the total propensity. This allows us to form
the following discrete distribution table
\begin{center}
\begin{tabular}{c|cccc}
$reaction~j$   & $1$ & $2$ & $\ldots$ & $m$ \\ \hline
$probability~p_j$ & $p_1$ & $p_2$ & $\ldots$ & $p_m$ \\
\end{tabular}
\end{center}
which tells us which reactions are more likely to occur and change the system states.

The basic assumption in the SSA is that, between time \( t \) and \( t + \tau \), exactly one reaction occurs and causes a change in the system's states by either \( 0 \), \( +1 \), or \( -1 \). Therefore, at any given time \( t \), the SSA includes three steps:

\begin{enumerate}

\item {\em when} the next reaction will occur (i.e. the value of $\tau$)
\item {\em which} reaction will occur
\item {\em Updating the state} of the system
\end{enumerate}
\noindent
The waiting time \( \tau \) for the next reaction is sampled from an exponential distribution with rate parameter $\lambda = a(\y)$:
$$
when \sim \Exp\left(a(\y)\right).
$$
A higher value of $a$ (a higher total propensity) corresponds to a shorter expected waiting time between events. 

Which reaction? The reactions with higher probabilities are more likely to occur. To determine which reaction takes place, we sample from the discrete distribution based on the probabilities \( p_j \) defined earlier:
$$
which \sim \mathcal{DD}\big([1,2,\ldots,m],[p_1,\ldots,p_m]\big).
$$
Once the reaction is determined, the final step is to update the system's state. We use the {\em state-change} vectors $\vv_j$ (also called the {\em stoichiometry} vectors), which represent how each reaction alters the state of the system. For example, in the SIR model, the vectors are as follows:
\begin{equation*}\label{stoch1}
\begin{array}{r r c l l}
1.&\emptyset &\xrightarrow{\;\quad \;} & S, &\quad \vv_1 = [\;1,\;\;\; 0,\;\;\; 0]\\
2.&S& \xrightarrow{\;\quad \;}& I, & \quad \vv_2 = [-1,\;\; 1,\;\; 0]\\
3.&I& \xrightarrow{\;\quad \;}& R,   & \quad \vv_3 = [0,\; -1,\;\; 1]\\
4.&S &\xrightarrow{\;\quad \;} & \emptyset ,&\quad \vv_4 = [-1,\;\;\, 0,\;\, 0]\\
5.&I &\xrightarrow{\;\quad \;} & \emptyset ,&\quad \vv_5 = [\; 0,\; -1,\;\, 0]\\
6.&R &\xrightarrow{\;\quad \;} & \emptyset ,&\quad \vv_6 = [\; 0, \;\; 0,\;-1 ]\\
\end{array}
\end{equation*} 
\vspace{.2cm}

These vectors describe how the number of susceptible, infected, or recovered individuals changes due to each specific reaction.

In the {\em which} phase if the reaction index \( k \) is sampled, then the state is updated based on the corresponding state-change vector \( {\vv}_k \), i.e.,
$$
\y(t+\tau) = \y(t) + \vv_k.
$$ 

Considering all these together, the Gillespie algorithm is presented below. 

%{{\small
%\begin{algorithm}[H]
%\SetKw{KwSet}{Set}
%\KwData{Initial state $\y=\y_0$, final time $t_{final}$, propensity functions $w_1,\ldots,w_m$ and state change vectors $\vv_1,\ldots,\vv_m$
%}
%Set $t =0$\;
%\While{$t<t_{final}$}{
%Compute $a(\y)=w_1(\y)+\cdots+w_m(\y)$ and $p_j(\y)=w_j(\y)/\a_j(\y)$\;
%Generate $\tau \sim \Exp(a(\y))$ \;
%Generate $k \sim  \mathcal{DD}\big([1,\ldots,m],[p_1,\ldots,p_m]\big)$\;
%Update $t \leftarrow t+\tau$\;
%Update $\y \leftarrow \y + \vv_k$
%}
%\caption{Gillespie algorithm}
%\end{algorithm}
%}}

\begin{algorithm}
\caption{Gillespie algorithm (SSA)}
\begin{algorithmic}
\Require Initial state $\y=\y_0$, final time $t_{final}$, propensity functions $w_1,\ldots,w_m$ and state change vectors $\vv_1,\ldots,\vv_m$
\Ensure State vector $\y(t)$ at final time $t = t_{final}$
\State $t \gets 0$
\While {$t\leq t_{final}$}
    \State Compute $a(\y)=w_1(\y)+\cdots+w_m(\y)$ and $p_j(\y)=w_j(\y)/a(\y)$
    \State Generate $\tau \sim \Exp(a(\y))$ 
    \State Generate $k \sim  \mathcal{DD}\big([1,\ldots,m],[p_1,\ldots,p_m]\big)$
    \State Update $t \gets t+\tau$
    \State Update $\y \gets \y + \vv_k$ 
\EndWhile
\end{algorithmic}
\end{algorithm}

\subsection{Python implementation}

In this section, we implement the SSA using Python. At each time step, the code calls the \texttt{RandExp} and \texttt{RandDisct} functions (Section \ref{sect:generation}) to sample the steplength and determine which reaction occurs.

\begin{shaded}
\vspace*{-0.3cm}
\begin{verbatim}
##  Gillespie algorithm (SSA)
def SSA(Initial, StateChangeMat, FinalTime):
    # Inputs:
     #  Initial: initial conditins of size (StateNo x 1)
     #  StateChangeMat: State-change matrix of size (ReactNo, StateNo)
     #  FinalTime: the maximum time we want the process be run
    # Outputs:
     #  AllTimes: the dict. of all selected time levels 
     #  AllStates: the dict. of all state values at corresponding time levels
    [m,n] = StateChangeMat.shape
    ReactNum = np.array(range(m))
    AllTimes = {}   # define a dict. for storing all time levels
    AllStates = {}  # define a dict. for storing all states at all time levels
    AllStates[0] = Initial
    AllTimes[0] = [0]
    k = 0; t = 0; State = Initial
    while True:
        w = PropensityFunc(State, m)     # propensities
        a = np.sum(w)
        tau = RandExp(a,1)               # WHEN the next reaction happens
        t = t + tau                      # update time
        if t > FinalTime:
            break
        which = RandDisct(ReactNum,w/a,1)             # WHICH reaction occurs
        State = State + StateChangeMat[which.item(),] # Uppdate the state
        k += 1
        AllTimes[k] = t
        AllStates[k] = State
    return AllTimes, AllStates
\end{verbatim}
\vspace*{-0.3cm}
\end{shaded}
The above \texttt{SSA} function works for a general stochastic model provided that the propensity functions and state-change vectors are given.  
Below, we will specify these functions for the SIR model discussed earlier.

\begin{shaded}
\vspace*{-0.3cm}
\begin{verbatim}
mu, bet, gam = 1e-4, 0.25, 0.05  # rates
Initial = [198,2,0]              # initial values
FinalTime = 120                  # final time 
def PropensityFunc(State, ReactNo):
    S,I,R = State 
    N = S + I + R;
    w = np.zeros(ReactNo)
    w[0] = mu * N              # birth (newborns)
    w[1] = bet/N * S * I       # infection
    w[2] = gam * I             # recovery
    w[3] = mu * S              # death of susceptible individuals
    w[4] = mu * I              # death of infected individuals
    w[5] = mu * R              # death of recovered individuals
    return w
StateChangeMat = np.array([
                    [+1,  0,  0],
                    [-1, +1,  0],
                    [ 0, -1, +1],
                    [-1,  0,  0],
                    [ 0, -1,  0],
                    [ 0,  0, -1]])
\end{verbatim}
\vspace*{-0.3cm}
\end{shaded}

Finally, we call the \texttt{SSA} function to run multiple simulations, for instance $N = 10$, and subsequently plot the results.

\begin{shaded}
\vspace*{-0.3cm}
\begin{verbatim}
N = 10   # number of simulations
plt.figure(figsize = (6, 4))
for k in range(N):
    Time, States = SSA(Initial, StateChangeMat, FinalTime) 
    n = len(Time)
    t = [Time[i][0] for i in range(n)]
    S = [States[i][0] for i in range(n)]
    I = [States[i][1] for i in range(n)]
    R = [States[i][2] for i in range(n)]
    plt.plot(t,S,linestyle = '-', color='blue')
    plt.plot(t,I,linestyle = '-', color='red')
    plt.plot(t,R,linestyle = '-', color='green')
plt.xlabel('Time');
plt.ylabel('Individuals');
plt.title('Stochastic solutions using SSA')
plt.legend(['$S$','$I$','$R$'],loc='center right')
plt.show()
\end{verbatim}
\vspace*{-0.3cm}
\end{shaded}

The results are given in Figure \ref{fig:ssa_stoch} for 10 simulations of the above SSA algorithm. This figure represents the stochastic nature of the SIR model. 

\begin{center}
\includegraphics[scale=.5]{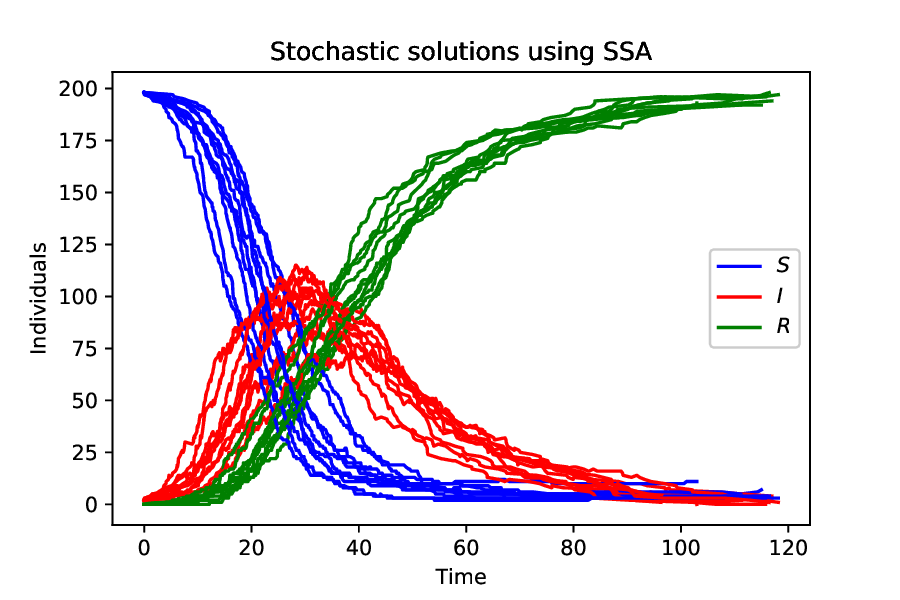}
\captionof{figure}{Stochastic simulation of the SIR model using the Gillespie algorithm: ten simulations}\label{fig:ssa_stoch}
\end{center}

One can execute the code and obtain many simulations (stochastic processes) and finally compute the mean of the solutions as a Monte Carlo solution. In many cases the mean is close to the deterministic solution of the counterpart ODE model. 

To speed up the SSA, Daniel Gillespie proposed another algorithm called {\em tau-leaping}\footnote{Daniel T. Gillespie, {\em Approximate accelerated stochastic simulation of chemically reacting systems}, The Journal of Chemical Physics. 115 (2001) 1716-1733.  
}, which replaces the variable waiting times between reactions with a fixed step size \( \tau \) for each iteration. Unlike the standard SSA, where only one reaction occurs per step, tau-leaping allows multiple reactions to occur in a single time step, and the state changes can be larger than just 1. 
For each reaction \( j \), the state change is sampled from a Poisson distribution with parameter \( \lambda = w_j \tau \), where \( w_j \) is the propensity of reaction \( j \), and \( \tau \) is the step length. Such Poisson sample represents the number of events expected in an interval of length \( \tau \) with propensity \( w_j \). 
This algorithm requires careful handling to avoid negative populations when some states are close to zero and state changes are large. Additionally, selecting an optimal step size \( \tau \) is a key challenge because steps that are too large can lead to inaccuracies, while steps that are too small reduce the performance benefit of the algorithm.

Finally, we note that the Gillespie algorithm is also available in the \texttt{GillesPy2} library, which provides an object-oriented framework for building and simulating the mathematical model. 
The methods include the SSA, the tau-leaping, and some numerical ODE solvers.
The library is optimized for performance, and is written in C++ and NumPy. For more details check this link\footnote{\url{https://gillespy2.readthedocs.io/en/latest/}}.

\subsection{Application to biochemical kinetics}

In biochemical systems, a finite number of particles present in living cells and the inherent randomness associated with molecular interactions and reaction rates lead to a need for models that account for discrete and stochastic dynamics rather than continuous and deterministic ones.
 
Consider a model comprised of \(n\) species, denoted as \(\{S_1, S_2, \ldots, S_n\}\), which interact through \(m\) chemical reactions, represented as \(\{r_1, \ldots, r_m\}\). As before, the state of the system can be described by the vector \({\y}(t) = (y_1(t), \ldots, y_n(t))\), with the initial condition specified as \({\y}(0) = {\y}_0\). 
The interactions are modeled by reactions, propensity functions, and state-change vectors. We assume that each reaction \(r_\ell\) is {\em elemental} and  is either {\em unimolecular} or {\em bimolecular}. For simplicity, let us consider a system with three species with state variable \({\y} = (y_1, y_2, y_3)\). The following table illustrates various reaction cases, where \(c\) denotes a constant rate:

\[
\begin{array}{|c|c|c|c|}
\hline
\text{Reaction} & \text{State-change vector} & \text{Propensity function} \\
\hline
y_1 \xrightarrow{\;c\;} y_2 & {\vv} = [-1, 1, 0] & w = c y_1 \\
y_1 \xrightarrow{\;c\;} y_2 + y_3 & {\vv} = [-1, 1, 1] & w = c y_1 \\
y_1 \xrightarrow{\;c\;} y_1 + y_2 & {\vv} = [0, 1, 0] & w = c y_1 \\
y_1 \xrightarrow{\;c\;} 2y_1 & {\vv} = [1, 0, 0] & w = c y_1 \\
y_1 \xrightarrow{\;c\;} \emptyset & {\vv} = [-1, 0, 0] & w = c y_1 \\
y_1 + y_2 \xrightarrow{\;c\;} y_3 & {\vv} = [-1, -1, 1] & w = c y_1 y_2 \\
2y_1 \xrightarrow{\;c\;} y_2 & {\vv} = [-1, 1, 0] & w = c {y_1(y_1 - 1)}/{2} \\
2y_1 \xrightarrow{\;c\;} y_1 & {\vv} = [-1, 0, 0] & w = c {y_1(y_1 - 1)}/{2} \\
\emptyset \xrightarrow{\; c\;} y_1 & {\vv} = [1, 0, 0] & w = ? \\
\hline
\end{array}
\]
\vsp 

The propensity functions reflect the concentration of the reactants and emphasize how molecular interactions dictate the rates of reaction.

As an example consider the {\em Michaelis-Menten system} which is a standard model for enzyme-catalyzed reactions. In this model, we consider a substrate \(S\), enzyme \(E\), the enzyme-substrate complex \(C\), and the product \(P\). The model below shows the series of reactions occurring in this system.

\begin{center}
\begin{tikzpicture}[>=latex,scale=1]
\node at (0.25,0) {$S+E$}; 
\draw[-to,black,thick] (1,.05) -- (1.7,.05);
\draw[to-,black,thick] (1,-.05) -- (1.7,-.05);
\node at (2,0) {$C$};

\draw[-to,black,thick] (2.25,0) -- (2.95,0);
\node at (3.6,0) {$P+E$};
\node at (1.35,.25) {$c_1$};
\node at (1.35,-.25) {$c_2$};
\node at (2.6,.2) {$c_3$};
\end{tikzpicture}
\end{center}

\begin{center}
\includegraphics[scale=.45]{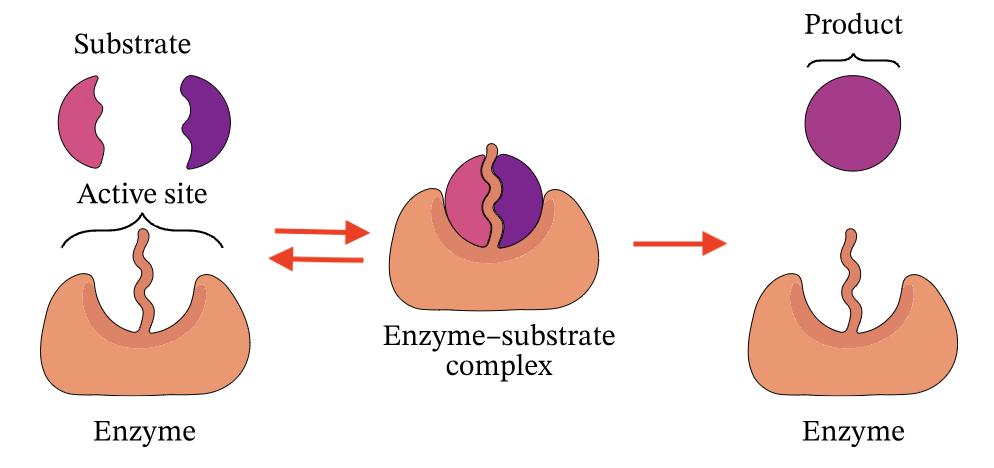}
\captionof{figure}{The Michaelis-Menten reactions (image from \texttt{www.nagwa.com/en/})}\label{fig:Michaelis-Menten}
\end{center}

If we denote the state vector by \({\y} = (S, E, C, P)\) then the reactions can be summarized as follows:

\begin{table}[!h]
\centering
\begin{tabular}{llll}
$r_1:$ & 
$S+E \xrightarrow{\;c_1\;} C$,
& $\vv_1 = [-1,-1,+1,0]$, & $w_1 = c_1SE = c_1y_1y_2$ \\
$r_2:$ & 
$S+E \xleftarrow{\;c_2\;} C$,
& $\vv_2 = [+1,+1,-1,0]$, & $w_2 = c_2C = c_2y_3$ \\
$r_3:$ & 
$C \xrightarrow{\;c_3\;} P+E$,
& $\vv_3 = [0,+1,-1,+1]$, & $w_3 = c_3C = c_3y_3$ 
\end{tabular}

\end{table}

The first reaction involves the formation of the enzyme-substrate complex. The second reaction is the reverse process, where the complex dissociates back into the substrate and enzyme. The third reaction describes the conversion of the complex into product and enzyme. For given constant rates $c_j$ and initial numbers of proteins, we can simply solve this model using the Gillespie algorithm. 
 \vsp 

\subsection{Lotka-Volterra models}

Other important examples are the {\em Lotka-Volterra models} which are used to describe predator-prey interactions. 
Predator and prey animals, such as hawks and mice or foxes and rabbits, interact in an ecosystem where predator eat prey. If the prey population is large, food becomes easily available for the predator and the population grows. This leads to a decrease in the prey population, which, in turn, leads to reduced access to food, resulting in a decrease in the predator population. This means that more prey survives, and so on. This {\em process} continues indefinitely. 

As an example of a simple Lotka-Volterra model assume that $R$ represents the number of prey, e.g. rabbits, and $F$ represents the number of predators, e.g. foxes. To analyse the predator-prey system, we consider the following set of reactions:

\begin{align*}
    R &\xrightarrow{\;\alpha\;} 2R \\
    R+F &\xrightarrow{\;\beta\;} 2F \\
    F &\xrightarrow{\;\gamma\;} \emptyset
\end{align*}

The first reaction incorporates prey reproduction into our model, with $\alpha$ denoting the reproduction rate of each prey. The second reaction incorporates predator reproduction per prey, with $\beta$ denoting the rate of interactions between prey and predator animals. In other words, the predator consumes prey and reproduces at a rate of $\beta$. The third and final reaction incorporates predator mortality, with $\gamma$ denoting the rate at which predators are removed from the ecosystem.
Now, we can apply the Gillespie algorithm to solve this simple model provided that the initial populations $F_0$ and $R_0$ and the model rates are given. This task is left as an exercise to the reader.

%% file: part_MCMC.tex
\section{Markov chain Monte Carlo (MCMC)}

In preceding sections we have typically generated iid random variables directly from the density of interest $f$. In this section the
{\em Markov chain Monte Carlo (MCMC)} is introduced as a powerful tool to {\em approximately} generate samples from an {\em arbitrary}
distribution\footnote{
The MCMC method was motivated by the pioneer work of Metropolis, et. al. in 1953:\newline
M. Metropolis, A. W. Rosenbluth, M. N. Rosenbluth, A. H. Teller, and E. Teller,
{\em Equations of state calculations by fast computing machines,} Journal of Chemical Physics, 21(1953) 1087-1092.\\
After that, many modifications were done on the original MCMC algorithm, notably the Hastings algorithm in 1970:\\
W. K. Hastings, {\em Monte Carlo sampling methods using Markov chains and their applications}, Biometrika, 57(1970) 92-109.\\
An alternative methodology then proposed by German and German in 1984 which is known as Gibbs sampler:\\
S. Geman and D. Geman, {\em Stochastic relaxation, Gibbs distribution and the Bayesian
restoration of images}, IEEE Transactions on PAMI, 6(1984) 721-741.
}.
The main idea behind the MCMC algorithms is to simulate a Markov chain such that its stationary distribution {approximately} coincides with the desired distribution $f$. 
%Here, we give couple of such algorithms namely the {\em Metropolis-Hastings} algorithm and the {\bf Gibbs sampler} algorithm.
One of the MCMC algorithms is the {\em Metropolis-Hastings} algorithm which is discussed in detail here. 
We follow \cite{Rubinnstein-Kroese:2017} in this section. 
\vsp 

\subsection{Metropolis-Hastings algorithm}
Given a target density $f$, we want to generate a Markov chain $\{X_t:t=0,1,\ldots\}$
with stationary distribution $f$.
For simplicity, we start with a general discrete distribution.
Assume that we want to generate a random variable $X$ which takes its values in sample space 
$$
\Ss = \{1,2,\ldots,m\}
$$ 
with target distribution 
$$\{\pi_1,\pi_2,\ldots,\pi_m\}.
$$ 
In the Metropolis-Hastings algorithm a Markov chain $\{X_t: t=0,1,\ldots\}$ on $\Ss$ is simulated based on a primary transient matrix $Q=(q_{ij})$. This matrix is used to approximate the actual transient matrix of the chain.
The algorithm contains the following steps:
\begin{itemize}
  \item[(1)] [{\em Variable generation}] Given $X_t=i$, generate a random variable $Y$ such that $\pr(Y=j)=q_{ij}$ for all $j\in \Ss$. On the other words, generate $Y$ from the $i$-th row of $Q$.
  \item [(2)] [{\em Accept or reject}] If $Y = j$ then accept $Y$ with probability $\alpha_{ij}$ and let $X_{t+1}=Y$, otherwise reject $Y$ and let $X_{t+1}=X_t$, where
      \begin{equation}\label{alpha_accept}
      \alpha_{ij} = \min\left\{\frac{\pi_j}{\pi_i}\frac{q_{ji}}{q_{ij}},1 \right\}.
      \end{equation}
\end{itemize}
Since $X_{t+1}$ is obtained from $X_t$ only, the chain is Markov. Besides, the transient matrix $P=(p_{ij})$ of the chain is
\begin{equation} \label{pij_metro}
  p_{ij} = \begin{cases} q_{ij}\alpha_{ij}, & i\neq j\\ \ds 1-\sum_{k\neq i}q_{ik}\alpha_{ik},& i=j  \end{cases},
\end{equation}
because, for $i\neq j$ we can write
\begin{align*}
  p_{ij} & =\, \pr(X_{t+1}=j\,|\,X_t = i) \\
   &=\, \pr(Y=j, Y\,\mbox{is accepted}\,|\, X_t=i)\\
   &=\, \pr(Y=j\,|\,X_t=i)\times \pr(Y\,\mbox{is accepted}\,|\, Y=j, X_t=i)\\
   &=\, q_{ij}\times \alpha_{ij}.
\end{align*}
In the third equality we have used the identity $\pr(A\cap B|C)=\pr(A|C)\pr(B|A\cap C)$. The case $i= j$ in \eqref{pij_metro} follows from the fact that each row of $P$ sums up to unity. Using \eqref{pij_metro} and the definition of $\alpha_{ij}$ in \eqref{alpha_accept} we have
\begin{equation}\label{detailed_balance}
  \pi_i p_{ij}= \pi_j p_{ji},\quad i,j\in \Ss
  \vsp
\end{equation}
which is the {\em detailed balance equation} for the Markov chain.
The proof is as follows. First we have
\begin{equation*}
  \pi_i p_{ij} = \pi_i q_{ij}\alpha_{ij} = \pi_i q_{ij} \frac{\pi_j}{\pi_i}\frac{q_{ji}}{q_{ij}} = \pi_jq_{ji}
\end{equation*}
provided that $\frac{\pi_j}{\pi_i}\frac{q_{ji}}{q_{ij}}\leqslant 1$. On the other hand
\begin{equation*}
  \pi_j p_{ji} = \pi_j q_{ji}\alpha_{ji} = \pi_jq_{ji}
\end{equation*}
because $\alpha_{ji}=1$ from the fact that $\frac{\pi_i}{\pi_j}\frac{q_{ij}}{q_{ji}}\geqslant 1$.

Taking summation on $j$ from both sides of equation \eqref{detailed_balance} gives $\pi_i= \sum_{j}\pi_j p_{ji}$ for $i\in \Ss$, or 
$\bpi=\bpi P$, which means that the chain has stationary probability $\bpi=[\pi_1,\ldots,\pi_m]$. This stationary distribution is also a limiting distribution if the Markov chain is irreducible and aperiodic.

An important advantage of the above Metropolis-Hastings algorithm is that the target distribution $\bpi$ needs to be only known up to a normalization constant $C$, because in the definition of $\alpha_{ij}$ (the only place $\{\pi_j\}$ is used) any constant $C$ will be cancelled in quotient $\frac{\pi_j}{\pi_i}$.

The above algorithm can be generalized to generate samples from an arbitrary multi-dimensional density $f(x)$ instead of $\{\pi_j\}$.
From here on $x$, $y$, $X$ and $Y$ with or without subscripts are $d$-dimensional variables. 
In this generalization, the non-negative transient kernel $q(x,y)$
will be used instead of the transient matrix $Q$. Since the transient kernel is a conditional pdf, we can also write $q(y|x)$ instead of
$q(x,y)$. This kernel is usually called the {\em proposal} function, and plays a role similar to proposal distribution $g$ in the acceptance-rejection method of section \ref{sect:accept-reject}.

The Metropolis-Hastings algorithm starts with an initial state $X_0$ and a target pdf $f(x)$, a proposal function $q(x,y)$ and the number of required samples $N$ as inputs and generates a Markov chain $X_1,X_2,\ldots,X_N$ approximately distributed according to $f(x)$.
The algorithm is given below.

\begin{algorithm}
\caption{Metropolic-Hastings Algorithm}\label{alg:metropolis-hastings}
\begin{algorithmic}
\Require Target distribution $f$, Proposal distribution $q$, Initial state $X_0$, Number of samples $N$
\For {$t = 1,2,\ldots, N$} 
\State 1. Given $X_t$, generate $Y\sim q(X_t,y)$
\State 2. Set $\ds \alpha =\min\left\{\frac{f(Y)}{f(X_t)} \frac{q(Y,X_t)}{q(X_t,Y)} ,1\right\} $
\State 3. Generate $U\sim \Uni(0,1)$
\State 4. If $U\leqslant \alpha$ accept $Y$ and set $X_{t+1}=Y$, otherwise reject $Y$ and set $X_{t+1}=X_t$.
\EndFor
\Ensure The sequence $X_1,X_2,\ldots,X_N$
\end{algorithmic}
\end{algorithm}

Starting by $X_0$, we continue this process until $X_N$ is generated. The sequence $X_1,\ldots,X_N$ is a set of dependent random variables and $X_t$ for large $t$ is approximately distributed according to $f(x)$.

The original Metropolis algorithm uses a proposal function $q$ with symmetrical property $q(x,y)=q(y,x)$, while the modified version by Hastings allows the nonsymmetric kernels as well. With a symmetric $q$, the probability $\alpha$ reduces to
\begin{equation}\label{alpha_prob_fg}
\alpha = \min\left\{\frac{f(Y)}{f(X_t)},1\right\}.
\end{equation}
This does not mean that $q$ is ruled out because $Y$ is still generated from $q$.
%\vsp 

%\subsubsection*{Different types of proposal function}
The simplest proposal function is to take
\begin{equation*}
q(x,y) = g(y)
\end{equation*}
for some pdf $g(y)$. Using this proposal function, in step (1) of the Metropolis-Hastings algorithm, $Y$ is generated from $g(y)$ independent of the current variable $X_t$. The acceptance probability $\alpha$ then is
$$
\alpha = \min\left\{\frac{f(Y)}{f(X_t)}\frac{g(X_t)}{g(Y)},1\right\}\vsp
$$
which depends on $X_t$. Thus the chain still produces dependent samples.

In a {\em random walk sampler} the current state $Y$ for a given state $x$ is given by $Y = x+Z$ where $Z$ is generated from a radially symmetric
distribution such as $\Nor(0,\Sigma)$. For this case $Y$ is indeed generated from $Y\sim \Nor(x,\Sigma)$, i.e.,
$$
q(x,y) = (2\pi)^{-d/2} \exp\left(-\tfrac{1}{2}(x-y)^T\Sigma (x-y)\right).
$$
Since the proposal function is symmetric the acceptance probability is reduced to \eqref{alpha_prob_fg}.
A Phyton code for MCMC with multidimensional normal random walk sampler is given here. Inputs are the desired probability density function \verb+pdf+ we aim to sample from, the initial state \verb+X0+, the covariance matrix of the proposal function $q$ and the number of samples we ask for. The output is a Markov chain \verb+X+. The \texttt{RandMultiNormal} function is given in Section \ref{sect:randvectorgen}.

\begin{shaded}
\vspace*{-0.3cm}
\begin{verbatim}
def McMcRandWalkGen(pdf, X0, SigmaWalk, N):
    dim = np.size(X0)
    X = np.zeros([dim,N])
    X[:,0] = X0
    for t in range(N-1):
        Z = RandMultiNormal(np.zeros(dim),SigmaWalk,1).T
        Y = X[:,t] + Z
        Xt = np.array([X[:,t]])
        alpha = min(pdf(Y)/pdf(Xt),1)
        U = np.random.rand()
        if U <= alpha:
            X[:,t+1] = Y
        else:
            X[:,t+1] = X[:,t]
    return X
\end{verbatim}
\vspace*{-0.3cm}
\end{shaded}
\vsp 

\begin{example}[Rubinnstein-Kroese:2017] \label{ex_metro-hasting}
Consider a random vector $X=[X_1,X_2]$ with the following bivariate pdf
\begin{equation}\label{fx1x2_def}
  f(x,y) = c \exp(-(x^2y^2+x^2+y^2-8x-8y)/2)
\end{equation}
where $c \doteq 1/20216.335877$ is the normalization constant. 
 \begin{center}
 \includegraphics[scale=0.30]{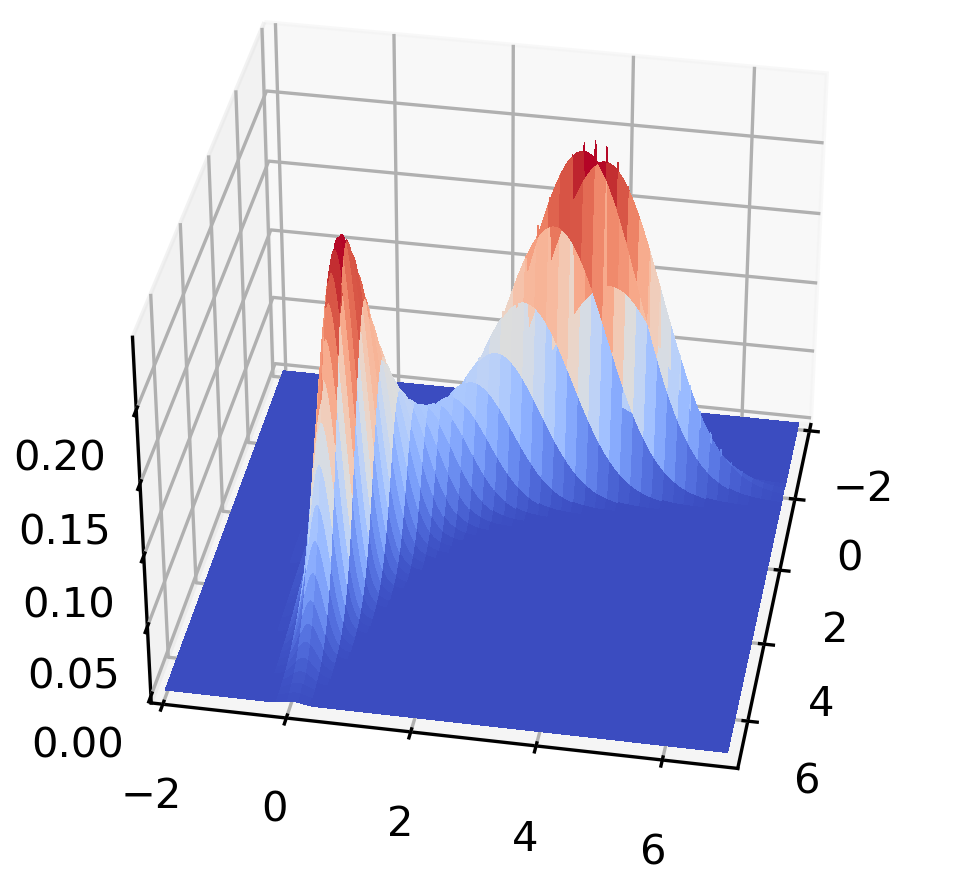}
 \captionof{figure}{Surface plot of bivariate density function $f(x,y)$.}\label{fig:pdf2dex}
 \end{center}
 
The surface graph of this density is shown in Figure \ref{fig:pdf2dex}.
We want to generate samples $\{X_t=(X_{t1},X_{t2}): t=1,2,\ldots,N\}$ from $f(x,y)$ using the Metropolis-Hastings algorithm with the random walk sampler as a transition kernel. We assume that $\Sigma = \mathrm{diag}[\sigma^2,\sigma^2]$ where the moderate value $\sigma = 2$ is chosen in our experiment.
We run the following script to produce $N=10^4$ samples from $f$. The contour plot of $f$ and the  samples are displayed in
Figure \ref{fig:mcmc_randwalk1}. We discarded an initial ``burn-in'' period (say the first $1000$ samples) to ensure the chain has reached a stable distribution. 
We observe that the correct region is sampled.

\begin{shaded}
\vspace*{-0.3cm}
\begin{verbatim}
def f(xy):         # Define the distribution 
    x,y = xy[0];  c = 1/20216.335877
    return c*np.exp(-(x**2*y**2+x**2+y**2-8*x-8*y)/2)
Sigma = 2*np.eye(2)
N = 10**4          # number of MH samples
X0 = [0,0]         # initial guess
X = McMcRandWalkGen(f, X0, Sigma, N)  # Call the MCMC function 
#      plot the results 
xeval = np.linspace(-1, 7, 1000)
[x,y] = np.meshgrid(xeval,xeval)
feval = f([[x,y]])
plt.figure(figsize = (5, 5))
plt.contour(x,y,feval)
plt.figure(figsize = (5, 5))
plt.plot(X[0,1000:], X[1,1000:], color = 'red',
         marker = 'o', markersize = 2, linestyle = '')
\end{verbatim}
\vspace*{-0.3cm}
\end{shaded}

 \begin{center}
 \includegraphics[scale=0.6]{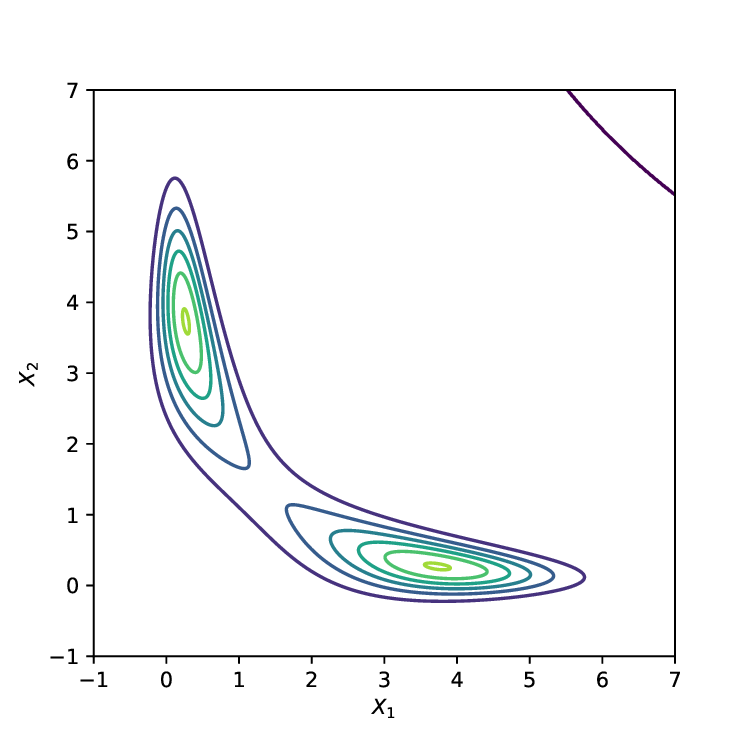}\includegraphics[scale=0.6]{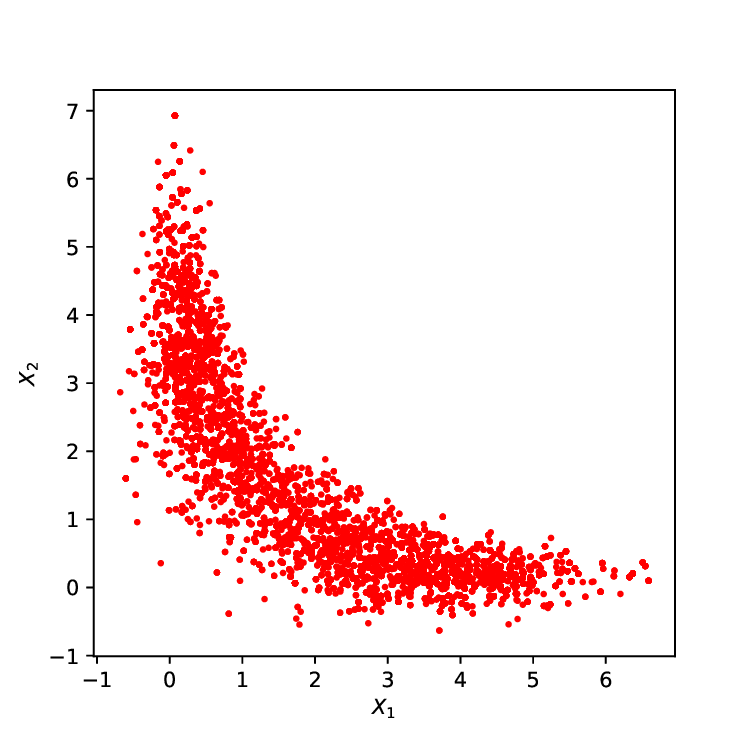}
 \captionof{figure}{The contour plot of bivariate distribution $f$ (left) and the samples using the Metropolis-Hastings algorithm (right). The first 1000 samples are discarded.}\label{fig:mcmc_randwalk1}
 \end{center}

The histogram of variable $X_1$ is also shown in Figure \ref{fig:mcmc_randwalk2}. It is close to the true marginal pdf. We use a numerical integration to compute the exact marginal distribution via formula
$$
f_{X_1}(x) = \int_{-\infty}^{\infty} f(x,y) dy.
$$

 \begin{center}
 \includegraphics[scale=0.7]{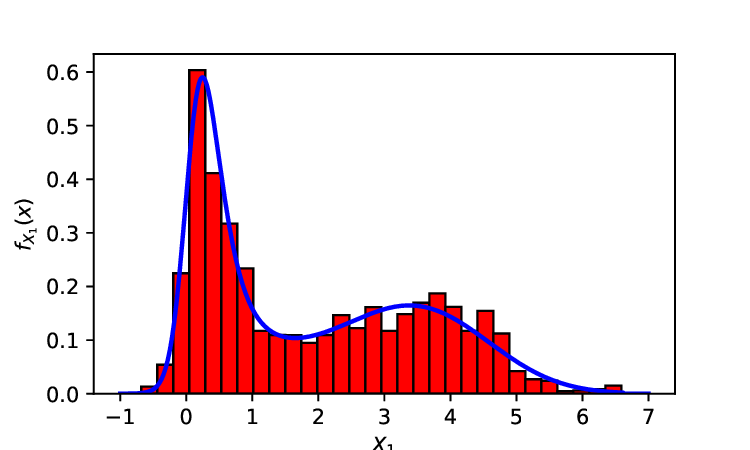}
 \captionof{figure}{The histogram of the marginal variable $X_1$, and the true corresponding marginal pdf (the blue curve).}\label{fig:mcmc_randwalk2}
 \end{center}

Now, assume that we want to estimate $\E_f(X_1)$ using the Monte Carlo method. 
It is enough to compute the mean of marginal samples $\{X_{t1}\}$. 
\begin{shaded}
\vspace*{-0.3cm}
\begin{verbatim}
MeanVals = np.zeros(4)
for k in range(4):
    N = 10**(k+3)
    X0 = [0,0]
    X = McMcRandWalkGen(f, X0, SigmaWalk, N)
    MeanVals[k] = np.mean(X[0,1000:])
print("MCMC estimates = ", np.round(MeanVals,4))
\end{verbatim}
\vspace*{-0.3cm}
\end{shaded}
The output of a run is
\begin{shaded}
\vspace*{-0.3cm}
\begin{verbatim}
MCMC estimates =  [1.6999 1.8169 1.8985 1.8584]
\end{verbatim}
\vspace*{-0.3cm}
\end{shaded}
which actually four estimations for $\E_f(X_1)$ with $N = 10^3,10^4,10^5$ and $10^6$. The exact value is $\E_f(X_1)\doteq 1.85997$.
\end{example}
\vsp 

\subsection{MCMC Bayesian parameter estimation}
One application of the Metropolis-Hastings algorithm arises in Bayesian inference for parameter estimation.
In Bayesian inference, we aim to estimate parameters 
\( \theta = (\theta^{(1)}, \theta^{(2)}, \ldots, \theta^{(d)}) \)
 of a model by calculating the {\em posterior distribution} \( p(\theta | \text{data}) \), which reflects our updated beliefs about the parameters after observing the data. The posterior distribution is typically computed as
\[
p(\theta | \text{data}) = p(\text{data} | \theta) \cdot p(\theta)
\vsp
\]
where \( p(\text{data} | \theta) \) is the {\em likelihood} of observing the data given the parameters, and \( p(\theta) \) is the {\em prior distribution}, which encodes our initial beliefs about \( \theta \) before seeing the data.

However, calculating \( p(\theta | \text{data}) \) explicitly can be difficult, especially if the likelihood is complex or the dimensionality of \( \theta \) is high. This is where the Metropolis-Hastings algorithm becomes useful. We assume that the proposal distribution $q$ is given and 
$f(\theta) = p(\text{data} | \theta) \cdot p(\theta)$.
We start with an initial guess \( \theta_0 = (\theta^{(1)}_0, \ldots, \theta^{(d)}_0) \) for the parameters and we use Algorithm \ref{alg:metropolis-hastings} to generate a chain of parameter samples.
After a sufficient number of iterations the generated samples \( \{ \theta_1, \theta_2, \ldots, \theta_N \} \) approximate the target posterior distribution \( p(\theta | \text{data}) \). 
We can now use the Metropolis-Hastings samples to
estimate the {\em posterior} mean, variance, and {confidence intervals} for each parameter.
\vsp

\begin{example}[Estimating Patient Recovery Rate]\label{ex:mcmc-patient1}
Suppose we observe the recovery times of $10$ patients after treatment, given in days as
\[
\text{data} = \{5, 8, 12, 7, 9, 10, 3, 6, 8, 11\}.
\]
Our goal is to estimate the posterior distribution of the recovery rate parameter \( \theta \), which represents the average rate at which patients recover.
We assume that recovery times follow an exponential distribution
  \[
  p(x | \theta) = \theta e^{-\theta x}.
  \]
The aim is to estimate the parameter $\theta$ by approximating its (posterior) distribution, i.e., generating samples from its distribution that is unknown to us. But we assume that \( \theta \) has a prior distribution. In this example we will use the {\em Gamma} distribution as a prior for \( \theta \), i.e.,
  \[
  p(\theta) = \frac{\beta^\alpha}{\Gamma(\alpha)} \theta^{\alpha - 1} e^{-\beta \theta}
  \vsp
  \]
where \( \alpha = 2 \) and \( \beta = 1 \), which encodes our prior belief that \( \theta \) is likely around 2 because the mean of the Gamma distribution is $\alpha/\beta$.

For the given independent recovery times \( x_j \) (our data), the likelihood of observing the data given \( \theta \) is
   \[
   p(\text{data} | \theta) = \prod_{j=1}^{10} \theta e^{-\theta x_j}.\vsp
   \]
Note that, in the likelihood $\theta$ is the parameter of the distribution while in the prior $\theta$ is the variable of the distribution. 

We start with an initial guess \( \theta_0= 1 \) and use the random walk sampler (normal distribution $\Nor(\theta_t,\sigma^2)$) with a small variance $\sigma^2 = 0.25$ at each time step $t$ to generate the new parameter $\theta_{t+1}$ using the Metropolis-Hastings algorithm. The code is given below. Note that to define the distribution $f(\theta) = p(\text{data} | \theta) \cdot p(\theta)$ we use the logarithm of likelihood and prior distributions to compute $\log(f(\theta))$ and finally return $f(\theta)=\exp(\log(f(\theta)))$. 

\begin{shaded}
\vspace*{-0.3cm}
\begin{verbatim}
import numpy as np
import matplotlib.pyplot as plt
from scipy.stats import gamma, expon   # Gamma and Exponential distributions

# Given data
data = [5,8,12,7,9,10,3,6,8,11]
# Define MCMC distribution = likelihood * priors
def mcmc_pdf(theta):
    log_likelihood = np.sum(expon.logpdf(data, scale = 1/theta))
    alpha = 2; beta = 1; 
    log_prior = gamma.logpdf(theta, alpha, scale = 1/beta)
    return np.exp(log_likelihood+log_prior)

Sigma = 0.5   # std value for random walk sampler (transient kernel)
theta0 = .5   # initial guess
N = 10**4     # length of MCMC chain
BurnIn = 1000 # burn-in period for MH

Theta = McMcRandWalkGen(mcmc_pdf, theta0, Sigma, N)  # Call MCMC algorithm
Theta = Theta[:,BurnIn:]          # Discard an initial burn-in period
MeanTheta = np.mean(Theta) 

err = 1.96*np.std(Theta)/np.sqrt(np.size(Theta))
print("Posterior mean of parameter is ${0} +/- {1} with 95% of probability"
      .format(np.round(MeanTheta,4), np.round(err,4)))
      
plt.hist(Theta[0,:], bins=30, density=True, alpha=0.6, color='red')
plt.axvline(MeanTheta, color='b', linestyle='dashed', linewidth=1.5, 
                      label=f"Mean recovery rate: {MeanTheta:.4f}")
plt.title("Distribution of recovery rate parameter")
plt.xlabel("Recovery rate")
plt.ylabel("Density")
plt.legend()
plt.show()
\end{verbatim}
\vspace*{-0.3cm}
\end{shaded}

The histogram of $\theta$ samples is plotted and together with its estimated mean value are shown in Figure 
\ref{fig:hist_bayesian1}. 

\begin{center}
\includegraphics[scale=.6]{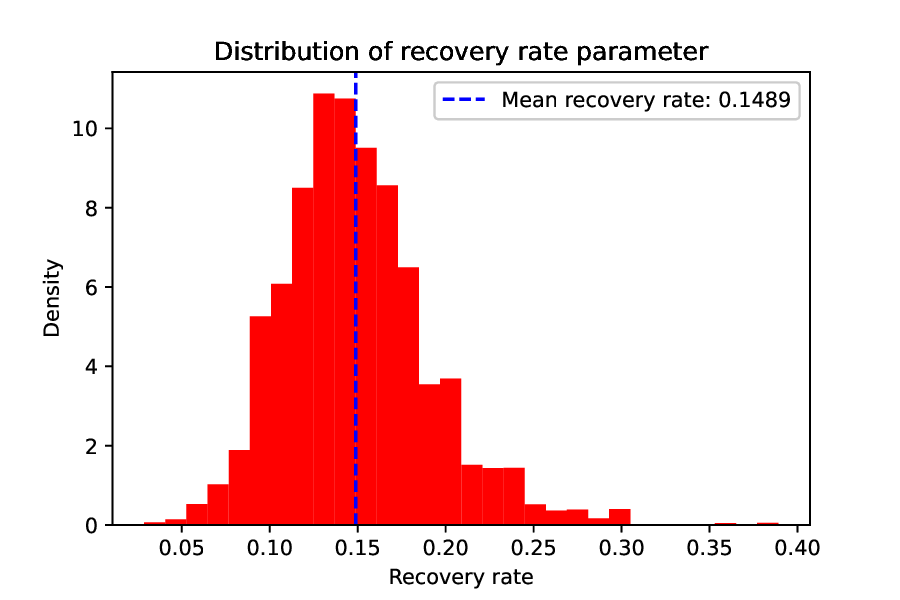}
\captionof{figure}{Posterior distribution of parameter $\theta$ sampled by the Metropolis-Hastings algorithm, and its estimated mean value.}\label{fig:hist_bayesian1}
\end{center}

We also calculated the standard deviation  of \( \theta \) values and a confidence interval with $95\%$ probability to get an estimate for the error of the recovery rate. The output is 
\begin{shaded}
\vspace*{-0.3cm}
\begin{verbatim}
Posterior mean of parameter is $0.1489 +/- 0.0009 with 95% of probability
\end{verbatim}
\vspace*{-0.3cm}
\end{shaded}

\end{example}
\vsp

\begin{example}[Estimating the probability of a large portfolio loss]

A financial analyst wants to estimate the likelihood that an investment portfolio valued at $\$1,000,000$ will lose more than $\$100,000$ over a 6-month period. This probability, often referred to as the {\em Value at Risk (VaR)} at a certain confidence level, is crucial for managing risk and setting capital reserves.

We assume the portfolio value $S_t$ follows a geometric Brownian motion 
  \[
  \td S_t = \mu S_t \, \td t + \sigma S_t \, \td W_t
  \]
 which is commonly used to model stock prices, as discussed in Example \ref{ex:stock_value}.
Here  \( \mu \) is the drift rate (expected return rate of the portfolio),
 \( \sigma \) is the volatility of returns, and
  \( W_t \) is a Wiener process. 
We learnt in Section \ref{sect:gauss_process_gen} how to generate a sample path $S_t$ using the Euler-Maruyama method. 

The goal here is to estimate the posterior distributions of parameters \(\theta = (\mu , \sigma) \) given historical data on portfolio returns. From these, then we can simulate future values of \( S_t \) and estimate the probability of a loss greater than $\$100,000$ over the next six months.

The analyst has monthly return rate data for the portfolio in the past five years, which we denote by
\[
\text{data} = \{r_1, r_2, \ldots, r_{60}\} \vsp
\]
where each \( r_j \) is the observed monthly return rate.
Assuming that monthly returns are normally distributed, we define the likelihood of observing a monthly return \( r_j \) given \( \mu \) and \( \sigma \) as
   \[
   p(r_j | \mu, \sigma) = \frac{1}{\sqrt{2 \pi \sigma^2}} \exp \left(-\frac{(r_j - \mu)^2}{2 \sigma^2}\right). \vsp
   \]
   For the entire dataset, the likelihood is
\[
   p(\text{data} | \mu, \sigma) = \prod_{j=1}^{60} p(r_j | \mu, \sigma).
\]

It is remained to specify prior distributions for \( \mu \) and \( \sigma \). 
For \( \mu \), a normal distribution centered at the historical average return, say \( \tilde \mu \), with standard deviation \( \tilde\sigma \) is used, i.e.,
     \[
     p_1(\mu) \sim \Nor(\tilde\mu, \tilde \sigma^2).
     \]
For \( \sigma \), an Inverse-Gamma distribution
     \[
     p_2(\sigma^2) \sim
     \mathcal{I}nv\mathcal{G}am(\alpha, \beta)
     \]
is used. The Inverse-Gamma distribution 
 is commonly used for variance and precision parameters in Bayesian statistics, and has the distribution 
$$
f(x) = \frac{\beta^\alpha}{\Gamma(\alpha)} (1/x)^{\alpha+1}\exp({-\beta /x}),\quad x >0, \quad \alpha,\beta >0.
$$
This distribution has parameters \( \alpha \) and \( \beta \) that can be chosen based on prior knowledge of the portfolio's volatility. Note that for $\alpha>1$ the mean of the Inverse-Gamma distribution is $\beta/(\alpha-1)$. 
The joint prior distribution is defined as
\[
p(\mu,\sigma^2) = p_1(\mu)\cdot p_2(\sigma^2). 
\]  
To apply the Metropolis-Hastings algorithm
we start with initial guesses for \( \mu \) and \( \sigma \), say \(\theta_0 =  (\mu_0 ,\sigma_0)\), and use a random walk sampler with covariance matrix 
$$
\Sigma = \begin{bmatrix}
\tau_\mu^2 & 0 \\ 0 & \tau_\sigma^2
\end{bmatrix}
$$
to propose a new parameter vector $\theta_{t+1} = (\mu_{t+1}, \sigma_{t+1})$
at each step $t=0,1,\ldots$ of the algorithm.
Here \( \tau_\mu \) and \( \tau_\sigma \) are small tuning parameters to control the step size of the sampler.
We run the Metropolis-Hastings algorithm for many iterations (e.g., \(N_{\text{MH}} = 10^4 \)), and discard the first few thousand samples for burn-in.

Finally, to estimate the VaR, we do the following steps:
\begin{enumerate}
\item (Simulate future portfolio values):
   Using the posterior samples of \( \mu \) and \( \sigma \), simulate the portfolio value after 6 months under the geometric Brownian motion model. For each sampled \( (\mu, \sigma) \) pair, generate $N_{\text{MC}}$ paths using 
function \texttt{DiffusionProcessGen} for drift coefficient $a(t,S_t) = \mu S_t$ and diffusion coefficient $b(t,S_t) = \sigma S_t $,  
   and calculate the final portfolio values \( S_{T} \) for $T = 0.5$ a year (6 months). Finally, take a mean to have one estimate for the portfolio value for each pair $(\mu,\sigma)$.

\item (Calculate loss probability):
   Compute the proportion of paths where the portfolio value is below $\$900,000$ (a $\$100,000$ loss from the starting value). This proportion estimates the probability of a large loss in the period.
\end{enumerate}

The Python code is given below where we use the input values as followings:
The initial portfolio value is \( S_0 = 1,000,000 \), the loss threshold is \( 900,000 \) (a $\$100,000$ loss), the 
time frame is 6 months, 
the historical mean return is \( \tilde \mu = 0.05 \) ($5\%$ per month), 
the historical volatility is \( \tilde\sigma = 0.1 \), 
parameters for Inverse-Gamma distribution are $\alpha = 2$ and $\beta = 0.0004$, 
the proposal variances for \( \mu \) and \( \sigma \) are \( \tau_\mu = 0.001 \) and \( \tau_\sigma = 0.001 \), 
the initial guess is $\theta_0=(\mu_0,\sigma_0)= (\tilde \mu,\tilde \sigma)$, 
the number of MCMC iterations is \( N_{\text{MH}} = 10^4 \), and 
the number of Monte Carlo paths per posterior sample is \(N_{\text{MC}} = 10^3 \).
The data values $r_j$ are given in the code. 

\begin{shaded}
\vspace*{-0.3cm}
\begin{verbatim}
import numpy as np
import matplotlib.pyplot as plt
from scipy.stats import invgamma, norm  # import Inv-Gamma and normal dists.

# Parameters
S0 = 1000000               # Initial portfolio value
LossThreshold = 900000     # Loss threshold (100,000 loss)
FinalTime = 0.5            # Time horizon in years (6 months)
N_MH = 10000               # Number of Metropolis-Hastings iterations
N_MC = 1000                # Monte Carlo paths for each posterior sample
BurnIn = 2000              # Burn-in period for MH

# Prior parameters
mu_prior = 0.05            # Historical mean return
sigma_prior = 0.1          # Historical volatility
alpha_prior, beta_prior = 2, 0.0004

# Metropolis-Hastings proposal variances
tau_mu = 0.001
tau_sigma = 0.001

# Given data for sixty months
data = np.array(
  [0.07,  0.13,  0.10,  0.17,  0.11,  0.03,  0.15,  0.09,  0.12,  0.12,
  -0.06,  0.07,  0.09, -0.01,  0.08,  0.08,  0.07,  0.19,  0.09,  0.12,  
   0.03,  0.16, -0.02,  0.2 ,  0.14,  0.05,  0.08,  0.06,  0.10, -0.07,
  -0.01, -0.07, -0.05,  0.21, -0.05,  0.02, -0.02,  0.15,  0.08,  0.02, 
  -0.03,  0.01,  0.08,  0.13,   0.16, -0.03, -0.13, 0.14,  0.11,  0.12, 
  -0.01, -0.07,  0.16,  0.27, -0.06,   0.01,  0.01, 0.01,  0.01,  0.16])     
\end{verbatim}
%\vspace*{-0.3cm}
\end{shaded}
\begin{shaded}
%\vspace*{-0.3cm}
\begin{verbatim}
# Define MCMC distribution = likelihood * priors
def mcmc_pdf(theta):
    mu, sigma = theta[0]
    log_likelihood = np.sum(norm.logpdf(data, loc=mu, scale=sigma))
    log_prior_mu = norm.logpdf(mu, loc=mu_prior, scale=sigma_prior)
    log_prior_sigma = invgamma.logpdf(sigma**2, a=alpha_prior, scale=beta_prior)
    return np.exp(log_likelihood + log_prior_mu + log_prior_sigma)

# Call the Matropolis-Hastings algorithm
theta0 = [mu_prior, sigma_prior]   # initial guess theta0 = [mu0,sigma0]
SigmaWalk = np.diag([tau_mu,tau_sigma])     # covariance of the sampler
Theta = McMcRandWalkGen(mcmc_pdf, theta0, SigmaWalk, N_MH)

Theta = Theta[:,BurnIn:]   # Discard burn-in samples

# Monte Carlo simulation for each posterior sample (mu, sigma)
def drift(t,x,*args):      # drift
    mu = args[0]
    return mu*x
def diffusion(t,x,*args):  # diffusion
    sigma = args[1]
    return sigma*x
TimeStep = 0.01         # Time step in Euler-Maruyama method
LossPr = []             # Loss probability vector
ST = np.empty(N_MC)     # Stock value vector at final time
for mu, sigma in zip(Theta[0,:], Theta[1,:]):  # loop over all samples
    for j in range(N_MC):                      # Monte Calro loop
        ST[j] = DiffusionProcessGen(drift, diffusion, [0,FinalTime], 
                                    TimeStep, S0, mu, sigma)[-1]
    LossPr.append(np.mean(ST < LossThreshold))

# Calculate the mean loss probability
MeanLossPr = np.mean(LossPr)

# Results
print(f"Estimated Probability of a $100,000 Loss: {MeanLossPr:.4f}")

# Plotting the results
\end{verbatim}
%\vspace*{-0.3cm}
\end{shaded}
\begin{shaded}
%\vspace*{-0.3cm}
\begin{verbatim}
plt.figure(figsize=(6, 4))
plt.hist(LossPr, bins=30, density=True, alpha=0.6, color='red')
plt.axvline(MeanLossPr, color='b', linestyle='dashed', linewidth=1.5,
                        label=f"Mean Loss Probability: {MeanLossPr:.4f}")
plt.title('Distribution of Loss Probabilities from MC Simulations')
plt.xlabel('Loss Probability')
plt.ylabel('Density')
plt.legend()
plt.show()
\end{verbatim}
\vspace*{-0.3cm}
\end{shaded}

We computed the mean probability of loss from all simulations and visualize the distribution of the estimated loss probabilities by drawing the histogram of the loss probability values. See Figure \ref{fig:loss_prob}. The estimated probability of a $\$100,000$ loss is about 
$1\%$. 
\begin{center}
\includegraphics[scale=.6]{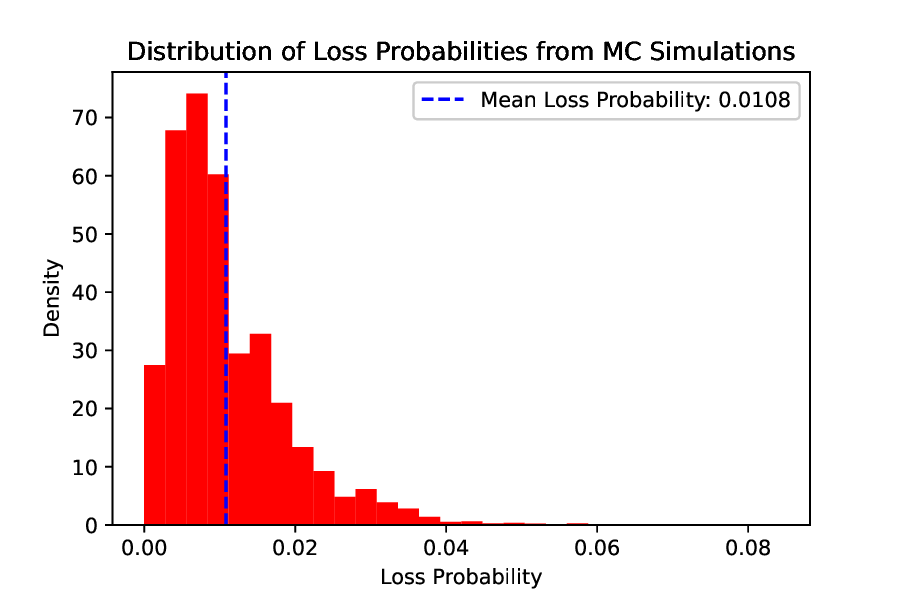}
\captionof{figure}{The distribution of the estimated loss probabilities and its estimated mean.}\label{fig:loss_prob}
\end{center}

\end{example}
\vsp 
\vsp

The main advantage of MCMC is that it can be used to generate random samples from any target distribution, regardless of its dimensionality and complexity. Main disadvantages are
(1) The resulting samples are often highly correlated, (2) Sometimes $N$ should be large so that the Markov
chain settles down to its steady state, (3) the estimates obtained by MCMC often have greater
variances than those obtained from iid samplings. During the past years,  modified versions and more efficient MCMC algorithms have been developed to overcome such disadvantages. See for example \cite{Robert-Casella:2004}.

%\subsection{Gibbs sampler algorithm}

%% file: workouts_4.tex
\section{Exercises}

This section includes a number of exercises designed to deepen your understanding of the lecture material. Some of the exercises are adapted from the references listed in the bibliography.
\vsp

\begin{workout}
The {\em Weibull distribution} is named after Swedish mathematician Waloddi Weibull, who described it in detail in 1939.
The pdf (probability density function) of this distribution is defined by two parameters $\lambda$ (scale parameter) and $\alpha$ (shape parameter), and is given by
$$
f(x) = \frac{\alpha}{\lambda} \left( \frac{x}{\lambda}\right)^{\alpha-1}e^{-\left({x}/{\lambda}\right)^\alpha}, \quad x\geq 0. \vsp 
$$
Demonstrate how the Inverse Transform Method (ITM) can be used to generate random numbers from a Weibull distribution. Write down all the steps. 
Then design a Monte Carlo algorithm for estimating the value of integral 
$$
I = \int_0^\infty \sqrt{1+x^2} f(x)\, dx 
$$
where $f(x)$ is a Weibull pdf with parameters $\alpha=2.5$ and $\lambda=1$. Be sure to not only return the estimation of $I$ but also an estimate of the error.
\end{workout}
\vsp

\begin{workout}
Consider the integral 
$$
\displaystyle I = \int_{0}^\pi (x+\sqrt x) \sin x \, \text dx. \vsp 
$$
Let us consider $I$ as an expectation value of a random variable $g(X)$ where $X$ has a sine-distribution with pdf
$$
f(x) = \begin{cases} \frac{1}{2} \sin(x),& 0\leq x\leq \pi\\ 0 , & otherwise \end{cases}.\vsp 
$$
We write $X\sim \sin$. 
 Given uniformly distributed random numbers $U\sim \Uni(0,1)$, write with details that how you convert these into sine distributed random numbers.
Then write $I$ as an expectation value using such a random variable $X\sim \sin$ and write a pseudocode that estimates the integral using a Monte Carlo method with a given number of $N$ samples. Be sure to not only return estimation of $I$ but also an estimate of the error.
\end{workout}
\vsp 

\begin{workout}
Let $X_1,\ldots, X_n$ be iid random variables with cdf $F$.
Assume that
$$
X_{(1)} := \min\{X_1,\ldots,X_n\}, \quad X_{(n)}:=\max\{X_1,\ldots,X_n\}.
$$
First prove that the cdf of $X_{(n)}$ is $F_n(x) = [F(x)]^n$ and the cdf of $X_{(1)}$ is $F_1(x)=1-[1-F(x)]^n$. Then show that
$$
X_{(n)} = F^{-1}(U^{1/n}), \quad X_{(1)} = F^{-1}(1-U^{1/n})
$$
where $U\sim \Uni(0,1)$. Random variables $X_{(1)}$ and $X_{(n)}$ are called {\em ordered statistics}.
Show how inverse transform method can be used to sample from  ordered statistics. 

\end{workout}
\vsp

\begin{workout}
How the inverse transform method can be applied to generate from Beta distributions ${\Bet}(\alpha,1)$ and 
${\Bet}(1,\beta)$. Derive the formulation and implement the Python code.   
\end{workout}
\vsp 

\begin{workout}
Explain how we can generate a random variable $X$ from the semicircular pdf
$$
f(x) = \frac{2}{\pi R^2}\sqrt{R^2-x^2}, \quad x\in [-R,R]
$$
using the acceptance-rejection algorithm. Then use the \verb+RandAcceptReject+ function for illustration.

\end{workout}
\vsp 

\begin{workout}%[mini-project]
Write a Python code for drawing normal samples using the acceptance-rejection method with exponential distribution as a proposal distribution. Refer to discussions at the end of subsection \ref{sect:accept-reject}.
Plot histograms for different number of samples.
\end{workout}
\vsp

\begin{workout}
Develop an algorithm for sampling from $\Bin(p,n)$ for large values of $n$ using the fact that the distribution of a binomial variable
$X\sim {\Bin}(p,n)$ is close to that of $Y\sim \Nor(np-1/2,np(1-np))$ for a large $n$.
\end{workout} 
\vsp 
\begin{workout} 
[Monty Hall Problem]
You are on a game show, being asked to choose between three doors. One door
has a car, and the other two have goats. The host, Monty Hall, opens one of the other doors, which he knows has a goat behind it. Monty then asks whether you would like to switch your choice of door to the other remaining door. Do you choose to switch or not to switch?
Solve it with Monte Carlo method. 

Answer: If you switch you win the car with $2/3$ probability. Image source: Wikipedia. 

\begin{center}
\includegraphics[scale=.15]{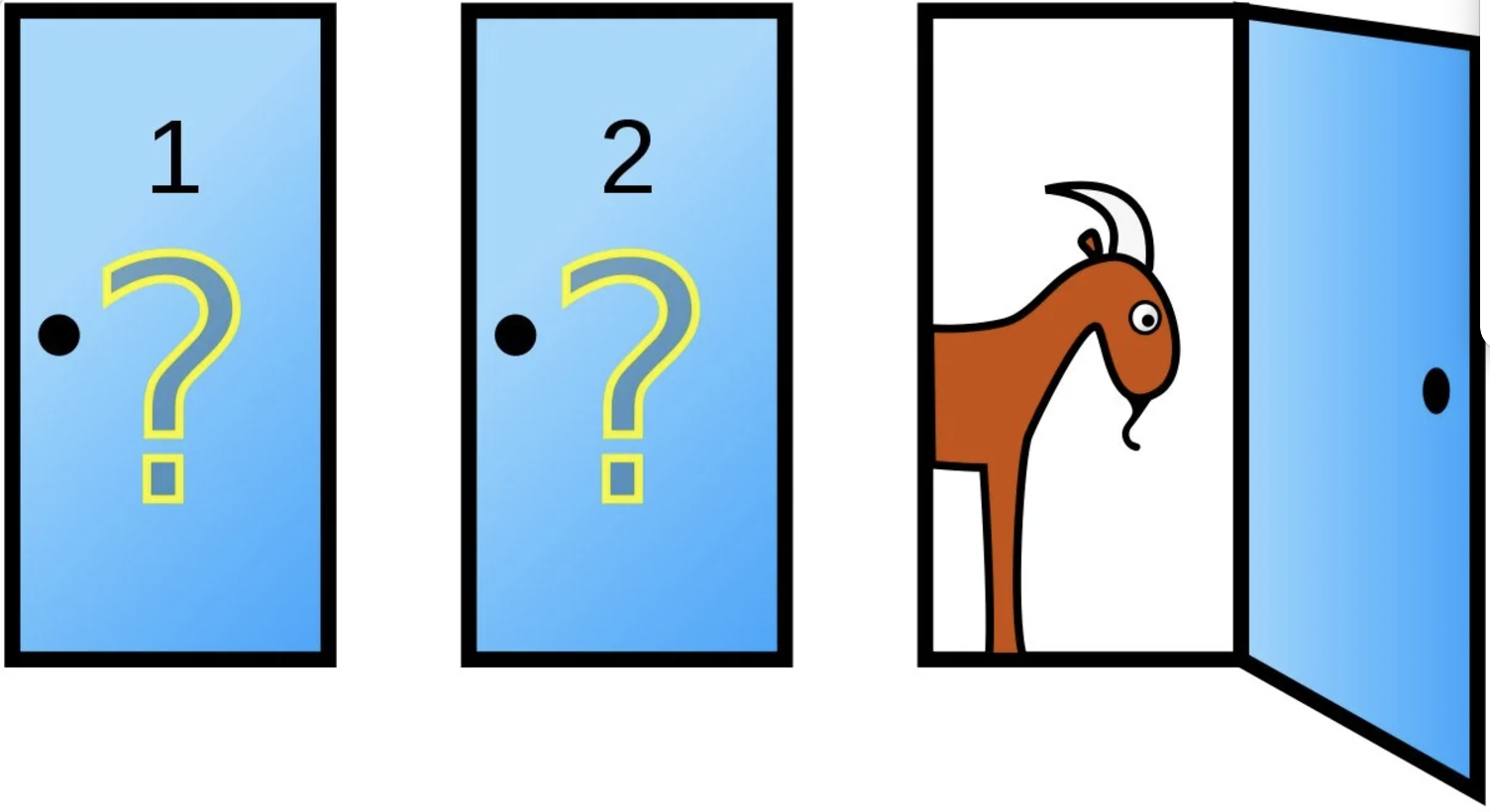}
\end{center}
\end{workout}
\vsp

\begin{workout}
For the normal-Cauchy Bayes estimator
$$
\delta(t) = \int_{-\infty}^{\infty}\frac{x}{1+x^2}\ee^{-(x-t)^2/2}\td x \Big/ \int_{-\infty}^{\infty}\frac{1}{1+x^2}\ee^{-(x-t)^2/2}\td x
$$
use the Monte Carlo integration based on normal simulations to estimate $\delta(t)$ for $t=0,2,4$. Monitor the convergence with the standard error of the estimate. Determine the minimum value $N$ to obtain three digits of accuracy with $0.95\%$ probability.
\end{workout}
\vsp

\begin{workout}
Write down the proof of Theorem \ref{thm-markov0}.
\end{workout}
\vsp

\begin{workout}
Suppose that the morning weather of a city in a time period can be only sunny or cloudy,
and the weather conditions on successive mornings form
a Markov chain with transition  matrix
$$
\begin{array}{rcc}
& \mathrm{sunny}    & \mathrm{cloudy}  \\
  \mathrm{sunny}  & \big\lceil 0.7 & 0.3\big\rceil \\
  \mathrm{cloudy} & \big\lfloor0.6 & 0.4\big\rfloor
\end{array}.
$$
\begin{itemize}
\item[1.] If it is cloudy on a given day, what is the probability
that it will also be cloudy the next day?
\item[2.] If it is sunny on a given day, what is the probability
that it will be sunny on the next two days?
\item[3.] If it is cloudy on a given day, what is the probability
that it will be sunny on at least one of the next three
days?
\item[4.] If it is sunny on a certain Wednesday, what is the
probability that it will be sunny on the following
Saturday?
\item[5.] If it is cloudy on a certain Wednesday, what is the
probability that it will be sunny on the following
Saturday?
\item[6.] If it is sunny on a certain Wednesday, what is the
probability that it will be sunny on both the following
Saturday and Sunday?
\item[7.] If it is cloudy on a certain Wednesday, what is the
probability that it will be sunny on both the following
Saturday and Sunday?
\item[8.] Suppose that the probability that it will be sunny
on a certain Wednesday is $0.2$ and the probability that it
will be cloudy is $0.8$.
Determine the probability that it will be cloudy on
the next day, Thursday.
\item[9.] With assumptions of item 8, determine the probability that it will be cloudy on
Friday.
\end{itemize}

\end{workout}
\vsp 
\begin{workout}
Suppose that a Markov chain has state space $\Ss$ with four states $\{1, 2, 3, 4\}$
and transition matrix
$$
\begin{array}{rcccc}
         & 1    & 2 & 3& 4 \\
  1  & \big\lceil 0.25 & 0.25 & 0.00 & 0.50 \big\rceil \\
  2  & \big|      0.00 & 1.00 & 0.00 & 0.00 \big| \\
  3  & \big|      0.50 & 0.00 & 0.50 & 0.00 \big| \\
  4  & \big\lfloor0.25 & 0.25 & 0.25 & 0.25 \big\rfloor
\end{array}\vsp
$$
If the chain is in state $3$ at a given time $t$, what is the
probability that it will be in state $2$ at time $t + 2$?
If the chain is in state $1$ at a given time $t$, what is the
probability that it will be in state $3$ at time $t + 3$?
\end{workout}
\vsp 

\begin{workout}
We want to use a finite Markov chain to model the probability of customers making purchases based on their past behavior.
Assume that customers can be in one of the following four states:
\begin{enumerate}
\item {\em Browsing}: The customer is browsing the website without purchasing.
\item  {\em Added to Cart}: The customer has added items to their cart but has not checked out.
\item {\em Checkout}: The customer is in the checkout process but has not completed the purchase.
\item {\em Purchased}: The customer has completed a purchase.
\end{enumerate}
We have observed customer behavior and estimated the following transition probabilities for moving from one state to another.
\begin{itemize}
\item[-] Customers in the {\em Browsing} state stay there with $60\%$ probability, move to {\em Added to Cart} with $30\%$ probability, and go directly to {\em Checkout} with $10\%$ probability.
\item[-] Customers who have added items to their cart have a $20\%$ chance of going back to {\em Browsing}, a $50\%$ chance of staying in {\em Added to Cart}, a $20\%$ chance of moving to {\em Checkout}, and a $10\%$ chance of making a {\em Purchase}.
\item[-] Customers in {\em Checkout} have a $60\%$ chance of staying in {\em Checkout}, a $30\%$ chance of making a {\em Purchase}, and a $10\%$ chance of returning to {\em Added to Cart}.
\item[-] Once in the {\em Purchased} state, customers stay there permanently ($100\%$ probability), as the process ends with a purchase.
\end{itemize}
%\[
%P = \begin{bmatrix}
%0.6 & 0.3 & 0.1 & 0.0 \\
%0.2 & 0.5 & 0.2 & 0.1 \\
%0.0 & 0.1 & 0.6 & 0.3 \\
%0.0 & 0.0 & 0.0 & 1.0 \\
%\end{bmatrix}
%\]

The goal is to estimate the probability that a customer will complete a purchase within 6 steps, starting from the {\em Browsing} state. 
Write down the transition matrix and use the \texttt{MarkovChainGen} function in a Monte Carlo loop to estimate the the purchase probability in 6 steps. 
Write a Python code and report the result. 
\end{workout}
\vsp 

\begin{workout}
You are running a startup company. Let $V_t$ denote the valuation of your company at time $t$, $t=0,1,\ldots$  (say, in months). If $V_t$, then the company goes bankrupt and stops operating; if $V_t = v_{\text{max}}$, then the company is acquired by a larger company, you receive a payout $v_{\text{max}}$, and the company stops operating. In each time period that the company operates, you invest additional fixed amount $c_{\text{invest}}$ in the company, and you also incur an operating cost $c_{\text{operate}}$.
If $0<V_t<v_{\text{max}}$  then
$$
V_{t+1} = \begin{cases} V_t + \delta, & \text{with probability } p\\ V_t - \delta, & \text{with probability } 1- p \end{cases}
$$                             
Here $\delta>0$ is a given parameter. The initial valuation $V_0$ is an integer multiple of $\delta$, as is $v_{\text{max}}$, so all $V_t$ are also integer multiples of $\delta$. With this model, you will eventually either go bankrupt or be acquired.
Consider this model with the following parameter instances: 
\begin{align*} 
&\delta = 2M \text{ sek}, \; \; \;V_0 = 10M \text{ sek}, \; \;\; v_\text{max} = 100M \text{ sek}, \\ &c_\text{operate} = 10K \text{ sek}, \; \;\; c_{\text{invest}} = 200K \text{ sek}, \;\; p = 0.6 
\end{align*}
First, identify the type of random process  and explain why. Then design a Monte-Carlo algorithm (write a pseudo-code) to estimate: 
\begin{itemize}
 \item[   a) ] the probability that the startup goes bankrupt, 

\item[    b) ] the expected time until the startup goes bankrupt or is acquired,

\item[    c) ] the expected profit, if the startup is acquired,

\item[    d) ] the expected loss, if the startup goes bankrupt.
\end{itemize}
Finally, implement your algorithm in a Python environment  (with $N = 5000$ Monte Carlo  simulations) and report all the above estimations.

Note that, {\em profit} is the payout, when the company is acquired, minus the initial value of company, minus the total operating cost, minus the total of any investments made.
{\em Loss} is the initial value of company plus the total operating cost plus the total of any investments made, when the company goes bankrupt. 
\end{workout}
\vsp 

\begin{workout}
Apply the Gillespie algorithm to solve the predator-prey model
\begin{align*}
     R &\xrightarrow{\;\alpha\;} 2R \\
    R+F &\xrightarrow{\;\beta\;} 2F \\
    F &\xrightarrow{\;\gamma\;} \emptyset
\end{align*}
with rates $\alpha = 1$, $\beta=0.005$, and  $\gamma = 0.6$ and initial value $(F_0,R_0) = [50,100]$ up to final time $t_{final} = 30$. There exists also the deterministic model 
\begin{align*}
    &\frac{dF}{dt} = \beta FR - \gamma F \\
    &\frac{dR}{dt} = \alpha R -\beta FR
\end{align*}
with $R(0)=R_0$ and $F(0)=F_0$ associated to this simple ecology. 
Use the ODE solver \texttt{solve\_ivp} from the \texttt{scipy.integrate} library to solve this ODE with the same input data and compare the results of stochastic and deterministic models.

\end{workout}
\vsp 

\begin{workout}

Consider the Michaelis-Menten model
\begin{center}
\begin{tikzpicture}[>=latex,scale=1]
\node at (0.25,0) {$S+E$}; 
\draw[-to,black,thick] (1,.05) -- (1.7,.05);
\draw[to-,black,thick] (1,-.05) -- (1.7,-.05);
\node at (2,0) {$C$};

\draw[-to,black,thick] (2.25,0) -- (2.95,0);
\node at (3.6,0) {$P+E$};
\node at (1.35,.25) {$c_1$};
\node at (1.35,-.25) {$c_2$};
\node at (2.6,.2) {$c_3$};
\end{tikzpicture}
\end{center}
which is the standard model for enzyme catalysts. 
Here, $c_1$ (forward rate), $c_2$ (reverse rate), and $c_3$ (catalytic rate) denote the constant rates of the reactions.  Our intention is to solve this model using the Gillespie's algorithm (SSA).

Write down propensity functions and state-change vectors for this model. Assume that for a certain enzyme the reactions rates are $c_1 = 0.002~ \text{mol}^{-1}\text{sec}^{-1}$, $c_2 = 0.1~ \text{sec}^{-1}$ and $c_3 = 0.75~ \text{sec}^{-1}$ where $\text{sec}$ stands for the unit of time and $\text{mol}$ is the unit for number of proteins. Furthermore, assume that at time $t=0.1,\, \text{sec}$ the number of proteins have been computed as $E(t) = 300~\text{mol}$, $S(t) = 200~\text{mol}$, $C(t) = 100~\text{mol}$, and $P(t) = 50~\text{mol}$. The task is to compute the number of proteins in the next time level $t+\tau$. To this aim, we have generated two uniform random numbers $u_1=0.64$ and $u_2=0.83$ from the $\Uni(0,1)$ distribution. The number $u_1$ must be used to determine the steplength $\tau$, and $u_2$ to determine the specific reaction that will occur. Given these conditions, proceed to compute the next time level and the number of proteins at this new time.
Write down all steps and details of your solution. 
\end{workout}
\vsp 

\begin{workout}
Estimate $\E_f[h(X_1,X_2)]$ using the Metropolis-Hastings algorithm where $g(x,y) = xy$ and $f$ is given in \eqref{fx1x2_def}. Use different values $N = 10^3,10^4,10^5$ and $10^6$.
\end{workout}
\vsp

\begin{workout}
In thins exercise, we extend the the scenario in Example \ref{ex:mcmc-patient1} to
a two-parameter case where we estimate the recovery rates for two different patient groups, \( \theta = (\theta^{(1)}, \theta^{(2)}) \). Each parameter, \( \theta^{(1)} \) and \( \theta^{(2)} \), stands for the recovery rate for a different group of patients and both are assumed to follow exponential distributions with independent Gamma priors.
We have two sets of observed recovery times
\begin{align*}
&\text{data}_1 = \{5, 8, 12, 7, 9, 10, 3, 6, 8, 11\} \\
&\text{data}_2 = \{10, 14, 7, 11, 13, 8, 15, 9, 10, 16\}.
\end{align*}
Each group has 10 observations, and we assume that recovery times for each group are independent and exponentially distributed:
   \[
   p(x | \theta^{(1)}) =\theta^{(1)} e^{-\theta^{(1)} x} \quad \text{and} \quad 
   p(x | \theta^{(2)}) = \theta^{(2)} e^{-\theta^{(2)} x}.
   \]
Both recovery rates, \( \theta^{(1)} \) and \( \theta^{(1)} \), are independent and follow Gamma prior
   \[
   \theta^{(i)} \sim \Gam(\alpha_i,\beta_i), \quad i=1,2.
   \]
Assume that the prior parameters are \( \alpha_1 = 2 \), \( \beta_1 = 1 \), \( \alpha_2 = 2 \), and \( \beta_2 = 1 \) which indicates that we expect similar rates for both groups.

For independent recovery times in each group, the likelihood function for each set of observations is
   \[
   p(\text{data}_i | \theta^{(i)}) = \prod_{j=1}^{10} \theta^{(i)} e^{-\theta^{(i)} x_j}, \quad 
   i = 1,2,\vsp 
   \]
and the joint likelihood of observing all data given \( \theta = (\theta^{(1)}, \theta^{(2)}) \) is
\[
p(\text{data} | \theta^{(1)}, \theta^{(2)}) = p(\text{data}_1 | \theta^{(1)})\cdot p(\text{data}_2 | \theta^{(2)}).
\]
Start with an initial values \( \theta_0= (\theta_0^{(1)}, \theta_0^{(2)}) = (1,1) \)  and use the random walk sampler with covariance matrix 
$$
\Sigma = \begin{bmatrix}
0.2 & 0 \\ 0& 0.2
\end{bmatrix}
$$
at each time step $t$ to generate the new parameter vector $\theta_{t+1}$ using the Metropolis-Hastings algorithm. Plot the histogram of generated samples, and compute the mean, variance and the confidence intervals for each parameter. 

\end{workout}
\vsp

%% file: Appendix.tex
\section{Appendix}\label{sect:appendixA}
In this appendix section we summarize some basic definitions and results from probability theory. For more details, see standard textbooks in the subject. As examples I refer you to \cite{DeGroot-Schervish:2007} and \cite{Rubinnstein-Kroese:2017}.
\vsp
\subsection{Random experiments}
An experiment whose outcome can not be determined in advance is called a {\em random experiment}.
The {\em sample space} of the random experiment is the set of all its possible outcomes. We denote the sample space by
$\Omega$.
For example, assume that a fair coin
is flipped three times. If $H$ and $T$ stand for `heads' and `tails', the sample space of this experiment is
$$
\Omega = \{HHH,HHT,HTH,THH,TTH,THT,HTT,TTT\}
$$
which contains eight possible outcomes. Here $THT$ means that the first flip lands tails, the second heads, and the third tails. Subspaces of the sample space are called {\em events}. For example the event $A$ that the second flip is heads is
$$
A = \{HHH,HHT,THH,THT\}.
$$
We say that event $A$ {\em occurs} if the outcome of the experiment is one of the elements of $A$.
Since events are sets, we can apply the usual set operations to them. For
example, the event
$$
A\cup B
$$
is the event that $A$ or $B$ or
both occur, and the event
$$
A\cap B
$$
is the event
that $A$ and $B$ both occur. Similar notation holds for unions and intersections of
more than two events.
For intersection of $n$ events $A_1,A_2,\ldots,A_n$, we usually use the abbreviation
$A_1A_2\ldots A_n = A_1\cap A_2\cap \cdots \cap A_n$, for simplicity.
The event
$
A^c
$
called the complement of $A$, is the event that
A does not occur.
Two events $A$ and $B$ are called {\em disjoint} if
their intersection is empty.
\vsp 

\begin{definition}\label{def:probability}
The probability $\pr$ is a rule that assigns a number $\pr(A)$ to each event $A\subseteq \Omega$
such that
\begin{itemize}
\item[1.] $0\leqslant \pr(A)\leqslant 1$,
\item[2.] $\pr(\Omega)=1$,
\item[3.] for any sequence $A_1,A_2,\ldots$ of disjoint events we have
\begin{equation}\label{sum_rule}
  \pr\left(\bigcup_{k}A_k  \right) = \sum_{k}\pr(A_k).
\end{equation}
\end{itemize}
\end{definition}
\vsp 

The item 1 states that the probability that the outcome of the experiment lies
within $A$ is some number between $0$ and $1$. The item 2 states that with probability
$1$ any outcome is a member of the sample space $\Omega$, and the item 3 states that for any
set of mutually disjoint events, the probability that at least one of these events
occurs is equal to the sum of their respective probabilities.

Since $A$ and $A^c$ are always mutually disjoint, and since $A\cup A^c=\Omega$,
we have from items 2 and 3 that
$$
1 = \pr(\Omega) = \pr(A \cup A^c) = \pr(A) + \pr(A^c)
$$
or equivalently
$$
\pr(A^c) = 1 - \pr(A).
$$
In other words, the probability that an event does not occur is $1$ minus the probability
that it does.

In the coin flipping experiment, since the coin is fair, the eight possible outcomes are equally likely to occur and
thus has probability $1/8$. For example
$\pr(\{HTH\})=1/8$.
Since each event $A$ is the union of the events $\{HHH\},\ldots,\{TTT\}$, and this events are disjoint, we have
$$
\pr(A)= \frac{|A|}{|\Omega|}
$$
where $|A|$ denotes the number of outcomes in $A$. For example, the probability of the event $A$ that the second flip is heads is
$\pr(A)=4/8=1/2$.
\vsp 
\subsection{Conditional probability and independence}
Assume that $B\subset \Omega$ is an event. Given that the outcome lies in $B$, the event $A$ will occur if and only if $A\cap B$ occurs and the relative chance of $A$ occurring is therefore $\pr(A \cap B)/\pr(B)$. This leads to
the definition of the {\em conditional probability} of $A$ given $B$:
\begin{equation}\label{conditional-prob}
  \pr(A|B) = \frac{\pr(A\cap B)}{\pr(B)}.
\end{equation}
For example, if a fair coin is flipped $3$ times, and $B$ is the event of total number of heads being $2$, then the probability of event $A$ that the second flip is heads given that $B$ occurs is
$\frac{2/8}{3/8}=\frac{2}{3}$ because
$$
B = \{HHT,HTH,THH\}, \quad \pr(B)=\frac{3}{8}
$$
and
$$
A = \{HHH,HHT,THH,THT\}, \quad A\cap B = \{HHT,THH\}, \quad \pr(A\cap B)=\frac{2}{8}.\vsp
$$

From \eqref{conditional-prob}, by changing the role of $A$ and $B$ we can write
\begin{equation}\label{product_rule0}
\pr(AB) = \pr(A)\pr(B|A),
\end{equation}
and this formula can be generalized for any sequence of events $A_1,A_2,\ldots, A_n$,
\begin{equation}\label{product_rule}
  \pr(A_1A_2\cdots A_n) = \pr(A_1)\pr(A_2|A_1) \pr(A_3|A_1A_2)\cdots \pr(A_n|A_1\cdots A_{n-1})
\end{equation}
which is known as the {\em product rule of probability}.

Assume that $B_1,B_2,\ldots, B_n$ are disjoin events and their union is $\Omega$. Then any event $A\subset \Omega$ can be written as
$A = \cup_{k=1}^n (A\cap B_k)$.
From the third property of
probability, i.e. \eqref{sum_rule}, we have $\pr(A)=\sum_{k=1}^{n}\pr(A\cap B_k)$. Then form \eqref{conditional-prob} we can write
\begin{equation}\label{total_probability}
  \pr(A) = \sum_{k=1}^{n}\pr(A|B_k)\pr(B_k),
\end{equation}
which is known as the {\em law of total probability}. Then, from the fact that $\pr(A)\pr(B_j|A)=\pr(AB_j) = \pr(A|B_j)\pr(B_j)$, we may write
\begin{equation}\label{bayes_rule}
  \pr(B_j|A) =\frac{\pr(A|B_j)\pr(B_j)}{\ds \sum_{k=1}^{n}\pr(A|B_k)\pr(B_k)}
\end{equation}
which is known as the {\em Bayes' rule}.

Two events $A$ and $B$ are called {\em independent} if $\pr(A|B)=\pr(A)$ which means that the occurrence of $B$ does not effect on the occurrence of $A$. An equivalent definition is: $A$ and $B$ are independent if and only if
$$
\pr(AB) = \pr(A)\pr(B).
$$
Thus definition can be extended to a sequence of events.
\begin{definition}
The events $A_1,A_2,\ldots,$ are called independent if for any $k$ and any distinct indexes $i_1,\ldots,i_k$ we have
$$
\pr(A_{i_1}A_{i_2}\cdots A_{i_k}) = \pr(A_{i_1})\pr(A_{i_2})\cdots \pr(A_{i_k}).
$$
\end{definition}
\vsp

\subsection{Random variables and distributions}
It might not always be feasible or necessary to provide a model for a random experiment through a detailed description of $\Omega$ and $\pr$. In practice, we are only interested in certain observations in the experiment.
These quantities of
interest that are determined by the results of the experiment are known as {\em random variables}, which are
usually denoted by capital letters $X$, $Y$ or $Z$ with or without subscripts.
\vsp 
\begin{example}
Consider an experiment in which a fair coin is tossed $n$ times. In this
experiment, the sample space $\Omega$ can be regarded as the set of outcomes consisting of
the $2^n$ different sequences of elementary outcomes, for instance $\{\underbrace{HTHHT\cdots T}_{n~\mathrm{times}}\}$. Assume that we are interested only in the number of heads in the observed outcome. Let $X$ be a
real-valued function defined on $\Omega$ that counts the number of heads in each outcome.
For example, for $n=10$ in the elementary outcome $s = HHTTTHTTTH$, we have $X(s)= 4$. For each possible
sequence $s$, the value $X(s)$ equals the number of
heads in the sequence. The possible values for the function $X$ are $\{0, 1, \ldots , 10\}$.
\end{example}
\vsp
\begin{definition}[Random Variable]
Let $\Omega$ be the sample space for an experiment. A real-valued function $X$
that is defined on $\Omega$ is called a random variable.
\end{definition}
\vsp

The {\em cumulative distribution function (cdf)}, or more simply the {distribution function},
$F$ of a random variable $X$ is defined for any real number $x$ by
\begin{equation*}
  F(x) = \pr(X\leqslant x).
\end{equation*}
A random variable that can take either a finite or at most a countable number of
possible values is said to be a {\em discrete random variable}. For a discrete random variable $X$ we define
its {\em probability mass function (pmf)}  $f(x)$ by
\begin{equation*}
  f(x) = \pr(X=x).
\end{equation*}
If $X$ is a discrete random variable that takes on one of the possible values $x_1, x_2, \ldots ,$
then, since $X$ must take one of these values, we have
\begin{equation*}
  \sum_{k=1}^{\infty}f(x_k)=1.
\end{equation*}

\begin{example}
Assume that a biased coin with $p$ the probability of heads is flipped $n$ times.
Suppose that we are interested only in number of heads in this experiment.
The number of heads is a random variable, let us denote it by $X$, and can take
any of the values in $\{0, 1, \ldots , n\}$.
Each elementary event $\{HTH \cdots T\}$ with exactly $k$
heads and $n-k$ tails has probability $p^k(1-p)^{n-k}$, and there are
${n\choose k}$ such events. Thus we have
\begin{equation*}
f(k) = \pr(X=k) = {n\choose k}p^k(1-p)^{n-k}, \quad k=0,1,\ldots,n.
\end{equation*}
This is the pmf of random variable $X$. The cdf of $X$ then is
\begin{equation*}
F(k) = \sum_{j=0}^{k}\pr(X=j) = \sum_{j=0}^{k}{n\choose j}p^j(1-p)^{n-j}, \quad k=0,1,\ldots,n.
\end{equation*}

\end{example}
\vsp

A random variable $X$ is said to have a continuous distribution if there exists a
positive function $f$ with total integral $1$, such that for all $a$ and  $b$,
\begin{equation*}
  \pr(a\leqslant X\leqslant b) = \int_{a}^{b} f(u)\td u.
\end{equation*}
The function $f$ is called the {\em probability density function (pdf)} of $X$. Note that in
the continuous case the cdf is given by
\begin{equation*}
  F(x) = \pr(X\leqslant x) = \int_{-\infty}^{x} f(u)\td u.
\end{equation*}
Differentiating both sides yields
$$
\frac{d}{\td x}F(x) = f(x). \vsp
$$
Note that in the discrete case we use the term probability mass function (pmf) for $f$ and in the continuous case the term
probability density function (pdf). In a more advance probability theory, both pmf and pdf can be viewed as
particular instances of a general notion called {\em probability density}. Therefore, from here on
we will call $f$ a pdf in both discrete and continuous cases.

We use the notation $X\sim f$ and $X\sim F$ to denote that $X$ has the pdf $f$ and cdf $F$, respectively. Sometimes we write $f_X$ to stress that $f$ is the distribution of random variable $X$. In Tables \ref{tb:discont-distributions} and \ref{tb:cont-distributions} the list of some well-known discrete and continuous distributions are given.

\begin{table}[ht!]
\centering
\caption{Some well-known discrete distributions}\label{tb:discont-distributions}
\begin{tabular}{llcll}
  \hline
  Name & Notation & $f(x)$ & domain of $x$ & Parameters\\
  \hline
  Bernoulli  & $\Ber(p)$  & $p^x(1-p)^{1-x}$ & $\{0,1\}$ & $0\leqslant p\leqslant 1$ \\
  Binomial & $\Bin(n,p)$ & ${n\choose x}p^x(1-p)^{n-x}$ & $\{0,1,\ldots,n\}$& $0\leqslant p\leqslant 1$ \\
  Discrete uniform & $\mathcal{DU}\{1,\ldots,n\}$ & $\frac{1}{n}$& $\{1,\ldots,n\}$& $n\in\{1,\ldots,n\}$ \\
  Geometric & $\Geo(p)$ & $p(1-p)^{1-x}$ & $\{1,2,\ldots\}$ & $0\leqslant p\leqslant 1$ \\
  Poisson & $\Poi(\lambda)$ & $\ee^{-\lambda}\frac{\lambda^x}{x!}$& $\mathbb N$& $\lambda >0$\\
  \hline
\end{tabular}

\end{table}

\begin{table}[ht!]
\centering
\caption{Some well-known continuous distributions}\label{tb:cont-distributions}
\begin{tabular}{llcll}
  \hline
  Name & Notation & $f(x)$ & domain of $x$ & Parameters\\
  \hline
  Uniform  & $\Uni(a,b)$            & $\frac{1}{b-a}$ & $[a,b]$ & $a<b$ \\
  Normal   & $\Nor(\mu,\sigma^2)$   & $\frac{1}{\sigma\sqrt{2\pi}}\ee^{-\frac{1}{2}\left(\frac{x-\mu}{\sigma}\right)^2} $& $\R$ & $\sigma>0,\, \mu\in\R$ \\
  Gamma    & $\Gam(\alpha,\beta)$ & $\frac{\beta^\alpha}{\Gamma(\alpha)} x^{\alpha-1}\ee^{-\beta x}$ & $\R_{+}$ & $\alpha,\beta>0$ \\
  Exponential & $\Exp(\lambda)$ & $\lambda \ee^{-\lambda x}$ & $\R_{+}$ & $\lambda >0$ \\
  Inv-Gamma    & $\mathcal Inv\Gam(\alpha,\beta)$ & $\frac{\beta^\alpha}{\Gamma(\alpha)} (1/x)^{\alpha+1}\ee^{-\beta/ x}$ & $\R_{+}$ & $\alpha,\beta>0$ \\
  Beta& $\Bet(\alpha,\beta)$ & $\frac{\Gamma(\alpha+\beta)}{\Gamma(\alpha)\Gamma(\beta)}x^{\alpha-1}(1-x)^{\beta-1}$ & $[0,1]$ & $\alpha,\beta>0$ \\
  Weibull& $\mathcal Weib(\alpha,\lambda)$ & $\frac{\alpha}{\lambda}(\frac{x}{\lambda})^{\alpha-1}\ee^{-(x/\lambda)^\alpha}$ & $\R_{+}$ & $\alpha,\lambda>0$ \\
  \hline
\end{tabular}

\end{table}

\subsection{Expectation and variance}
The distribution of a random variable $X$ contains all of the probabilistic information
about $X$. Summaries of the distribution, such as expected value and variance, can be useful for giving people some information about $X$ without trying to describe the entire distribution.

The intuitive idea of the {\em expectation}  or {\em mean} of a random variable is that it is the weighted
average of the possible values of the random variable with the weights equal to the
probabilities.
\vsp
\begin{definition}\label{def:expectation}
Let $X$ be a random variable with pdf $f$. The
expectation (or expected value or mean) of $X$, denoted by $\E[X]$ (or sometimes $\mu$),
is defined by
$$
\E[X]=\begin{cases}
\ds \sum_x xf(x), & \mbox{discrete case}\\
\ds \int_{-\infty}^{\infty}xf(x)\td x,&  \mbox{continuous case},
\end{cases}
$$
provided that the sum and integral in the definition are finite.
\end{definition}
\vsp

\begin{example}Let $X$ be a random variable with Bernoulli distribution with parameter $p$,
that is, assume that $X$ takes only the two values $0$ and $1$ with $\pr(X = 1) = p$. Then the
mean of $X$ is
$$
\E[X] = 0\times(1-p)+1\times p = p.
$$
If $X$ has exponential distribution with parameter $\lambda$ then
$$
\E[X] = \int_0^\infty \lambda x e^{-\lambda x}dx  = \frac{1}{\lambda}.
$$
\end{example}
\vsp
It is important to note that although
$\E[X]$ is called the expectation of $X$, it depends only on the distribution of $X$. Every
two random variables that have the same distribution will have the same expectation. So, we refer
to the expectation of a distribution even if we do not have in mind a random variable
with that distribution.

If $X$ is a random variable then a function of $X$ such as $X^2$ or $\sin(X)$ is another random variable. The expectation of a function of $X$, say $g(X)$, is simply defined as
\begin{equation}\label{def:expectation_g}
\E[g(X)]=\E_f[g(X)]=\begin{cases}
\ds \sum_x g(x)f(x), & \mbox{discrete case}\\
\ds \int_{-\infty}^{\infty}g(x)f(x)\td x,&  \mbox{continuous case},
\end{cases}
\end{equation}
provided that the sum and integral are finite. Here $f$ is the pdf of $X$.

Another useful quantity is the {variance} which measures the spread or dispersion of the distribution.
\vsp
\begin{definition}
The {\bf variance} of a random variable $X$ is denoted by $\Var(X)$ (or $\sigma^2$) and is defined by
$$
\Var(X) = \E[(X-\mu)^2] = \E[X^2]-(\E[X])^2.
$$
The square root of the variance is called {\em standard deviation}.
\end{definition}
\vsp
The variance tells us how much $X$ deviates from its mean value. 
If $X$ has infinite mean or if the mean of $X$ does not exist, we say that $\Var(X)$ does
not exist.
In Table \ref{tb:expectvar} the expectations and variances of some well-known distributions are summarized.
\vsp

\begin{table}[ht!]
\centering
\caption{Expectations and variances of some well-known distributions}\label{tb:expectvar}
\begin{tabular}{llcc}
  \hline
  Name       & Notation        & $\E(X)$              & $\Var(X)$ \\
  \hline
  Bernoulli  & $\Ber(p)$       & $p$                  &  $p(1-p)$ \\
  Binomial   & $\Bin(n,p)$     & $np$                 &  $np(1-p)$\\
  Geometric  & $\Geo(p)$       & ${1}/{p}$       &  ${(1-p)}/{p^2}$ \\
  Poisson    & $\Poi(\lambda)$ & $\lambda$            &  $\lambda$\\

  Uniform  & $\Uni(a,b)$ & ${(a+b)}/{2}$ & ${(b-a)^2}/{12}$ \\
  Normal   & $\Nor(\mu,\sigma^2)$ & $\mu$  & $\sigma^2$ \\
  Gamma    & $\Gam(\alpha,\beta)$ & ${\alpha}/{\beta}$   & ${\alpha}/{\beta^2}$\\
  Exponential & $\Exp(\lambda)$ & ${1}/{\lambda}$ & ${1}/{\lambda^2}$ \\
  Inv-Gamma    & $\mathcal Inv\Gam(\alpha,\beta)$ & $\frac{\beta}{\alpha-1},\, \alpha>1$   & $\frac{\beta^2}{(\alpha-1)^2(\alpha-2)},\, \alpha>2$\\ 
  Beta        & $\Bet(\alpha,\beta)$ &  $\frac{\alpha}{\alpha+\beta}$   & $\frac{\alpha\beta}{(\alpha+\beta)^2(1+\alpha+\beta)}$ \\
   Weibull& $\mathcal Weib(\alpha,\lambda)$ & $\lambda\Gamma(1+1/\alpha)$ & $\lambda^2\left[\Gamma(1+2/\alpha)-\Gamma(1+1/\alpha)^2 \right]$  \\
  \hline
\end{tabular}

\end{table}
\vsp

\subsection{Joint distribution}\label{sect:jointdist}
Let $X_1,X_2,\ldots,X_d$ be random variables describing some random experiments. We can collect them to a random vector
$X = (X_1,X_2,\ldots,X_d)$. Then the joint distribution of $X_1,\ldots,X_d$ is specified by the joint cdf
$$
F(x_1,\ldots,x_d) = \pr(X_1\leqslant x_1,\ldots, X_d\leqslant x_d)
$$
and the joint pdf is given in the discrete case by
$$
f(x_1,\ldots,x_d) = \pr(X_1=x_1,\ldots,X_d=x_d),
$$
and in the continuous case, $f$ is such that
$$
\pr(X\in B) = \int_B f(x_1,\ldots,x_d)\td x_1\ldots \td x_d
$$
for any measurable set $B$ in $\R^d$. The marginal pdfs can be recovered from the joint pdf by integration or summation. For example in the case of continuous random vector $(X,Y)$ with joint pdf $f(x,y)$, the pdf $f_{X}$ of $X$ is obtained as
$$
f_{X}(x) = \int f(x,y) \td y.
$$
Suppose that $X$ and $Y$ are both discrete or both continuous, with joint pdf $f$,
and suppose that $f_X(x) > 0$. Then the conditional pdf of $Y$ given $X = x$ is given
by
$$
f_{Y|X}(y|x) =\frac{f(x, y)}{f_X(x)},\quad \mbox{for all}\; y.
$$
The corresponding conditional expectation is (in the continuous case)
$$
\E[Y |X = x] =\int y f_{Y |X}(y | x) \td y .
$$
Note that $\E[Y|X=x]$ is a function of $x$, say $h(x)$. The corresponding random
variable $h(X)$ is written as $\E[Y |X]$.

When the conditional distribution of $Y$ given $X$ is identical to that of $Y$, $X$ and
$Y$ are said to be independent. More precisely:
\begin{definition}[Independent Random Variables] The random variables
$X_1,\ldots, X_d$ are called independent if we have
$$
\pr(X_1\in A_1,\ldots,X_d\in A_d) = \pr(X_1\in A_1)\times \cdots \times \pr(X_d\in A_d).
$$
for all events $\{X_i\in A_i\}$ with $A_i\subset \R$,
$i = 1,\ldots, d$.
\end{definition}
A direct consequence of the definition above for independence is that the random
variables $X_1,\ldots, X_d$ with a joint pdf $f$ (discrete or continuous) are independent if
and only if
\begin{equation}\label{independence-def}
f(x_1,\ldots , x_d) = f_{X_1} (x_1)\cdots f_{X_d}(x_d)
\end{equation}
for all $x_1,\cdots, x_d$, where $f_{X_k}$ are the marginal pdfs.

An infinite sequence $X_1, X_2,\ldots$ of random variables is called independent
if for any finite choice of parameters $i_1, i_2,\ldots,i_d$ (none of them the same)
the random variables $X_{i_1},\ldots, X_{i_d}$ are independent.
Random variables $X_1,X_2,\ldots$ that are {\em independent and identically distributed (iid)}
are frequently appeared in probabilistic models.

The expectation of  any real-valued function $g$ of variables $X_1,\ldots, X_d$ is defined as
$$
\E[g(X_1,\ldots,X_d)] = \int \cdots \int g(x_1,\ldots,x_d)f(x_1,\ldots,x_d) \td x_1\cdots \td x_d.
$$
A direct consequent of the definitions of expectation and independence (equation \eqref{independence-def}) is
\begin{equation*}
  \E[a+b_1X_1+b_2X_2+\cdots+b_dX_d] = a+b_1\E[X_1]+b_2\E[X_2] + \cdots + b_d\E[X_d]
\end{equation*}
for any sequence of independent random variables $X_1,\ldots,X_d$ and constant $a$, $b_1$, $\ldots$, $b_d$.
For independent random variables we also have (prove!)
$$
\E[X_1X_2\cdots X_d] = \E[X_1]\E[X_2] \cdots\E[X_d].
$$

Sometimes it is useful to have a summary of how much the two random variables depend on each other.
The {\em covariance} and {\em correlation} are additional notions to measure that dependence, but they
only capture a linear dependence.
\vsp
\begin{definition}
The covariance of two random variables $X$ and $Y$ with expectations $\mu_X=\E[X]$ and $\mu_Y=\E[Y]$ is defined as
$$
\Cov(X,Y) = \E[(X-\mu_X)(Y-\mu_Y)].
$$
A scaled version of the covariance is given by correlation coefficient
$$
\rho(X,Y) = \frac{\Cov(X,Y)}{\sigma_X\sigma_Y}
$$
where $\sigma_X^2=\Var(X)$ and $\sigma_Y^2=\Var(Y)$.
\end{definition}
\vsp

It can be shown that the correlation coefficient is always lies between $-1$ and $1$. Values close to $0$ show a more dependency while values close to $1$ and $-1$ stand for a linear independency between two variables.
Some important properties of the variance and covariance are listed bellow.
\begin{enumerate}
  \item $\Var(X)=\E[X^2]-(\E[X])^2$,
  \item $\Var(aX+b)=a^2\Var(X)$,
  \item $\Cov(X,Y)=\E[XY]-\E(X)\E(Y)$,
  \item $\Cov(X,Y)=\Cov(Y,X)$,
  \item $\Cov(aX+bY,Z)=a\Cov(X,Z)+b\Cov(Y,Z)$,
  \item $\Var(X+Y)=\Var(X)+\Var(Y)+2\Cov(X,Y)$,
  \item $X$ and $Y$ independent $\Rightarrow$ $\Cov(X,Y)=0$.
\end{enumerate}

For random vectors, such as $X = (X_1,\ldots, X_d)^T$, it is convenient to write the
expectations and covariances in vector notation as
\begin{align*}
 & \E(X)    = (\E[X_1],\ldots,\E[X_d])^T = (\mu_1,\ldots,\mu_d)^T =\mu\\
 & \Sigma  = \E[(X-\mu)(X-\mu)^T] = [\Cov(X_i,X_j)]_{i,j=1,\ldots,d}.
\end{align*}
Note that any covariance matrix $\Sigma$ is a symmetric positive semidefinite matrix. 
\vsp 
\subsection{Functions of random variables}
Suppose that $X$ is a random variable and $Z=g(X)$ for some monotonically increasing function $g$. To find the pdf of $Z$
from that of $X$ we first write
$$
F_Z(z) = \pr(Z\leqslant z) = \pr(X\leqslant g^{-1}(z)) = F_X(g^{-1}(z)).
$$
Differentiation with respect to $z$ then gives
$$
f_Z(z) = f_X(g^{-1}(z))\frac{d}{\td z}g^{-1}(z) = \frac{f_X(g^{-1}(z))}{g'(g^{-1}(z))}=\frac{f_X(x)}{g'(x)}.
$$
For monotonically decreasing functions, $\frac{d}{\td z}g^{-1}(z)$ in the first equation need to be replaced with its negative value.
For example let $Z = aX+b$ where $a\neq 0$ and suppose that $a>0$. Since $g^{-1}(z)=\frac{1}{a}(z-b)$ and $g'(z)=a$, we have
$f_Z(z)=\frac{1}{a}f_X(\frac{1}{a}(z-b))$. Similarly, for $a<0$ we have  $f_Z(z)=\frac{1}{-a}f_X(\frac{1}{a}(z-b))$. Thus in general
$$
f_Z(z)=\frac{1}{|a|}f_X\left(\frac{z-b}{a}\right).\vsp
$$

Now consider a random vector $X = (X_1, \ldots, X_d)^T\in\R^d$, and let $Z = AX$ where $A$ is an $m\times d$ matrix.
Then $Z$ is a random vector in $\R^m$. If we know the joint distribution
of $X$, then we can derive the joint distribution of $Z$. First, we can show that
$$
\mu_{Z} = A\mu_{X}, \quad  \Sigma_{Z} = A\Sigma_{X} A^T
$$
because
$\mu_{Z} = \E[Z]=\E[AX]=A\E[X]=A\mu_{X}$ and
$\Sigma_{Z}=\E[(Z-\mu_{Z})(Z-\mu_{Z})^T]=\E[A(X-\mu_{X})(A(X-\mu_{X}))^T]=A\E[(X-\mu_{X})(X-\mu_{X})^T]A^T
= A\Sigma_{X}A^T$. Then by using some tools from calculus, we can prove that
\begin{equation}\label{pdf_multi0}
f_{Z}(z) = \frac{1}{|\mathrm{det}(A)|} f_{X}(x).
\end{equation}
In a more general case for a fixed $x\in\R^m$, let $z = g(x)$ where $g:\R^d\to\R^m$ is invertible. Then
$$
f_{Z}(z) = \frac{1}{|\mathrm{det}J_{g}(g(x))|} f_{X}(x),
$$
where $J_g$ is the Jacobian matrix of $g$.
\vsp 
\subsection{Joint normal random variables}\label{sect:multi-normal}
Assume that $X$ has standard normal distribution, i.e., $X\sim  \Nor(0, 1)$. Then $X$ has density $f_X$ given by
$$
f_X(x) = \frac{1}{\sqrt{2\pi}} \exp\left({-\frac{x^2}{2}}\right), \quad x\in \R. \vsp
$$
If $Z=\mu + \sigma X$ then $\E[Z]=\mu$ and $\Var(Z)=\sigma^2$, i.e., $Z\sim\Nor(\mu,\sigma^2)$, thus $Z$ has density function
$$
f_Z(z) = \frac{1}{\sigma\sqrt{2\pi}}\exp\left(-\frac{1}{2}\left(\frac{z-\mu}{\sigma}\right)^2\right),\quad z\in\R.
$$
We can also state this as follows: if $Z\sim\Nor(\mu,\sigma^2)$ then
$$
\frac{Z-\mu}{\sigma}\sim \Nor(0,1),
$$
which is called {\em standardization}.
This can be generalized to $d$ dimensions. Let $X_1,\ldots,X_d$ be independent and
standard normal random variables. The joint pdf of $X$ is given
by
\begin{equation*}
f_{X}(x) = (2\pi)^{-n/2}\exp\left({-\frac{1}{2}x^Tx}\right), \quad x\in \R^d.
\end{equation*}
If $Z = \mu + BX$ for some $m\times d$ matrix $B$ then $Z$ has expectation vector $\mu$
and covariance matrix $\Sigma=BB^T$. The random vector $Z$ is said to have a {\em jointly normal} or {\em multivariate normal} distribution. We write
$$
Z \sim \Nor(\mu, \Sigma).
$$
In a special case when $B$ is $d\times d$ and invertible, the pdf of $Z$ can be derived as follows.
If $Y = Z-\mu$ then from \eqref{pdf_multi0} the pdf of $Y$ is given by
$$
f_{Y}(y) = \frac{(2\pi)^{-n/2}}{|\mathrm{det(B)}|}\exp\left({-\frac{1}{2}(B^{-1}y)^T(B^{-1}y)}\right).
$$
Using the facts that $|\mathrm{det(B)}|=\sqrt{|\mathrm{det(\Sigma)}|}$ and $(BB^T)^{-1}=\Sigma^{-1}$, we have
$$
f_{Y}(y) = \frac{1}{\sqrt{(2\pi)^d|\mathrm{det(\Sigma)}|}}\exp\left({-\frac{1}{2}y^T\Sigma^{-1}y}\right).
$$
Since $Z=Y+\mu$ for a constant vector $\mu$, we have
$f_{Z}(z) = f_{Y}(z-\mu)$, and therefore
\begin{equation*}
  f_{Z}(z) = \frac{1}{\sqrt{(2\pi)^d|\mathrm{det(\Sigma)}|}}\exp\left({-\frac{1}{2}(z-\mu)^T\Sigma^{-1}(z-\mu)}\right).
\end{equation*}
Conversely, given a covariance matrix $\Sigma=(\sigma_{ij})$, there exists a unique lower
triangular matrix
$$
B =\begin{bmatrix}
     b_{11} & 0 & \cdots & 0 \\
     b_{21} & b_{22} & \cdots & 0 \\
     \vdots & \vdots & \ddots & \vdots \\
     b_{n1} & b_{n2} & \cdots & b_{nn}
   \end{bmatrix}
$$
such that $\Sigma = BB^T$. Recall that every covariance matrix is positive semi-definite. The lower triangular matrix $B$ can be obtained efficiently via the {\em Cholesky factorization.}
\vsp

\subsection{Generating normal random variables}\label{sect:appendix_normal_gen}
In this section we use the {\em inverse transform method} to generate normal random variables. Before start reading this part, read Section \ref{sect:generation} if you are not familiar with the inverse transform method. 
The problem is that the inverse of the normal cdf is not available explicitly, and it is numerical expensive to compute.
A suitable transformation helps to make things simpler.
%\begin{example}[Sampling from normal distribution]\label{ex:normal_gen}
Assume that $X\sim \Nor(\mu,\sigma^2)$. The pdf of $X$ is
\begin{equation}\label{pdf_normal}
f(x)=\frac{1}{\sigma\sqrt{2\pi}}\exp\left(-\frac{1}{2}\left(\frac{x-\mu}{\sigma}\right)^2\right),\quad x\in\R.
\vsp
\end{equation}
It is enough to generate from the standard normal distribution $\Nor(0,1)$ because any random variable $Z\sim \Nor(\mu,\sigma^2)$ can be written as $Z = \sigma X+\mu$ where $X\sim\Nor(0,1)$. The cdf of the standard normal distribution then is
$$
F(x) = \frac{1}{\sqrt{2\pi}}\int_{-\infty}^{x} \exp\left(-\frac{1}{2}u^2\right)\td u = \frac{1}{2}\left[1+\mathrm{erf}\left(\frac{x}{\sqrt 2}\right)\right]
$$
where $\mathrm{erf}$ is the error function. The error function should be computed via its series representation or numerical integration.
Computing an accurate inverse for $F$ is numerically expensive thus the inverse transform method is not practically applicable.

One of the earliest method for generating from the standard normal distribution was developed by Box and Muller as follows:
Assume that $X\sim \Nor(0,1)$ and $Y\sim \Nor (0,1)$ are two independent random variables. The variable $(X,Y)$ has joint distribution
$$
f_{X,Y}(x,y) = \frac{1}{2\pi}\exp\left(-\frac{1}{2}(x^2+y^2)\right), \quad (x,y)\in \R^2.
$$
By transferring the variables to polar coordinates $(r,\theta)$ with change of variables
\begin{equation}\label{polar_change}
  x = r\cos\theta, \quad y = r\sin\theta,\quad r\geqslant 0,\; \theta\in[0,2\pi),
\end{equation}
the joint distribution of the transferred variables $(R,\Theta)$ becomes
$$
f_{R,\Theta}(r,\theta) = \frac{r}{2\pi}\exp\left(-\frac{r^2}{2}\right), \quad r\geqslant0,\; \theta\in[0,2\pi),
$$
where the factor $r$ behind the exponential term comes from the determinant of the Jacobian matrix of the transformation.
Now, we have
$$
F_{R,\Theta}(r,\theta) = \int_{0}^{\theta}\int_{0}^{r} \frac{r'}{2\pi}\exp\left(-\frac{r'^2}{2}\right)dr'd\theta' = -\frac{\theta}{2\pi}\exp\left(-\frac{r^2}{2}\right),
$$
which shows that $R$ and $\Theta$ are independent with cdf's
$$
F_\Theta(\theta)=-\frac{\theta}{2\pi},\;\theta\in[0,2\pi),  \quad F_R(r) = \exp\left(-\frac{r^2}{2}\right),\; r\geqslant0.
$$
Consequently, $\Theta \sim \Uni (0,2\pi)$ and $R$ has the same distribution as $\sqrt Q$ with $Q\sim {\Exp}(1/2)$ because
$$
F_{\sqrt Q}(r) = \pr(\sqrt Q\leqslant r) =\pr(Q\leqslant r^2) = \exp(-r^2/2).
$$
Both $R$ and $\Theta$ are easy to generate. First we generate $U_1\sim \Uni(0,1)$ and $U_2\sim \Uni(0,1)$. Then we set $\Theta = 2\pi U_1$ and $R = \sqrt{-2\ln U_2}$ from
\eqref{exp-gen-for}. Finally, we transfer $R$ and $\Theta$ back to $X$ and $Y$ via  transformation \eqref{polar_change}:

\begin{align*}
  X & = R\cos(\Theta) = \sqrt{-2\ln U_2}\cos(2\pi U_1), \\
  Y & = R\sin(\Theta) = \sqrt{-2\ln U_2}\sin(2\pi U_1).
\end{align*}
Note that, using this approach we generate two standard random variables $X$ and $Y$ from two uniform variables $U_1$ and $U_2$.
A Python code is given here.

\begin{shaded}
\vspace*{-0.3cm}
\begin{verbatim}
def RandNormal(mu, sigma2, N):
    U1 =  np.random.rand(N)
    U2 =  np.random.rand(N)
    X = np.sqrt(-2*np.log(U1))*np.cos(2*np.pi*U2)
    X = mu + np.sqrt(sigma2)*X
    return X
\end{verbatim}
\vspace*{-0.3cm}
\end{shaded}

The code snippet below generates $N$ random normal points and plots the corresponding histograms.
Outputs are shown in Figure \ref{fig:normal_gen} for $N = 500$ and $5000$.

\begin{shaded}
\vspace*{-0.3cm}
\begin{verbatim}
import numpy as np
import matplotlib.pyplot as plt
plt.figure(figsize = (5,3))
N = 500
X = RandNormal(0, 1, N)
plt.hist(X, bins = 30, histtype = 'bar', color = 'red', density = 'true')
x = np.linspace(-4,4,200)
f = 1/(np.sqrt(2*np.pi))*np.exp(-x**2/2)
plt.plot(x,f,linestyle = '-', color = 'blue')
plt.title('Histogram of $X$ and the pdf $f(x)$')
plt.xlabel('$X$'); plt.ylabel('Frequency %')
\end{verbatim}
\vspace*{-0.3cm}
\end{shaded}

 \begin{center}
 \includegraphics[scale=0.68]{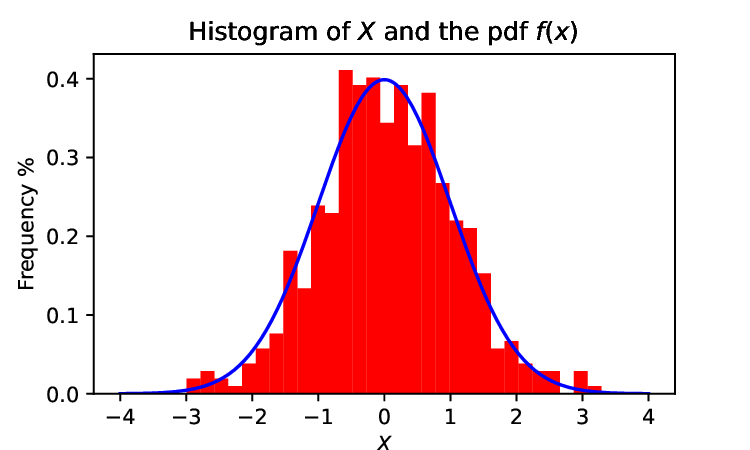}\includegraphics[scale=0.68]{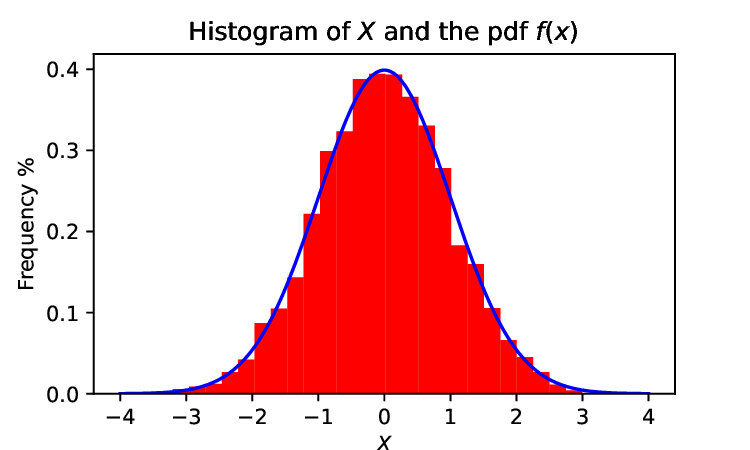}
 \captionof{figure}{Histograms of generated random points from the standard normal distribution with $N=500$ (left), $N=5000$ (right)}\label{fig:normal_gen}
 \end{center}
\vsp 

\subsection{Generating from multivariate distributions}\label{sect:randvectorgen}
Before start reading this section, you may need to take a look at Section \ref{sect:generation} for random point generation algorithms for univariate distributions. 

To generate a random vector $X = (X_1,\ldots,X_d)$ for a given $d$ dimensional distribution with pdf $f(x)=f(x_1,\ldots,x_d)$, we can use the product rule \eqref{product_rule}
\begin{equation*}
  f(x_1,\ldots,x_d) = f_1(x_1) f_2(x_2|x_1)\cdots f_d(x_d|x_1,\cdots,x_{d-1})
\end{equation*}
where $f_1(x_1)$ is the marginal pdf of $X_1$, $f_k(x_k|x_1,\ldots,x_{k-1})$ is the conditional pdf of $X_k$ given $X_1=x_1, \ldots, X_{k-1}=x_{k-1}$. If $X_1,X_2,\ldots, X_d$ are independent then $f_k(x_k|x_1,\ldots,x_{k-1})=f_k(x_k)$, and
the techniques of univariate distributions can be simply applied to each component
individually.
To generate $X$ in the case that variables are dependent, one can first generate $X_1$ from pdf $f_1(x_1)$, then given $X_1=x_1$ generate $X_2$ from $f_2(x_2|x_1)$, and so on. To run this approach, a knowledge on conditional distributions is required. Markov models provide
a feasible and simple way to obtain such a knowledge and to generate from a general joint distribution. See section \ref{sect:stoch_proc_gen}.

%\subsubsection*{Acceptance-rejection methods}

\subsubsection*{Drawing from the multi-normal distribution}

Assume that $X \sim \Nor (\mu,\Sigma)$, where $\mu= (\mu_1,\ldots,\mu_d)^T$ is the mean vector and $\Sigma\in \R^{d\times d}$ is the covariance matrix of $X$. Then we have $X = \mu + BZ$ where $Z$ is a vector of iid random variables with distributions $\Nor(0,1)$, and $B$ is the Cholesky factorization of the covariance matrix $\Sigma$. The following Python code can be used to generate
\verb+N+
 multi-normal random variables with expectation vector \verb+mu+ and covariance matrix \verb+Sigma+.
\begin{shaded}
\vspace*{-0.3cm}
\begin{verbatim}
def RandMultiNormal(mu, Sigma, N):
    dim = np.size(mu)
    Z = X = np.zeros([dim,N])
    if dim > 1:
       B = np.linalg.cholesky(Sigma)
    else:
       B = [np.sqrt(Sigma)]
    for d in range(dim):
        Z[d,:] = np.random.normal(0, 1, N)
    X = np.matlib.repmat(mu, N, 1).T + np.matmul(B,Z)
    return X
\end{verbatim}
\vspace*{-0.3cm}
\end{shaded}
\vsp 

\subsection{Limit theorems}\label{sect:limit-thm}
In this part we review the {\em law of large numbers}
and the {\em central limit theorem}. Both are associated with sums of independent
random variables.

Let $X_1,X_2, \ldots$ be iid random variables with expectation $\mu$ and variance $\sigma^2$. For
each $n$, let
$$
S_n = X_1+X_2+\cdots + X_n.
$$
Since $X_1,X_2, \ldots$ are iid, we have $\E[S_n] = n \E[X_1] = n\mu$ and $\Var(S_n) = n \Var(X_1) = n\sigma^2$.
The law of large numbers states that $S_n/n$ is close to $\mu$ for large $n$. Here is the
more precise statement.
\vsp 
\begin{theorem}[Strong Law of Large Numbers] If $X_1, \ldots , X_n$ are iid with
expectation $\mu$, then
$$
\pr \left(\lim_{n\to\infty} \frac{S_n}{n}=\mu\right) = 1.
$$
\end{theorem}
\vsp 
The central limit theorem describes the limiting distribution of $S_n$ (or $S_n/n$),
and it applies to both continuous and discrete random variables. Loosely, it states
that the random sum $S_n$ has a distribution that is approximately normal, when $n$
is large. The more precise statement is given next.
\vsp
\begin{theorem}[Central Limit Theorem] If $X_1,\ldots, X_d$ are iid with expectation
$\mu$ and variance $\sigma^2 <\infty$, then for all $x\in\R$,
\begin{equation*}
\lim_{n\to\infty}\pr \left( \frac{S_n-n\mu}{\sigma\sqrt n}\leqslant x\right) = \Phi(x)
\end{equation*}
where $\Phi$ is the cdf of the standard normal distribution.
\end{theorem}
\vsp
In other words, $S_n$ has a distribution that is approximately normal, with expectation
$n\mu$ and variance $n\sigma^2$. These theorem are valid independent of the type of distribution of
$X_1,X_2,\ldots$.

There is also a central limit theorem for random vectors. The multidimensional
version is as follows: Let $X_1,\ldots ,X_n$ be iid random vectors with expectation vector
$\mu$ and covariance matrix $\Sigma$. Then for large values of $n$ the random vector $X_1+\cdots + X_n$ has
approximately a multivariate normal distribution with expectation vector $n\mu$ and
covariance matrix $n\Sigma$.
\vsp